\documentclass{article}
\usepackage{graphicx} 
\usepackage[letterpaper,top=2cm,bottom=2cm,left=3cm,right=3cm,marginparwidth=1.75cm]{geometry}
\usepackage{cite, graphicx, amsmath, array, url, hyperref, xstring, amsthm, amsfonts, mathtools, comment, booktabs, xcolor}
\usepackage{amssymb}
\usepackage{subcaption}
\usepackage{algpseudocode}
\usepackage{algorithm}
\usepackage{listings}
\usepackage{hyperref}
\usepackage{multirow}
\DeclareMathOperator*{\argmax}{arg\,max}
\DeclareMathOperator*{\argmin}{arg\,min}
\usepackage{natbib} 

\usepackage{csquotes}

\providecommand{\keywords}[1]

\newtheorem{theorem}{Theorem}[section]
\newtheorem{lemma}{Lemma}[section]
\newtheorem{proposition}{Proposition}[section]

\theoremstyle{EX}
\newtheorem{remark}{Remark}[section]

\newcommand{\abs}[1]{\left\lvert#1\right\rvert}

\newcommand{\Halmos}{\hfill$\square$}

\DeclareMathOperator{\dom}{dom}
\DeclareMathOperator{\intt}{int}
\DeclareMathOperator{\conv}{conv}

\title{Distributionally Robust Airport Ground Holding Problem under Wasserstein Ambiguity Sets}
\author{Haochen Wu\footnote{Department of Aerospace Engineering, University of Michigan, Ann Arbor, MI, USA \texttt{\href{mailto:haocwu@umich.edu}{haocwu@umich.edu}}}\;, Alexander S. Estes\footnote{Institute for Systems Research, Department of Decisions, Operations, and Information Technology, University of Maryland, College Park, MD, USA \texttt{\href{mailto:aestes@umd.edu}{aestes@umd.edu}}}, Max Z. Li\footnote{Department of Aerospace Engineering, Department of Civil and Environmental Engineering, Department of Industrial and Operations Engineering, University of Michigan, Ann Arbor, MI, USA \texttt{\href{mailto:maxzli@umich.edu}{maxzli@umich.edu}}}}
\date{}

\begin{document}

\maketitle
\begin{abstract}
Ground Delay Programs (GDPs) mitigate demand--capacity imbalances by holding flights on the ground when an airport's arrival capacity is reduced, thereby reducing costly airborne holding. A central challenge is that day-to-day demand-capacity balancing relies on accurate predictions of airport capacities. However, these predictions are deeply uncertain: forecast errors, operational disruptions, and climate change-driven shifts in weather severity can induce distribution shifts in capacity outcomes. Thus, policies optimized for a single predicted distribution may be brittle out of sample.

We develop a distributionally robust formulation of the single-airport ground holding problem (dr-SAGHP) using a $1$-Wasserstein ambiguity set with an $\ell_1$ transportation metric over the full physical capacity support. We exploit the
problem's structural properties to develop an exact dynamic programming separation algorithm within delayed constraint generation
(DCG-DP). We also develop K-IL, which combines Kelley's cutting-plane method with an integer L-shaped framework to solve a
a restricted-support approximation, rs-dr-SAGHP. Delayed constraint generation with MILP separation (DCG-MILP) serves as the
computational baseline. Computational experiments demonstrate substantial runtime reductions for DCG-DP and K-IL relative to this baseline.

Airport capacity distributions are constructed using Gaussian process regression based on historical weather and operational data at Denver International Airport from 2021 to 2023. We develop a chronological expanding-window calibration procedure to select the Wasserstein radius and the results show that the calibrated dr-SAGHP approaches reduce mean operational cost by approximately $2.6\%$ and the Conditional Value at Risk $\mathrm{CVaR}_{90}$ by approximately $8\%$ relative to the stochastic SAGHP benchmark, while avoiding the excessive conservatism observed in a classical robust optimization benchmark.
\end{abstract}

{
  \small	
  \textbf{\textit{Funding:}}
}
\begin{keywords}
The NASA University Leadership initiative (grant \#80NSSC24M0068) provided funds to assist the authors with their research, but this article solely reflects the opinions and conclusions of its authors.
\end{keywords}

{
  \small	
  \textbf{\textit{Keywords:}}
}

\section{Introduction}\label{sec:Intro}
Demand--capacity imbalances in the air transportation system are a major operational challenge with significant negative impacts on airlines, passengers, and the environment. There are a variety of potential bottlenecks within the air transportation system, with airport arrival and departure capacities being major components \citep{c22}. Typically, if airport capacity constraints can be identified proactively, air traffic managers prefer to delay flights before they are airborne, as airborne delay costs typically drastically outweigh delay costs incurred on the ground \citep{c1}. Within the U.S. National Airspace System (NAS), such demand--capacity balancing actions take the form of Traffic Management Initiatives (TMIs) known as Ground Delay Programs (GDPs) \citep{c23}.

However, the efficacy of these strategies is increasingly challenged by shifting environmental baselines that disrupt traditional capacity predictions. This shift is driven by the fact that climate change is not only an outcome of aviation activity, but also a driver of operational risk that aviation must adapt to \citep{ryley2020climate}. Climate change is expected to influence aviation through changes in temperature and precipitation patterns, sea-level rise, wind changes, and the impacts of more frequent and/or more intense extreme weather events \citep{ryley2020climate}. These hazards can directly affect airport and air traffic management performance via runway configuration changes, convective and storm-related constraints, and disruption from heavy precipitation and flooding, all of which translate into time-varying, uncertain airport capacities. In this sense, climate change can amplify both the frequency and the complexity of weather-driven capacity shortfalls, strengthening the need for decision models that remain reliable under imperfect predictive information.

Given known parameters such as airport capacities, nominal flight times between airports, and flight schedule information, the family of optimization models known as \emph{Ground Holding Problems} (GHPs) can be solved to obtain optimal rescheduling decisions to minimize incurred airborne and ground delay costs \citep{c1}. 

Previous works have addressed the formulation of both single airport and muti-airport variants of GHPs, referred to as the Single Airport Ground Holding Problem (SAGHP) and Multi-Airport Ground Holding Problem (MAGHP) respectively \citep[e.g.][]{richetta1993saghp,c1,c2,c3,mukherjee2007dsaghp,estes2020saghp}. In the deterministic SAGHP and MAGHP, the airport capacity (or capacities) is assumed to be known with certainty \citep[e.g.][]{c1,c2,c3}. However, in practice, airport capacity is often uncertain and can vary significantly over time. Recognizing the role that uncertainty plays in designing realistic ground holding policies, stochastic versions of the SAGHP and MAGHP (s-SAGHP and s-MAGHP, respectively) have also been examined: techniques such as stochastic programming \citep{atm1, c6, mukherjee2007dsaghp,estes2020saghp} and chance-constrained programming \citep{c7} have been used to model uncertainty in airport capacities. The results of these models demonstrate that stochastic formulations can balance robustness and cost, and often reduce airborne costs relative to deterministic policies.

Despite these advances, a key practical vulnerability remains: stochastic GHPs rely on a probabilistic model of future airport capacities, typically derived from weather forecasts, runway configuration predictions, and historical data. Forecast uncertainty, model estimation error, incomplete information available in operational decision support systems, and climate-driven nonstationarity can all lead to \emph{distribution shift}, i.e., realized capacity trajectories that deviate systematically from the assumed scenario probabilities \citep{c19,c20,c4}. When the assumed distribution is misspecified, a policy optimized for the nominal distribution can perform poorly out of sample, especially under rare but operationally consequential capacity shortfalls. This motivates an \emph{uncertainty-of-uncertainties} perspective: even the estimated capacity distribution is itself uncertain and may change over time, so robust decision support should hedge not only against capacity realizations but also against errors in the distributional inputs.

The goal of this paper is to address these uncertainty-of-uncertainties challenges by formulating and solving a stochastic GHP that is \emph{distributionally robust} to inaccuracies in the probabilistic airport capacity scenarios. Specifically, we propose a \emph{distributionally robust optimization (DRO)} approach that protects against plausible distributional perturbations around an estimated capacity distribution, thereby producing ground holding policies with improved reliability under distribution shift.

\subsection{Contributions}\label{subsec:contrib}
This paper makes the following contributions:
\begin{itemize}
    \item \textbf{A distributionally robust formulation for the single airport ground holding problem.}
    We develop a two-stage distributionally robust formulation of the single airport ground holding problem (dr-SAGHP) using a $1$-Wasserstein ambiguity set defined over the full physical support of airport capacity with the $\ell_1$ distance as its transportation metric. The proposed formulation allows the adversarial distribution to vary over all physically feasible capacity realizations within the prescribed support, thereby accounting for sampling error and potential distributional misspecification in airport capacity predictions. This provides a principled framework for hedging against uncertainty beyond the nominal predictive distribution while preserving the operational structure of the single airport ground holding problem.
     \item \textbf{An efficient exact solution algorithm for dr-SAGHP.}
        We develop a dynamic programming separation algorithm that exploits the structural properties of the single-airport ground holding problem and solves the separation problem exactly. Incorporated into delayed constraint generation, this method yields an exact solution for dr-SAGHP, referred to as DCG-DP. Computational experiments demonstrate that DCG-DP substantially reduces runtime relative to the DCG-MILP baseline.

    \item \textbf{A support-restricted approximation and an efficient solution algorithm.}
        We consider an approximation of dr-SAGHP that restricts adversarial distributions to the support of the estimated distribution, which we call rs-dr-SAGHP. To solve this problem, we develop K-IL, which combines Kelley's cutting-plane method with an integer L-shaped framework. Computational experiments show that solving rs-dr-SAGHP with K-IL provides substantial computational savings. After independent calibration of the Wasserstein radii, rs-dr-SAGHP achieves out-of-sample operational performance comparable to that of dr-SAGHP.

     \item \textbf{Computational experiments with data-driven capacity predictions.}
    We construct predictive airport capacity distributions using Gaussian Process Regression based on historical weather and operational data and develop a chronological expanding-window procedure to calibrate the Wasserstein radius and other model parameters without using test information. The proposed models are evaluated on held-out severe-weather days from real operational data and compared with stochastic and classical robust SAGHP benchmarks. The results show that dr-SAGHP improves both mean operational performance and tail-risk protection while avoiding the excessive conservatism associated with classical robust optimization. In addition, sensitivity analyses under simulated mean capacity reductions and covariance inflation demonstrate the robustness of the proposed framework under more severe distributional shifts caused by climate change.
\end{itemize}
In section~\ref{sec:lit_review}, we provide a literature review. In section~\ref{sec:Method}, we define variants of the SAGHP, including dr-SAGHP, and present a deterministic equivalent for the d-SAGHP. Our proposed solution approach is given in section~\ref{sec:solution_approach}. We then conduct numerical experiments to evaluate the out-of-sample performance of the proposed distributionally robust ground holding policies against those derived from s-SAGHP and robust optimization approach in section~\ref{sec:Experiment}.
\section{Literature Review} \label{sec:lit_review}
As discussed in the introduction, our work builds on the literature on the stochastic SAGHP. While several variants of this problem and their properties have been studied in the literature \cite[e.g.]{c5,c21,c6,c7,mukherjee2007dsaghp,estes2020saghp}, we are the first to solve the distributionally robust variant. As we will show, this provides better performance in the presence of distributional drift or misspecification. 

\subsection{Stochastic Air Traffic Flow Management}

More broadly, there exists work on other air traffic flow management (ATFM) problems under uncertainty beyond GHPs. These studies mainly include stochastic programming, chance-constrained optimization, robust optimization, and distributionally robust optimization, which differ in how uncertainty is represented and controlled.

In stochastic programming, uncertainty is represented by a probability distribution, usually approximated through a finite set of scenarios with associated probabilities, and decisions are optimized based on expected performance. For example, \citep{atm1,atm2} develop two-stage and multistage stochastic integer programs, respectively, in which tactical decisions such as ground delay, airborne holding, and flight rerouting are optimized against sampled capacity scenarios. Chance-constrained approaches also rely on a probabilistic description of uncertainty, but instead control risk by limiting the probability of capacity violations. Such formulations have been applied to ATFM under disruptive weather by bounding the probability that sector capacities are exceeded \citep{atm4}.

Robust optimization takes a different approach by assuming that uncertain parameters belong to a prescribed uncertainty set and optimizing against the worst-case realization in that set. For example, \citep{atm3} propose a robust and adaptive ATFM approach under weather-driven capacity uncertainty. This approach is pointwise robust, whereas the DRO framework considered in this paper is distributionally robust. Specifically, robust optimization minimizes the worst-case realized cost over a set of capacity realizations, while DRO minimizes the worst-case expected cost over a set of plausible probability distributions. Compared with stochastic programming, DRO therefore accounts explicitly for uncertainty in the underlying probability distribution rather than relying on a single estimated distribution.

Collectively, this literature demonstrates the value of explicitly modeling capacity uncertainty in ATFM while also highlighting the different assumptions made about the underlying uncertainty. The distributionally robust framework developed in this paper may therefore also be extended to other ATFM problems in which the capacity distribution itself is uncertain.

\subsection{Distributionally Robust Optimization}
Distributionally robust optimization (DRO) provides another way to model uncertainty when the true probability distribution is not known with confidence. Instead of optimizing with respect to a single estimated distribution, DRO considers a set of plausible distributions, referred to as an \emph{ambiguity set}, and optimizes the worst-case expected value over this set \citep{c8,c9,c10,c18}. Different ambiguity sets have been proposed in the literature, including moment-based and distance-based sets. Among them, Wasserstein ambiguity sets have received increasing attention because they can be constructed directly from data and often lead to tractable reformulations with good out-of-sample performance \citep{c9,c18}.

For two-stage problems, \citet{duque2022distributionally} study distributionally robust two-stage stochastic programs with linear recourse under Wasserstein and optimal transport ambiguity sets and develop cutting-plane methods for the resulting formulations. Similarly, \citet{gamboa2021decomposition} consider two-stage data-driven Wasserstein DRO problems with right-hand-side uncertainty and develop several decomposition methods, including column-and-constraint generation and Benders-type approaches. These studies show that Wasserstein DRO can be applied effectively to two-stage recourse problems and also highlight the importance of decomposition methods for solving large-scale formulations.


In this paper, we adopt a Wasserstein ambiguity set to develop a distributionally robust formulation of the single-airport ground holding problem. The resulting ground holding policies hedge against misspecification of airport capacity distributions and
distribution shifts arising from forecast errors and potentially changing climate conditions. We exploit the structure of SAGHP to develop a dynamic programming separation algorithm that, when incorporated into delayed constraint generation, solves dr-SAGHP exactly and more efficiently than the baseline approach using MILP separation. We also develop K-IL, an decomposition approach for solving an restricted-support
approximation of dr-SAGHP efficiently.

\subsection{Climate Change Impacts}
A growing body of attribution and process-based evidence indicates that anthropogenic climate change is already altering the frequency and intensity characteristics of several high-impact weather hazards. Event-attribution syntheses find especially robust signals for increased likelihood and/or intensity of heat extremes and heavy rainfall, with more mixed or regionally dependent signals for drought and some storm types \citep{clarke2022extreme,faranda2022climate}. For tropical cyclones, recent work links observed global increases in rapid intensification to thermodynamic environments that have become more favorable under anthropogenic warming, increasing the probability of conditions conducive to very fast storm strengthening \citep{bhatia2022potential}. In the hydrologic domain, a central implication is the breakdown of stationarity assumptions used in risk analysis and design: climate-driven non-stationarity and deep uncertainty imply that return levels/return periods of floods and extreme precipitation should be treated as time-varying quantities rather than fixed constants \cite{mondal2016hydrologic}. Empirically and in downscaled projections, multiple studies show that short-duration rainfall extremes are intensifying and that standard intensity--duration--frequency (IDF) relationships can shift materially over time, motivating periodic updating of design curves and climate-informed frequency estimation for resilient infrastructure planning \citep{maity2022changing,tamm2023intensification,doulabian2023non,nagai2026spatiotemporal}. Together, these findings support a planning paradigm in which weather extremes are both intensifying and becoming less well-described by historical distributions, increasing the value of decision methods that remain robust under distributional change and imperfect probabilistic information.

Climate change affects transportation systems through both gradual shifts in baseline conditions and changes in extreme-event regimes that drive disruptions, damages, and cascading failures across interdependent networks. Reviews of road transport infrastructure document multiple climate-sensitive failure and degradation pathways---including flooding and scour, coastal inundation and sea-level rise impacts on low-lying assets, landslides, and heat-related deterioration of pavements and structures---and emphasize that effective adaptation typically combines engineering hardening, improved drainage and design standards, operational measures, and risk-informed asset management \citep{moretti2018climate,de2022climate}. The resilience literature further highlights that transportation performance impacts are system-level (not just asset-level): weather shocks propagate via network connectivity, capacity reductions, and recovery constraints, so evaluating resilience requires moving beyond static risk notions toward robustness, redundancy, rapidity of recovery, and adaptive capacity \citep{markolf2019transportation,zhou2019resilience,jaroszweski2014impact}. 
\section{The Distributionally Robust Single Aiport Ground Holding Problem}\label{sec:Method}
Distributionally Robust Optimization is a variation of stochastic optimization techniques that can handle distributional uncertainties by defining an ambiguity set that represents a set of possible distributions for stochastic input parameters.  In this section, we present the deterministic SAGHP \citep{c3}, define stochastic and distributionally robust variants of the SAGHP, and  derive a deterministic equivalent form of the distributionally robust SAGHP under the Wasserstein ambiguity set. 
\subsection{Deterministic and Two-Stage Stochastic formulations of SAGHP}
We begin by presenting a deterministic model for the SAGHP, which we refer to as d-SAGHP. d-SAGHP is a special case of a multi-airport model proposed in \cite{c1} and referred to as the VBO model in \cite{c3}. A fundamental assumption of this model is that airport congestion arises solely from insufficient arrival capacity, while departure capacity is assumed to be infinite and does not impose a bottleneck on operations. While we focus on the single airport case, similar methods could be applied to the multi-airport case.

The objective of d-SAGHP is to minimize total ground holding delay by strategically shifting delays from the airborne phase to the ground phase. The operational environment is governed by three categories of constraints: 
\begin{enumerate}
    \item \textbf{Capacity constraints~\eqref{eq:dSAGHP_1b}:} bound the number of arrivals within a given time period based on airport capacity;
    \item \textbf{Assignment constraints~\eqref{eq:dSAGHP_1c}:} ensure each flight $f$ is assigned a unique actual arrival time $t \in \mathcal{T}_f$
    \item \textbf{Coupling constraints~\eqref{eq:dSAGHP_1d}:} preserve the temporal dependencies between preceding and successive connecting flights $(f_1, f_2) \in \mathcal{C}$ based on the buffer $S_{f_1, f_2}$
\end{enumerate}
The mathematical formulation of d-SAGHP is shown as follow:
\begin{subequations}
\label{eq:dSAGHP}
\begin{align}
\min_{x} \quad
    & \sum_{f\in \mathcal{F}} C_g
    \left(
        \sum_{t\in \mathcal{T}_f} t\,x_{f,t} - r_f
    \right)
    \label{eq:dSAGHP_objfunc} \\
\text{s.t.} \quad
    & \sum_{f\in \mathcal{F}} x_{f,t} \le K_t,
    && \forall t\in\mathcal{T},
    \label{eq:dSAGHP_1b} \\
    & \sum_{t\in \mathcal{T}_f} x_{f,t} = 1,
    && \forall f\in \mathcal{F},
    \label{eq:dSAGHP_1c} \\
    & \begin{aligned}[t]
        & \sum_{t\in \mathcal{T}_{f_1}} t\,x_{f_1,t}
          - r_{f_1} + S_{f_1,f_2} \\
        &\le
          \sum_{t\in \mathcal{T}_{f_2}} t\,x_{f_2,t} - r_{f_2}
      \end{aligned}
    && \forall (f_1,f_2)\in\mathcal{C},
    \label{eq:dSAGHP_1d} \\
    & x_{f,t}\in\{0,1\},
    && \forall f\in \mathcal{F},\ \forall t\in \mathcal{T}.
    \label{eq:dSAGHP_1e}
\end{align}
\end{subequations}
We note that the definitions of parameters are in Table~\ref{tab:notation}.
\begin{table}[H]
\centering
\caption{Sets, Parameters, and Decision Variables}
\label{tab:notation}
\begin{tabular}{ll}
\toprule
\textbf{Notation} & \textbf{Description} \\ \midrule
\textit{Sets and Indices} & \\
$\mathcal{T}$ & Set of time periods, indexed by $t$ \\
$\mathcal{F}$& Set of flights destined for the impacted airport, indexed by $f$ \\
$\mathcal{C}$ & Set of connecting flight pairs $(f_1, f_2)$ \\
$\mathcal{T}_f$ & Set of available actual arrival times for flight $f$, $\mathcal{T}_f = \{r_f, \dots, T\}$\\
$\widehat{\Xi}$ & Set of realized capacity trajectories \\
$\Xi$ & Set of physically feasible capacity trajectories\\
$\mathcal{M}(\widehat{\Xi})$ & Set of probability distributions supported on $\widehat{\Xi}$ \\
$\mathcal{M}(\Xi)$ & Set of probability distributions supported on $\Xi$ \\
$\mathcal{P}_{\varepsilon}(\widehat{\mathbb{P}})$
& Wasserstein ambiguity set with center $\widehat{\mathbb{P}}$, radius $\varepsilon$, and support $\Xi$ \\

$\mathcal{\widehat{P}}_{\varepsilon}(\widehat{\mathbb{P}})$
& Wasserstein ambiguity set with center $\widehat{\mathbb{P}}$, radius $\varepsilon$, and the support of estimated distribution $\widehat{\Xi}$ \\
 \midrule
\textit{Parameters} & \\
$\overline{K}$ & Upper bound on airport capacity  \\
$K_{t}$ & Airport capacity at time $t$ \\
$r_f$ & Scheduled arrival time of flight $f$ \\
$C_g$ & Unit ground holding cost \\
$C_h$ & Unit airborne holding cost \\
$S_{f_1, f_2}$ & Maximum delay for $f_1$ without affecting successive flight $f_2$ \\ 
$\widehat{\mathbb{P}}$ & Estimated distribution supported on $\widehat{\Xi}$\\
\midrule
\textit{Decision Variables} & \\
$x_{f,t}$ & Binary variable; 1 if flight $f$ lands at time $t$, 0 otherwise \\
$y_{t}(\xi)$ & Number of holding flights at time $t$ in scenario $\xi$\\ 
\bottomrule
\end{tabular}
\end{table}
We next present a two-stage stochastic variant of the SAGHP, which we refer to as s-SAGHP. In the s-SAGHP, the capacity of the impacted airport to accommodate arriving aircraft is uncertain. The central modeling assumption of the s-SAGHP is that departure times of flights are chosen at the beginning of the planning horizon, before airport capacity is revealed. It is then possible that the number of flights attempting to land in some time period is in excess of capacity. In such a case, it is necessary for some flights to be held in the air until the next time period. The goal is to minimize an expected total cost of ground and air holding.

Let $\widehat{\Xi}$ represent a set of possible scenarios or realizations of the uncertain elements that affect the problem. We then model the airport capacity as a $|\mathcal{T}|$-dimensional random variable $K$, with $K_t(\widehat{\xi})$ giving the number of flights that may arrive at the airport in time period $t$ under scenario $\widehat{\xi}$. The s-SAGHP requires some estimated distribution $\widehat{\mathbb{P}}$ supported on $\widehat{\Xi}$. Then, s-SAGHP is given as follows:
\begin{subequations}
\begin{align}
\min_{x}&\left\{\sum_{f\in \mathcal{F}} C_g \left(\sum_{t \in \mathcal{T}_f}tx_{f,t} - r_f\right) + \mathbb{E}_{\widehat{\xi}\sim\widehat{\mathbb{P}}}\left[Q(x,\widehat{\xi})\right]\right\}\label{eq:first_stage_objfunc} \\
\textup{s.t. }&\eqref{eq:dSAGHP_1c}-\eqref{eq:dSAGHP_1e}\notag.
\end{align}
\label{eq:s-SAGHP_first_stage}
\end{subequations}
where the recourse $Q(x, \xi)$ is formulated as:
\begin{subequations}
\label{eq:s-SAGHP_second_stage}
\begin{align}
\min \quad
    & \sum_{t=0}^{T} C_h\, y_t(\xi)
    \label{eq:second_stage_objfunc} \\
\text{s.t.} \quad
    & y_t(\xi) \ge y_{t-1}(\xi)
    + \sum_{f\in \mathcal{F}} x_{f,t} - K_t(\xi),
    && \forall t \in \mathcal{T},
    \label{eq:second_stage_1b} \\
    & y_0(\xi) \ge \sum_{f\in \mathcal{F}} x_{f,0} - K_0(\xi),
    \label{eq:second_stage_1a} \\
    & y_t(\xi) \ge 0,
    && \forall t \in \mathcal{T}.
    \label{eq:second_stage_1c}
\end{align}
\end{subequations}
The objective function \eqref{eq:first_stage_objfunc} is the sum of ground holding delays for all flights and the expected air delays in the second stage. The first stage constraints of d-SAGHP consist of the same constraints in s-SAGHP except for the airport capacity constraints, which are imposed in the second stage. 
The second stage decision variable $y_t(\xi)$ represents the number of arrival flights joining the arrival queue or conducting airborne holding around the terminal area in scenario $\xi$. The second stage objective~\eqref{eq:second_stage_objfunc} is to minimize the airborne delay cost. As for the constraints of the second stage problem, \eqref{eq:second_stage_1b} is the second stage capacity constraint. Intuitively, the expression $y_{t-1}(\xi)+\sum_{f\in \mathcal{F}}x_{f,t}$ gives the number of flights in the terminal airspace at time $t$. At most $K_t(\xi)$ such flights will be able to land, and any remaining flights must be held to the next time period. Constraint \eqref{eq:second_stage_1a} is second stage capacity constraint at time period $0$. This constraint follows a similar format as constraint~\eqref{eq:second_stage_1b}, with the assumption that no flights are held in the air prior to the initiation of the planning horizon. The magnitude of the unit cost of airborne delays is typically greater than that for ground delays.
\subsection{Wasserstein ambiguity set for airport capacity distributions}

Let
$\widehat{\Xi} = \{\widehat{\xi}^{1}, \ldots, \widehat{\xi}^{n}\}$
denote the finite support of the estimated capacity distribution, where
$\widehat\xi^{i} \in \mathbb{R}^{|\mathcal{T}|}$ represents an observed or sampled
airport capacity trajectory over the planning horizon $T$.
The estimated distribution is given by
$
    \widehat{\mathbb{P}}
    =
    \sum_{i=1}^{n} \widehat{p}_{i}
    \delta_{\widehat\xi^{i}},
$
where $\widehat{p}_{i}$ denotes the Dirac measure associated with
capacity trajectory $\widehat\xi^{i}$, satisfying
$\widehat{p}_{i} \geq 0$ and $\sum_{i=1}^{n}\widehat{p}_{i}=1$.

We distinguish the support $\widehat{\Xi}$ of the estimated distribution from the physical
uncertainty support $\Xi$. Let $\overline{K}$
denote the upper bound on airport capacity,
respectively. The physical uncertainty support is defined as
\begin{equation}
    \Xi =
    \left\{
        \xi \in \mathbb{R}^{|\mathcal{T}|}:
        0 \leq \xi_t \leq \overline{K},
        \quad \forall t \in \mathcal{T}
    \right\}.
\label{eq:physical_support}
\end{equation}
While the support of the estimated distribution $\widehat{\Xi}$ contains a finite number of
observed capacity trajectories, the physical uncertainty support $\Xi$
contains all capacity trajectories satisfying the prescribed bounds and thus $\widehat{\Xi} \subseteq \Xi$.


We use the estimated capacity distribution $\widehat{\mathbb{P}}$  as the nominal distribution in a Wasserstein ambiguity set. In our computational experiments, we construct
this estimate using Gaussian Process Regression models, as described in
Section~\ref{sec:ScenarioGeneration}. The worst-case distribution is not restricted to the support $\widehat{\Xi}$ of the estimated distribution: probability mass may be transported from each
nominal capacity trajectory $\widehat\xi^{i}$ to any physically feasible trajectory $\xi\in\Xi$, subject to the Wasserstein transportation budget.

Let $\mathcal{M}(\Xi)$, $\mathcal{M}(\widehat{\Xi})$ denote the set of probability measures supported on
$\Xi$, $\widehat{\Xi}$ respectively. For the estimated probability measure $\widehat{\mathbb{P}} \in \mathcal{M}(\widehat{\Xi})$ and any probability measures
$\mathbb{Q} \in \mathcal{M}(\Xi)$, let
$\Pi(\widehat{\mathbb{P}},\mathbb{Q})$ denote the set of joint probability
measures on $\widehat{\Xi} \times \Xi$ whose first and second marginals are
$\widehat{\mathbb{P}}$ and $\mathbb{Q}$, respectively.

Throughout this paper, we use the $1$-Wasserstein distance with the $\ell_1$ norm as the ground metric. The $1$-Wasserstein distance measures capacity perturbations linearly in expectation, consistent with the operational measurement of airport capacity changes. The $\ell_1$ norm further makes the transportation
penalty separable across time periods, enabling our DP separation algorithm.
The corresponding 1-Wasserstein distance is defined as
\begin{equation}
    W_1(\widehat{\mathbb{P}},\mathbb{Q})
    =
    \inf_{\pi \in \Pi(\widehat{\mathbb{P}},\mathbb{Q})}
    \int_{\widehat{\Xi} \times \Xi}
    \|\widehat{\xi}-\xi\|_1
    \,\mathrm{d}\pi(\widehat{\xi},\xi).
    \label{eq:Wasserstein_distance}
\end{equation}
Accordingly, the Wasserstein ambiguity set centered at the nominal
distribution is defined as
\begin{equation}
    \mathcal{P}_{\varepsilon}\left(\widehat{\mathbb{P}}\right)
    =
    \left\{
        \mathbb{Q} \in \mathcal{M}(\Xi):
        W_1(\widehat{\mathbb{P}},\mathbb{Q})
        \leq \varepsilon
    \right\},
    \label{eq:Wasserstein_ambiguity_set}
\end{equation}
where $\varepsilon \geq 0$ denotes the Wasserstein radius.

The dr-SAGHP is similar to the s-SAGHP. However, rather than minimizing the expected ground and air holding cost with respect to the estimated distribution $\widehat{\mathbb{P}}$, the dr-SAGHP minimizes the worst case expected costs amongst the distributions in a Wasserstein ball centered at $\widehat{\mathbb{P}}$. This is formulated as follows:
\begin{align}
\min_{x} &\left\{\sum_{f\in \mathcal{F}} C_g \left(\sum_{t \in \mathcal{T}_f}tx_{f,t} - r_f\right) + \max_{\mathbb Q \in \mathcal{P}_\varepsilon\left(\widehat{\mathbb{P}}\right)} \mathbb{E}_{p}\left[Q(x,\xi)\right]\right\}
\label{eq:dr_first_stage_objfunc}\\
\textrm{s.t.}\quad
&\eqref{eq:dSAGHP_1c} - \eqref{eq:dSAGHP_1e}. \notag
\end{align}
\label{eq:dr-SAGHP_first_stage}

Following standard reformulations of Wasserstein distributionally robust
optimization problems, such as \cite{c9}, we can obtain the exact
deterministic equivalent of dr-SAGHP.

Let \(\mathcal X\) denote the first-stage feasible set defined
by~\eqref{eq:dSAGHP_1c}--\eqref{eq:dSAGHP_1e}. Given a first-stage solution $x \in \mathcal{X}$, let \(A_t(x)\) denote the number of flights scheduled to arrive in
period \(t\), and let \(A(x)\) denote the corresponding vector:
\begin{equation}
A_t(x)
:=
\sum_{\substack{f\in\mathcal F:\\t\in\mathcal T_f}}x_{f,t},
\qquad
A(x):=\bigl(A_t(x)\bigr)_{t\in\mathcal T}.
\label{eq:scheduled_arrivals}
\end{equation}

The dual feasible region of the second-stage recourse
problem~\eqref{eq:s-SAGHP_second_stage} is
\begin{equation}
\mathcal V
=
\left\{
\nu\in\mathbb R^{|\mathcal T|}:
\begin{aligned}
&\nu_t\geq0,
&&t\in\mathcal T,\\
&\nu_t-\nu_{t+1}\leq C_h,
&&t=0,\ldots,T-1,\\
&\nu_T\leq C_h.
\end{aligned}
\right\}
\label{eq:recourse_dual_region}
\end{equation}
We denote its set of extreme points by
\(\operatorname{ext}(\mathcal V)\).

\begin{proposition}[Deterministic equivalent formulation of  dr-SAGHP]
\label{prop:dr_saghp_full_support}

Under the physical uncertainty support \(\Xi\) defined in~\eqref{eq:physical_support} and the second-stage recourse problem~\eqref{eq:s-SAGHP_second_stage}, dr-SAGHP admits the following deterministic equivalent formulation:
\begin{subequations}\label{eq:dr_saghp_full_support}
\begin{align}
\min_{x,\lambda,\alpha}\quad
&
\sum_{f\in \mathcal{F}} C_g
\left(
\sum_{t\in \mathcal{T}_f}t x_{f,t}-r_f
\right)
+\lambda\varepsilon
+\sum_{i=1}^{n}\widehat p_i\alpha_i
\label{eq:dr_saghp_full_support_obj}
\\
\textnormal{s.t.}\quad
& \alpha_i
  \geq
  A(x)^\top\nu +
  \min_{\|z\|_\infty\leq\lambda}
  \left\{
    z^\top\widehat\xi^{i}
    +
    \sum_{t\in\mathcal T}
    \overline K_t[-\nu_t-z_t]_+
  \right\},
\notag\\[-1mm]
&
\forall i\in[n],\quad
\forall\nu\in\operatorname{ext}(\mathcal V),
\label{eq:dr_saghp_full_support_dro}
\\
&
x\in\mathcal X,
\label{eq:dr_saghp_full_support_first_stage}
\\
&
\lambda\geq0, \\
&\alpha_i\in\mathbb R,
\quad \forall i\in[n],
\label{eq:dr_saghp_full_support_domain}
\end{align}
\end{subequations}
where $[n]\coloneqq\{1,\ldots,n\}$.
\end{proposition}

The proof of Proposition~\ref{prop:dr_saghp_full_support} is provided in
Appendix~\ref{sec:dr_saghp_reformulation}.

\section{Decomposition Approach}\label{sec:solution_approach}

Although \eqref{eq:dr_saghp_full_support} provides an exact finite-dimensional reformulation, solving it directly is generally impractical. Constraint \eqref{eq:dr_saghp_full_support_dro} must be imposed for every sample \(i\in[n]\) and every extreme point \(\nu\in\operatorname{ext}(\mathcal V)\), resulting in \(n|\operatorname{ext}(\mathcal V)|\) constraints. The number of extreme points of \(\mathcal V\) can grow exponentially with the length of the planning horizon, so explicitly enumerating these constraints may produce a prohibitively large optimization problem. However, for a fixed, possibly infeasible solution, the search for violated constraints decomposes across samples. This motivates a delayed constraint-generation (DCG) method that begins with a restricted master problem (RMP), solves an independent separation problem for each sample, and iteratively adds only the violated extreme-point constraints.

\subsection{Delayed constraint generation approach} \label{subsec:delayed_constraint_generation}
In this section, we propose a solution approach using delayed constraint generation (also referred to as row generation in \cite{bertsimas1997introduction}) to handle the LP relaxation of dr-SAGHP in \eqref{eq:dr_saghp_full_support}. Section~\ref{subs:l-shaped} describes how this approach is incorporated into a branch-and-cut framework to solve the original dr-SAGHP with integral first-stage solutions.

Our delayed constraint generation method maintains a set of dual extreme points for each sample $i$. Let  define \(\mathcal G_i^k\subseteq\operatorname{ext}(\mathcal V)\) denote the set of dual extreme points generated for sample \(i\) generated after $k$ iterations of the algorithm. The restricted master problem (RMP) is defined as:
\begin{subequations}\label{eq:dcg_rmp}
\begin{align}
\min_{x,\lambda,\alpha}\quad
&
\sum_{f\in \mathcal{F}}C_g
\left(
\sum_{t\in \mathcal{T}_f}t x_{f,t}-r_f
\right)
+\lambda\varepsilon
+\sum_{i=1}^{n}\widehat p_i\alpha_i
\label{eq:dcg_rmp_obj}
\\
\textnormal{s.t.}\quad
&
\alpha_i
\geq
A(x)^{\top}\nu + \Psi(\nu, \lambda, \widehat\xi^{i}),
\\ \notag &\forall \nu\in\mathcal G_i^k,\quad
\forall i\in[n],
\label{eq:dcg_rmp_cuts}
\\
&
x\in\Bar{\mathcal{X}}, \\
&\lambda\geq0, \\
&\alpha_i \in \mathbb R,
\forall i\in[n].
\label{eq:dcg_rmp_domains}
\end{align}
\end{subequations}
where
\begin{equation}
\Psi(\nu, \lambda, \widehat\xi^{i}) = \min_{\|z_i\|_\infty\leq\lambda}
\left\{
z_i^\top\widehat\xi^{i}
+
\overline K
\sum_{t\in\mathcal T}
[-\nu_t-z_{it}]_+
\right\},
\end{equation}
and $\Bar{\mathcal{X}}$ denotes the relaxed first-stage feasible region, $\Bar{\mathcal{X}} = \{x: \sum_{f \in \mathcal{F}}x_{f,t} \leq K_t\quad\forall t \in \mathcal{T}, \sum_{t \in \mathcal{T}_f}x_{f,t}=1\quad\forall f \in \mathcal{F}, x_{f,t} \in [0,1] \quad \forall f \in \mathcal{F}, t \in \mathcal{T}_f \}$.
Then, for each $x \in \Bar{\mathcal{X}}$ and $\lambda \geq 0$, the separation problem is defined as follows:
\begin{equation}
\operatorname{SEP}_{i}(x, \lambda) = \max_{\nu\in\mathcal V}
\left\{
A(x)^\top\nu
+
\Psi(\nu, \lambda, \widehat\xi^{i})
\right\}.
\label{eq:dcg_separation_problem}
\end{equation}

With the RMP and separation problem defined above, we now present the DCG algorithm in Algorithm~\ref{alg:dr_saghp_dcg}, which iteratively solves the RMP to identify a first-stage solution $x \in \bar{\mathcal{X}}$ and dual multiplier $\lambda \geq 0$, and then adds constraints by solving the separation problems $\operatorname{SEP}_{i}(x,\lambda)$. Since \(\mathcal V\) is a bounded polyhedron,
\(\operatorname{ext}(\mathcal V)\) is finite. Then, standard results on delayed constraint generation, such as \citet[Section~6.3]{bertsimas1997introduction}, imply that DCG converges in a finite number of iterations.
\begin{remark}[Finite convergence of DCG]
Assuming that the master and separation problems are solved
exactly, DCG terminates after finitely many iterations and returns
a solution $x \in \Bar{\mathcal{X}}$ whose objective value is within $\eta$ of the optimal solution to the LP relaxation of dr-SAGHP.
\end{remark}
Solving the separation problem~\eqref{eq:dcg_separation_problem} is nontrivial because it is a nonconvex optimization problem involving the maximization of a piecewise-linear convex function. To solve the separation problem, \citep{gamboa2021decomposition} propose an MILP reformulation for general two-stage Wasserstein DRO models with rectangular uncertainty support. We adopt a similar MILP reformulation as a baseline separation method within our DCG framework and denote the resulting approach by DCG-MILP. For dr-SAGHP, the bounded box support with zero lower bounds further simplifies the separation problem: an optimal adversarial capacity realization can be obtained analytically once an optimal dual solution is known. Moreover, the generated constraints can be indexed solely by dual extreme points, allowing us to use DCG without explicitly generating capacity scenarios, which would be necessary under a direct application of the column-and-constraint generation (CCG) algorithm of \citep{gamboa2021decomposition}. The MILP formulation and its derivation are provided in Appendix~\ref{sec:milp_separation}.

In the following section, we develop an efficient separation algorithm using a dynamic programming (DP) approach.
\begin{algorithm}
\caption{Delayed constraint generation for the LP relaxation of dr-SAGHP}
\label{alg:dr_saghp_dcg}
\begin{algorithmic}[1]

\Require Estimated distribution
\(\{(\widehat\xi^{i},\widehat p_i)\}_{i=1}^{n}\),
Wasserstein radius \(\varepsilon\), and tolerance \(\eta>0\)

\Ensure Solution \((x^*,\lambda^*,\alpha^*)\)

\State Set \(k\gets0\), \(LB\gets-\infty\),
\(UB\gets+\infty\), and
\(\mathcal G_i^0\gets\emptyset\) for all \(i\in[n]\)

\While{\(UB-LB>\eta\)}
    \State Solve RMP~\eqref{eq:dcg_rmp} and obtain
    \((\bar x^k,\bar\lambda^k,\bar\alpha^k)\) and \(LB\)

    \For{\(i=1,\ldots,n\)}
        \State Solve the separation problem
        \(\operatorname{SEP}_i(\bar x^k,\bar\lambda^k)\)
        \State Let \(\nu_i^k\) be an optimal solution and compute
        \(\Delta_i^k
        \gets
        \operatorname{SEP}_i(\bar x^k,\bar\lambda^k)
        -\bar\alpha_i^k\)
    \EndFor

    \For{\(i=1,\ldots,n\)}
        \If{\(\Delta_i^k>0\)}
            \State
            \(\mathcal G_i^{k+1}
            \gets
            \mathcal G_i^k\cup\{\nu_i^k\}\)
        \Else
            \State
            \(\mathcal G_i^{k+1}\gets\mathcal G_i^k\)
        \EndIf
    \EndFor

    \State Set
    $
        UB\gets
        \min\left\{
            UB,\,
            LB+\sum_{i=1}^{n}\widehat p_i
            \max\{0,\Delta_i^k\}
        \right\}
    $
    \State \(k\gets k+1\)
\EndWhile

\end{algorithmic}
\end{algorithm}

\subsection{DP separation algorithm}
\label{subsec:dp_separation}

In this section, we exploit the structure of the dual feasible
region $\mathcal V$ to solve the separation problem using
dynamic programming.

For fixed $(\bar x,\bar\lambda)$, define the reward at time $t$ as
\begin{equation}
    r_{i,t}(\nu_t)
    =
    \left(A_t(\bar x)-\widehat\xi^{i}_t\right)\nu_t
    +
    \widehat\xi^{i}_t[\nu_t-\bar\lambda]_+.
\end{equation}
The following lemma expresses the separation objective as a
sum of these rewards.

\begin{lemma}\label{lemma:separation_reformulation}
For a fixed first-stage solution
$\bar x\in\bar{\mathcal X}$ and $\bar\lambda\geq0$,
\begin{equation}
    \operatorname{SEP}_{i}(\bar x,\bar\lambda)
    =
    \max_{\nu\in\mathcal V}
    \sum_{t\in\mathcal T}r_{i,t}(\nu_t).
\end{equation}
\end{lemma}

We next characterize the structure of an optimal solution.
Define
\begin{equation}
    \mathcal E
    =
    \left\{
        \nu\in\mathbb R^{|\mathcal T|}:
        \begin{aligned}
            &\nu_t\in\{0,\nu_{t+1}+C_h\},\\
            &\forall t\in\mathcal T
        \end{aligned}
    \right\},
\end{equation}
where $\nu_{T+1}=0$.
We can restrict our attention to elements of $\mathcal E$,
as formalized in Proposition~\ref{prop:dp_separation}.

\begin{proposition}
\label{prop:dp_separation}
For a fixed first-stage solution
$\bar x\in\bar{\mathcal X}$ and $\bar\lambda\geq0$,
\begin{equation}
    \operatorname{SEP}_{i}(\bar x,\bar\lambda)
    =
    \max_{\nu\in\mathcal E}
    \sum_{t\in\mathcal T}r_{i,t}(\nu_t).
    \label{eq:structured_separation_problem}
\end{equation}
\end{proposition}

The proof of Proposition~\ref{prop:dp_separation} is provided
in Appendix~\ref{sec: proof-dp-separation}.

By Proposition~\ref{prop:dp_separation}, an optimal solution
$\nu^*$ can be found by searching over $\mathcal E$.
For each $\nu\in\mathcal E$, the component $\nu_t$ is an
integer multiple of $C_h$ between $0$ and $(T-t+1)C_h$.
For each $t\in\mathcal T$, define the set of possible values as
\begin{equation}
    \mathcal S_t
    \coloneqq
    \left\{
        jC_h:j=0,\ldots,T-t+1
    \right\}.
\end{equation}

Given the current component $\nu_t$, the optimal choice of
the remaining components does not depend on the earlier
components. This allows us to construct a DP recurrence.
For $t=0,\ldots,T-1$, let $V_{i,t}(\nu_t)$ denote the largest
total reward achievable over periods $t,\ldots,T$, given
the current dual variable $\nu_t$:
\begin{equation}
\begin{aligned}
    V_{i,t}(\nu_t)
    ={}&
    \max_{\nu_{t+1},\ldots,\nu_T}
    \sum_{\tau=t}^{T}r_{i,\tau}(\nu_\tau)
    \\
    \text{s.t.}\quad
    &\nu_\tau\in\{0,\nu_{\tau+1}+C_h\},
    \\
    &\tau=t,\ldots,T-1,
    \\
    &\nu_T\in\{0,C_h\}.
\end{aligned}
\end{equation}
The value function satisfies the Bellman recursion
\citep{bellman1954theory}
\begin{equation}
\begin{aligned}
    V_{i,t}(\nu_t)
    ={}&r_{i,t}(\nu_t)
    \\
    &+
    \begin{cases}
        \displaystyle
        \max_{\nu_{t+1}\in\mathcal S_{t+1}}
        V_{i,t+1}(\nu_{t+1}),
        & \nu_t=0,
        \\[6pt]
        V_{i,t+1}(\nu_t-C_h),
        & \nu_t>0,
    \end{cases}
\end{aligned}
\label{eq:dp_bellman_recursion}
\end{equation}
for $t=0,\ldots,T-1$, with terminal condition
\begin{equation}
\begin{aligned}
    V_{i,T}(\nu_T)&=r_{i,T}(\nu_T),
    \\
    \nu_T&\in\mathcal S_T=\{0,C_h\}.
\end{aligned}
\end{equation}
The separation objective is obtained as
\begin{equation}
    \operatorname{SEP}_{i}(\bar x,\bar\lambda)
    =
    \max_{\nu_0\in\mathcal S_0}V_{i,0}(\nu_0).
\end{equation}
The resulting backward recursion solves the separation
problem exactly in $O(T^2)$ time.

\begin{theorem}
\label{prop:dp_separation_complexity}
For a fixed first-stage solution
$\bar x\in\bar{\mathcal X}$ and $\bar\lambda\geq0$,
the separation problem
$\operatorname{SEP}_{i}(\bar x,\bar\lambda)$
can be solved exactly by dynamic programming in $O(T^2)$ time.
\end{theorem}

\begin{proof}[Proof sketch]
By Proposition~\ref{prop:dp_separation}, optimizing over
$\mathcal E$ is equivalent to solving the separation problem.
The Bellman recursion maximizes the total reward over all
dual sequences in $\mathcal E$ and therefore solves the
separation problem exactly.

At each period $t<T$, every positive state has the unique
successor $\nu_{t+1}=\nu_t-C_h$, while the zero state requires
a maximization over $\mathcal S_{t+1}$.
Thus, evaluating all states at period $t$ takes
\[
    O\bigl(|\mathcal S_t|+|\mathcal S_{t+1}|\bigr)
\]
time. Since $|\mathcal S_t|=T-t+2$, summing over all periods
gives a total running time of $O(T^2)$.
\end{proof}

The corresponding DP separation algorithm is provided in
Algorithm~\ref{alg:dp_separation} in
Appendix~\ref{sec:alg_dp_separation}.

\section{Restricted support approach}\label{subs:cutting-plane}
Computational results show that the DCG-DP method developed in
Section~\ref{sec:solution_approach} solves dr-SAGHP exactly with
substantial runtime reductions relative to the DCG-MILP baseline.
To obtain further computational savings, we develop a
restricted-support approximation, referred to as rs-dr-SAGHP,
that restricts adversarial distributions to the support
$\widehat{\Xi}$ of the estimated distribution
$\widehat{\mathbb{P}}$. We solve the resulting model using
Kelley's cutting-plane method, providing a faster approximate
approach to dr-SAGHP.
\subsection{Preliminaries}
Our approach to solving rs-dr-SAGHP can be applied to any distributionally robust mixed-integer linear program wherein there is a fixed, finite set of scenarios that on which all elements of the ambiguity set are supported. Let $\mathcal{\widehat{P}}_{\varepsilon}(\widehat{\mathbb{P}})$ be the elements of the Wasserstein ball of radius $\varepsilon$ centered on $\widehat{\mathbb{P}}$ that are supported on $\widehat{\Xi}$, i.e., $\mathcal{\widehat{P}}_{\varepsilon}(\widehat{\mathbb{P}}) = \widehat{P}_{\varepsilon}(\widehat{\mathbb{P}}) \cap \mathcal{M}(\widehat{\Xi})$.
Consider the following general setup:
\begin{equation}\label{tsdrmilp: formulation_1}
\begin{aligned}
    \min_{\mathbf{x}}&\quad\mathbf{c}^{\top}\mathbf{x} + \max_{p \in \mathcal{\widehat{P}}_{\varepsilon}(\widehat{\mathbb{P}})} \; \mathbb{E}_{\xi \sim p}[h(\mathbf{x},\xi)] \\
    \text{s.t.}\quad &\mathbf{B}\mathbf{x}\leq\mathbf{b} \\
    &\mathbf{x} \in \{0,1\}^{N}
\end{aligned}
\end{equation}
where $\mathbf{c}$ is an $N$-dimensional vector, $\mathbf{B}$ is an $M \times N$ matrix, $\mathbf{b}$ is an $M$-dimensional vector, and $h$ is a recourse function $h\colon \mathbb{R}^{N} \times \widehat{\Xi} \to \mathbb{R}$. The assumption that $\mathbf{x}$ is a vector of pure binary variables is made solely for convenience; a similar approach applies to general integer or mixed variables. 

In addition, we define
\begin{equation*}
    f(\mathbf{x}) = \max_{p \in \mathcal{\widehat{P}}_{\varepsilon}(\widehat{\mathbb{P}})}\mathbb{E}_{\xi\sim p} \left[h(\mathbf{x},\xi)\right]
\end{equation*} and assume that $h(\cdot,\xi)$ is convex, which implies that $f$ will likewise be convex; for each fixed distribution $p \in \mathcal{\widehat{P}}_{\varepsilon}(\widehat{\mathbb{P}})$ the function $\mathbf{x} \mapsto \mathbb{E}_p\left[h(\mathbf{x},\xi)\right] = \sum_{i=1}^{|\widehat{\Xi}|}p_ih(\mathbf{x},\widehat\xi^{i})$ is convex, as it is sum of convex functions multiplied by non-negative values, and taking the maximum of a collection of convex functions preserves convexity. The recourse $h(\cdot,\xi)$ is typically true when the second stage is the objective value of a linear program, as in the SAGHP variants. However, this is not typically true when the second stage contains some integer variables, in which case adjustments may be required to our approach. We leave such an extension as an avenue of future work.

Finally, we assume that $\mathbb{E}_{p}\left[h(\mathbf{x},\xi) \right]$ is finite for each $p \in \mathcal{\widehat{P}}_{\varepsilon}(\widehat{\mathbb{P}})$ and for each $\mathbf{x} \in \{0,1\}^{N}$ satisfying $\mathbf{B}\mathbf{x} \leq \mathbf{b}$. This is similar to the assumption of relatively complete recourse, often made in stochastic programming (see, e.g. chapter 5 of \cite{birge1997introduction}). In the case that the second stage is the objective value of a linear program, our method could be readily extended to the case where relatively complete recourse does not apply through the addition of standard Benders' feasibility cuts (e.g., section 6.5 of \citealt{bertsimas1997introduction}).
We first develop a decomposition algorithm to solve the linear programming relaxation of Problem~\eqref{tsdrmilp: formulation_1}. Our approach is a specialization of Kelley's cutting plane algorithm \citep{kelly1960cutting} that leverages the structure of the distributionally robust optimization problem to efficiently produce cuts.
We recast the LP relaxation of formulation~\eqref{tsdrmilp: formulation_1} in its epigraph form:
\begin{subequations}
\begin{align}
    \min_{\mathbf{x},\theta}\quad&\mathbf{c}^{\top}\mathbf{x} + \theta \\
    \text{s.t.}\quad &\mathbf{B}\mathbf{x}\leq\mathbf{b} \\
    & \theta \geq f(\mathbf{x}),\label{eq: epigraph}\\
    & \mathbf{x} \in [0,1]^{N}
\end{align}
\end{subequations}\label{tsdrmilp: formulation_2}
Let $L$ be a constant such that $f(\mathbf{x}) \geq L$ for all $\mathbf{x} \in [0,1]^{N}$, which must exist due to the convexity of $f$. Since $f$ is convex, then, as in \cite{kelly1960cutting}, we can approximate the constraints~\eqref{eq: epigraph} with a collection of constraints of the form:
\begin{equation}
\begin{aligned}
    \theta &\geq L\\
    \theta &\geq f(\mathbf{x}^k)+(\mathbf{s}^k)^\intercal(\mathbf{x} - \mathbf{x}^k) \textup{ for } k \in \{1,\ldots, K\} 
\end{aligned}\label{eq:epigraph-approx}
\end{equation}
where $\{\mathbf{x}^1,\ldots,\mathbf{x}^K\}$ are a collection of feasible test points, and where each vector $\mathbf{s}^k$ is a corresponding subgradient of $f$ at $\mathbf{x}^k$. The points $\mathbf{x}^1,\ldots,\mathbf{x}^K$ are generated iteratively. In particular, once $K$ points have been generated, the point $\mathbf{x}^{K+1}$ is the optimal solution to the problem produced by the first $K$ cuts.
Replacing constraints~\eqref{eq: epigraph} with the approximation~\eqref{eq:epigraph-approx} produces a relaxation of the problem. Thus, if we let $\widehat{z}^K$ be the optimal objective value of the relaxation with $K$ cuts added and we let $z^*$ be the optimal objective value of problem~\eqref{tsdrmilp: formulation_1}, we have that $\widehat{z}^K\leq z^*$. Furthermore,
since $\mathbf{x}^K$ is feasible, the objective value $f(\mathbf{x}^K) \geq z^*$. Thus, in each iteration, the difference $f(\mathbf{x}^K) - \widehat{z}^K$ provides an optimality gap, and the algorithm may be terminated once this gap is sufficiently small. See Algorithm~\ref{alg:kelly}. For any fixed optimality gap, the algorithm converges within a finite number of iterations 
\citep{kelly1960cutting}.
\begin{algorithm}
\caption{Kelley's cutting plane algorithm.}\label{alg:kelly}
    \begin{algorithmic}[1]
        \State Input: tolerated optimality gap $\eta$.

        \For{$K=1,2,\ldots$}
        \State Solve:
        \Statex \hspace{\algorithmicindent}
        $\begin{aligned}
            \min_{\mathbf{x},\theta}\quad&\mathbf{c}^{\top}\mathbf{x} + \theta \\
    \text{s.t.}\quad &\mathbf{B}\mathbf{x}\leq\mathbf{b} \\
    \theta &\geq L\\
    \theta &\geq f(\mathbf{x}^k)+(\mathbf{s}^k)^\intercal(\mathbf{x} - \mathbf{x}^k)\quad\forall k \in \{1,\ldots, K-1\}\\
    \mathbf{x} &\in [0,1]^{N}\\
    \end{aligned}$
    
    \State Let $\mathbf{x}^K$ be the solution and $\widehat{z}^K$ be the optimal value.
    \If{$\mathbf{c}^{\top}\mathbf{x} + f(\mathbf{x}^K) - \widehat{z}^K \leq \eta$}
    \State\Return $\mathbf{x}^K$
    \EndIf
        \State Calculate a subgradient $\mathbf{s}^K \in \partial f(\mathbf{x}^K)$

        \EndFor
    \end{algorithmic}

\end{algorithm}

Our method leverages the structure of the function $f$ arising from the distributionally robust setting in order to efficiently calculate the subgradients. We apply Danskin's Theorem to relate the subdifferential of $f$ to the subdifferentials of the recourse function in the scenarios. Given a feasible first-stage solution $\mathbf{x}$, let $\mathcal{\widehat{P}}^*_{\varepsilon}(\mathbf{x})$ be the set of maximizers of $\max_{p \in \mathcal{\widehat{P}}_{\varepsilon}(\widehat{\mathbb{P}})} \mathbb{E}_{\xi \sim p}\left[h(\mathbf{x},\xi)\right]$.
\begin{proposition}\label{prop: subgrad}
    Let $\mathbf{x}$ be any point in $\mathbb{R}^{N}$ such that $\mathbf{B}\mathbf{x} \leq \mathbf{b}$ and $\mathbf{x} \geq 0$. Then,
    \begin{align*}
        \partial f(\mathbf{x}) = \left\{\sum_{i=1}^{|\widehat{\Xi}|} p^*_i \mathbf{s}^{i} : p^* \in \mathcal{\widehat{P}}^*_{\varepsilon}(\mathbf{x}), \mathbf{s}^{i} \in \partial_{\mathbf{x}}h(\mathbf{x},\widehat\xi^{i})\right\}
    \end{align*}
    where $\partial_{\mathbf{x}} h(\mathbf{x},\xi)$ denotes the subdifferential of $h(\cdot, \xi)$ at $\mathbf{x}$.
\end{proposition}
A proof of Proposition~\ref{prop: subgrad} is given in Appendix~\ref{sec: proof-subgrad}. In essence, Proposition~\ref{prop: subgrad} states that calculating a subgradient in $\partial f(\mathbf{x})$ can be accomplished by calculating the corresponding worst-case distribution $p^*(\mathbf{x})$, calculating a subgradient of the recourse function in $\partial_{\mathbf{x}} h(\mathbf{x},\widehat\xi^{i})$ for each scenario, and then combining the two.

We assume that a subgradient in $\partial_{\mathbf{x}} h(\mathbf{x},\widehat\xi^{i})$ can be computed without excessive computational effort. In cases such as the SAGHP variants wherein the recourse function $h(\mathbf{x},\widehat\xi^{i})$ is the objective value of a second stage linear program, it follows from standard results in stochastic programming that the subgradient can be calculated from the dual optimal variables of the second stage (see, for example, chapter 2 of \citealt{shapiro2021lectures}). 

For dr-SAGHP, let $\nu^*(\mathbf{x},\xi)$ denote an optimal
solution to the dual recourse problem. A subgradient
$\mathbf{s}(\xi)$ of $h(\cdot,\xi)$ at $\mathbf{x}$ is given by
$s_{f,t}(\xi)=\nu_t^*(\mathbf{x},\xi)$ for all
$f\in\mathcal F$ and $t\in\mathcal{T}$.
Once a worst-case distribution
$p^*(\mathbf{x}^k)\in
\mathcal{\widehat{P}}_{\varepsilon}^*(\mathbf{x}^k)$
is identified, the cut generated at the test point
$\mathbf{x}^k$ is given by:
\begin{align*}
\theta &\geq \sum_{i=1}^{|\widehat{\Xi}|}p_i^{*}(\mathbf{x}^k)h(\mathbf{x}^k,\widehat\xi^{i})\\
&+ \sum_{i=1}^{|\widehat{\Xi}|}p_i^{*}(\mathbf{x}^k)\sum_{t=0}^{T}\nu^*_{t}(\mathbf{x}^k,\widehat\xi^{i})\left(\sum_{f\in \mathcal{F}}(x_{f,t} - x^k_{f,t})\right)
\end{align*}
By strong duality, $h(\mathbf{x}^k,\widehat\xi^{i}) = \sum_{t=0}^{T}
\nu^*_t(\mathbf{x}^k,\widehat\xi^{i})\left(\sum_{f\in \mathcal{F}}x^k_{f,t}- \widehat\xi^{i}_t\right)$, so this can be simplified to:
\begin{align*}
\theta \geq \sum_{i=1}^{|\widehat{\Xi}|}p_i^{*}(\mathbf{x}^k)\sum_{t=0}^{T}\nu^*_{t}(\mathbf{x}^k, \widehat\xi^{i})\left(\sum_{f\in \mathcal{F}}x_{f,t}-\widehat\xi^{i}_t\right),
\end{align*}

\subsection{Calculating the worst-case distribution}\label{subs:worst-case-dist}
Using the fact that the set $\widehat{\Xi}$ is discrete, we provide an efficient method for calculating the inner maximization problem $\max_{p \in \mathcal{\widehat{P}}_{\varepsilon}(\widehat{\mathbb{P}})} \mathbb{E}_{\xi \sim p}\left[h(\mathbf{x},\xi) \right]$ for a fixed, feasible first stage decision $\mathbf{x}$. Throughout this section, we treat the first stage decision $\mathbf{x}$ as fixed. It follows from the definition of the Wasserstein ambiguity set that this problem can be formulated as follows:
\begin{subequations}\label{eq:dr-general-inner_integration}
\begin{align}
\max_{\pi \in [0,1]^{|\widehat{\Xi}|\times |\widehat{\Xi}|}} &\; \sum_{i=1}^{\abs{\widehat{\Xi}}}\sum_{j=1}^{\abs{\widehat{\Xi}}}\pi_{i,j}\, h(\mathbf{x},\widehat{\xi}^{(j)})
\label{eq:dr_general_inner_integration_objfunc}\\
\text{s.t.}\quad
&\sum_{i=1}^{\abs{\widehat{\Xi}}}\sum_{j=1}^{\abs{\widehat{\Xi}}}\pi_{i,j}\mathbf{d}_{i,j} \leq \varepsilon,
\label{eq:dr_general_inner_integration_1a}\\
&\sum_{j=1}^{\abs{\widehat{\Xi}}}\pi_{i,j} = \widehat{p}_{i}, \quad \forall i,
\label{eq:dr_general_inner_integration_1b}\\
&\pi_{i,j} \ge 0, \quad \forall i,j,
\label{eq:dr_general_inner_integration_1c}
\end{align}
\end{subequations}
Given a solution $\pi^*$, the corresponding optimal distribution is given by $p^*_j = \sum_{i=1}^{\abs{\widehat{\Xi}}} \pi^*_{i,j}$ for each $j \in \{1,\ldots, |\widehat{\Xi}|\}$.
The dual of problem~\eqref{eq:dr-general-inner_integration} follows:
\begin{subequations}\label{eq:dr-general-inner_dual}
\begin{align}
\min_{\lambda \in \mathbb{R},\alpha \in \mathbb{R}^{|\widehat{\Xi}|}}\quad
& \lambda\varepsilon + \sum_{i=1}^{\abs{\widehat{\Xi}}}\alpha_{i}\widehat{p}_{i}
\label{eq:dr_general_inner_dual_objfunc}\\
\text{s.t.}\quad
& \alpha_{i} + \lambda \mathbf{d}_{i,j} \ge h(\mathbf{x},\widehat{\xi}^{(j)}),
\quad \forall i, j,
\label{eq:dr_general_inner_dual_1a}\\
& \lambda \ge 0.
\label{eq:dr_general_inner_dual_1b}
\end{align}
\end{subequations}
We can reduce the dual problem~\eqref{eq:dr-general-inner_dual} to a univariate optimization problem. First, note that for any $\lambda \geq 0$, there exists $\alpha \in \mathbb{R}^{|\widehat{\Xi}|}$ such that $(\lambda, \alpha)$ is feasible for \eqref{eq:dr-general-inner_dual}. In particular, if we let each $\alpha_i$ be any value no smaller than $\max_{j \in \{1,\ldots, \abs{\widehat{\Xi}}\}}
(h(\mathbf{x},\widehat{\xi}^{(j)}) - \lambda \mathbf{d}_{i,j})$, then $(\lambda, \alpha)$ is feasible. Furthermore, since the coefficient of $\alpha_i$ is given by $\widehat{p}_i$, which is non-negative, the optimal choice of $\alpha_i$ given fixed $\lambda$ is given by $\alpha^*_i(\lambda) := \max_{j \in \{1,\ldots,|\widehat{\Xi}|\}}(h(\mathbf{x},\widehat{\xi}^{(j)}) - \lambda \mathbf{d}_{i,j})$. Thus, the problem reduces to that of minimizing the function $\phi:\mathbb{R} \to \mathbb{R}$ given by:
\begin{align}
\phi&(\lambda)
:=~\lambda\varepsilon + \sum_{i=1}^{\abs{\widehat{\Xi}}} \alpha_i^*(\lambda) \widehat{p}_i\notag\\
&=
\lambda \varepsilon
+ \sum_{i=1}^{\abs{\widehat{\Xi}}} \widehat{p}_{i}\,
\max_{j\in\{1,\ldots,\abs{\widehat{\Xi}}\}}
\left(
h(\mathbf{x},\widehat{\xi}^{(j)})
- \lambda \mathbf{d}_{i,j}\right).
\label{eq:one_dimensional_inner_dual}
\end{align}

\begin{proposition}\label{prop:1-dim-dual}
Let $\lambda^*$ be a solution to $\min_{\lambda \geq 0}\phi(\lambda)$. Then $(\lambda^{*}, \alpha^{*}(\lambda^*))$ is an optimal solution to problem \eqref{eq:dr-general-inner_dual}
\end{proposition}
The function $\phi$ is piecewise-linear convex, as each function $\lambda \mapsto  \max_{j \in \{1,\ldots, |\widehat{\Xi}|\}}
(h(\mathbf{x},\widehat{\xi}^{(j)}) - \lambda \mathbf{d}_{i,j})$
is a piecewise-linear convex function and each $\widehat{p}_i \geq 0$.
A proof of Proposition~\ref{prop:1-dim-dual} is given in Appendix~\ref{sec:proof-1-dim-dual}.

We therefore develop a bisection scheme, given in Algorithm~\ref{alg:dual-bisection}, to compute the optimal dual multiplier $\lambda^{*}$. The algorithm requires as a starting point lower and upper bounds $\lambda_{\mathrm{hi}}$ and $\lambda_{\mathrm{lo}}$ such that there exists a minimizer $\lambda^*$ in the interval $[\lambda_{\mathrm{lo}}, \lambda_{\mathrm{hi}}]$. The algorithm also requires a prespecified tolerance $\delta$, and produces an solution that is within $\delta$ of an optimal solution.

At each iteration, Algorithm~\ref{alg:dual-bisection} chooses a test point $\lambda = (\lambda_{\mathrm{hi}}+\lambda_{\mathrm{lo}})/2$ at the middle of the interval and evaluates a subgradient $g(\lambda)\in\partial\phi(\lambda)$. If the subgradient is positive, this implies that there exists an optimal value $\lambda^*$ such that $\lambda^* \leq \lambda$, whence we can update the upper bound $\lambda_{\mathrm{hi}}$ to $\lambda$. Similarly, if the subgradient is negative, the lower bound $\lambda_{\mathrm{lo}}$ may be updated to $\lambda$. If the subgradient is identically zero, then $\lambda$ is optimal. This continues until the lower bound and upper bound within the specified tolerance, or an optimal solution is found.
\begin{proposition}\label{prop: dual-converge}
    Algorithm~\ref{alg:dual-bisection} stops after at most $\lceil\log_2(\frac{\lambda_{\mathrm{hi}}-\lambda_{\mathrm{lo}}}{\delta})\rceil$ iterations and returns a value $\lambda$ such that $(\lambda-\delta,\lambda+\delta) \cap \argmin_{\lambda \geq 0} \phi(\lambda)$ is non-empty. Here, $\lambda_{\mathrm{hi}}$ and $\lambda_{\mathrm{lo}}$ refer to the initial values of these variables input into Algorithm~\ref{alg:dual-bisection}.
\end{proposition}
A proof of Proposition~\ref{prop: dual-converge} is given in Appendix~\ref{sec:proof-dual-converge}.
In order to implement Algorithm~\ref{alg:dual-bisection}, we must identify the initial valid lower bounds and upper bounds. By definition, 0 is a valid lower bound; we can provide a valid upper bound that can be readily computed.
\begin{algorithm}
\caption{Dual Bisection Algorithm}
\label{alg:dual-bisection}
\begin{algorithmic}[1]
    \State Inputs: first stage feasible solution $\mathbf{x} \in [0,1]^{N}$; tolerance $\delta > 0$; values $\lambda_{\mathrm{lo}}$ and $\lambda_{\mathrm{hi}}$ such that $[\lambda_{\mathrm{lo}}, \lambda_{\mathrm{hi}}] \cap \argmin_{\lambda \geq 0} \phi(\lambda)$ is non-empty.
    \While{$\lambda_{\mathrm{hi}} - \lambda_{\mathrm{lo}} > 2\delta$}
        \State $\lambda \gets \tfrac{1}{2}(\lambda_{\mathrm{lo}} + \lambda_{\mathrm{hi}})$
    \State For each $i \in \{1,\ldots, |\widehat{\Xi}|\}$, choose $j^*(i,\lambda) \in \argmax_{j\in\{1,\ldots,|\widehat{\Xi}|\}}\{h(\mathbf{x},\widehat{\xi}^{(j)})-\lambda \mathbf{d}_{i,j}\}$
        \If{$\varepsilon - \sum_{i=1}^{\abs{\widehat{\Xi}}} \widehat{p}_i\, \mathbf{d}_{i,j^*(i,\lambda)} > 0$}
            \State $\lambda_{\mathrm{hi}} \gets \lambda$ 
        \ElsIf{$\varepsilon - \sum_{i=1}^{\abs{\widehat{\Xi}}} \widehat{p}_i\, \mathbf{d}_{i,j^*(i,\lambda)} < 0$}
            \State $\lambda_{\mathrm{lo}} \gets \lambda$
        \Else
            \State Return $\lambda$
        \EndIf
    \EndWhile
    \State $\lambda^{{*}} \gets \tfrac{1}{2}(\lambda_{\mathrm{lo}} + \lambda_{\mathrm{hi}})$
    \State \textbf{return} $\lambda^{{*}}$
\end{algorithmic}
\end{algorithm}

\begin{proposition}\label{thm:optimal_interval}
If the Wasserstein ambiguity radius $\varepsilon > 0$, then
all optimal solutions to $\min_{\lambda \geq 0}\phi(\lambda)$ are less than or equal to
\[
  \frac{
    \max_{j\in \{1,\ldots,|\widehat{\Xi}|\}} h(\mathbf{x},\widehat{\xi}^{(j)}) - \min_{i\in \{1,\ldots,|\widehat{\Xi}|\}} h(\mathbf{x},\widehat\xi^{i})
  }{
    \min_{i,j:\, \mathbf{d}_{i,j}\neq 0} \mathbf{d}_{i,j}
  }
\]
\end{proposition}

While Algorithm~\ref{alg:dual-bisection} yields an optimal dual multiplier $\lambda^{*}$, the cutting planes applied in constraint~\eqref{eq:epigraph-approx} require the primal optimal solution to $\max_{p \in \mathcal{\widehat{P}}_{\varepsilon}(\widehat{\mathbb{P}})}\mathbb{E}_{\xi \sim p}[h(\mathbf{x},\xi)]$. For this reason, we develop a primal recovery procedure, given in Algorithm~\ref{alg:primal-recovery}. Let $\widehat{\pi}$ denote the transportation plan constructed by Algorithm ~\ref{alg:primal-recovery} using $\lambda^{*}$. Intuitively, the transportation plan $\widehat{\pi}$ is constructed in order to satisfy complementary slackness with the dual solution $(\lambda^*, \alpha^*(\lambda^*))$. The probability mass at some scenario $\widehat\xi^{i}$ is transported only to scenarios $\widehat{\xi}^{(j)}$ for which the value $h(\mathbf{x},\widehat{\xi}^{(j)}) - \lambda^*\mathbf{d}_{i,j}$ is maximal; this is done to satisfy the complementary slackness condition $\widehat{\pi}_{i,j}(\alpha^{*}_{i}(\lambda^{*}) + h(\mathbf{x},\widehat{\xi}^{(j)}) - \lambda^* \mathbf{d}_{i,j}) = 0$. If $\lambda^* = 0$ then ties are broken by moving mass to the closest point. If $\lambda^* > 0$, then the transportation plan is carefully chosen to ensure that the resulting transportation costs are exactly $\varepsilon$, in order to satisfy the complementary slackness condition $\lambda^*(\varepsilon -\sum_{i=1}^{|\widehat{\Xi}|}\sum_{j=1}^{|\widehat{\Xi}|}\widehat{\pi}_{i,j}\mathbf{d}_{i,j}) = 0$. In either case, we can show that the resulting plan is feasible and satisfies complementary slackness with $(\lambda^*,\alpha^*(\lambda^*))$. 

\begin{theorem}\label{thm:primal_recovery}
 Let $\lambda^*$ be an optimal solution to the dual problem~\eqref{eq:one_dimensional_inner_dual}. Then $\widehat{\pi}$ is a primal optimal solution to~\eqref{eq:dr-general-inner_integration}.
\end{theorem}

\begin{algorithm}[htbp]
\caption{Primal recovery from an optimal dual multiplier $\lambda^{*}$}
\label{alg:primal-recovery}

\begin{algorithmic}[1]
    \State \textbf{Input:} Feasible first-stage solution $\mathbf{x} \in [0,1]^{N}$, optimal dual multiplier $\lambda^{*}\ge 0$.
    \State \textbf{Output:} A transportation plan $\widehat{\pi}\in[0,1]^{|\widehat{\Xi}|\times |\widehat{\Xi}|}$.

    \Statex
    \State \textbf{Step 1.}
    For each $i\in\{1,\dots,|\widehat{\Xi}|\}$, define
    \Statex \hspace{\algorithmicindent}
    \[
    \mathcal{A}_i(\lambda^{*}) :=
    \argmax_{j\in\{1,\dots,|\widehat{\Xi}|\}}
    \left\{
    h(\mathbf{x},\widehat{\xi}^{(j)})-\lambda^{*}\mathbf{d}_{i,j}
    \right\}.
    \]

    \State Select two representatives:
    \Statex \hspace{\algorithmicindent}
    \[
    \underline{j}(i)\in
    \argmin_{j\in \mathcal{A}_i(\lambda^{*})} \mathbf{d}_{i,j},
    \qquad
    \overline{j}(i)\in
    \argmax_{j\in \mathcal{A}_i(\lambda^{*})} \mathbf{d}_{i,j}.
    \]

    \Statex
    \State \textbf{Step 2.}
    Define the lower and upper achievable transportation costs:
    \Statex \hspace{\algorithmicindent}
    \[
    \underline{\gamma}
    :=
    \sum_{i=1}^{|\widehat{\Xi}|}\widehat{p}_i\,
    \mathbf{d}_{i,\underline{j}(i)},
    \qquad
    \overline{\gamma}
    :=
    \sum_{i=1}^{|\widehat{\Xi}|}\widehat{p}_i\,
    \mathbf{d}_{i,\overline{j}(i)}.
    \]

    \State Set
    \Statex \hspace{\algorithmicindent}
    \[
    \beta :=
    \begin{cases}
    \dfrac{\varepsilon-\underline{\gamma}}
    {\overline{\gamma}-\underline{\gamma}},
    & \text{if }\lambda^{*}>0
    \text{ and }\overline{\gamma}>\underline{\gamma},\\[10pt]
    0,
    & \text{if }\lambda^{*}=0
    \text{ or }\overline{\gamma}=\underline{\gamma}.
    \end{cases}
    \]

    \Statex
    \State \textbf{Step 3.}
    For each $i \in \{1,\ldots,|\widehat{\Xi}|\}$ and
    $j \in \{1,\ldots,|\widehat{\Xi}|\}$, set
    \Statex \hspace{\algorithmicindent}
    \[
    \widehat{\pi}_{i,j}:=
    \begin{cases}
    (1-\beta)\widehat{p}_i,
    & \text{if }\lambda^{*}>0
    \text{ and } j=\underline{j}(i),\\
    \beta\,\widehat{p}_i,
    & \text{if }\lambda^{*}>0
    \text{ and } j=\overline{j}(i),\\
    \widehat{p}_i,
    & \text{if }\lambda^{*}=0
    \text{ and } j=\underline{j}(i),\\
    0,
    & \text{otherwise}.
    \end{cases}
    \]

    \State \Return $\widehat{\pi}$.
\end{algorithmic}
\end{algorithm}

Algorithm~\ref{alg:dual-bisection} and Algorithm~\ref{alg:primal-recovery} together provide an approach to find a distribution $p^* \in \mathcal{\widehat{P}}^*_{\varepsilon}(\mathbf{x})$: first, Algorithm~\ref{alg:dual-bisection} is used to find an optimal dual multiplier $\lambda^*$, and second, the primal optimal solution is reconstructed using Algorithm~\ref{alg:primal-recovery}. In particular, if $\widehat{\pi}$ is the output of Algorithm~\ref{alg:primal-recovery}, then the desired distribution $p^*$ is given by $p^*_j = \sum_{i=1}^{|\widehat{\Xi}|}\widehat{\pi}_{i,j}$ for each $j \in \{1,\ldots,|\widehat{\Xi}|\}$. We apply this in combination with Proposition~\ref{prop: subgrad} to quickly find subgradients of $f$ for use in Kelley's cutting plane algorithm (Algorithm~\ref{alg:kelly}). 

\section{Branch-and-cut implementation}
\label{subs:l-shaped}

The delayed constraint generation and Kelley's cutting-plane procedures
developed in Sections~\ref{subsec:delayed_constraint_generation}
and~\ref{subs:cutting-plane}, respectively, can be incorporated into
a variant of the integer L-shaped method~\citep{LAPORTE1993133}
and implemented within a branch-and-cut framework. We specify
a global optimality-gap tolerance $\epsilon_{\mathrm{opt}}$ and
a tolerance $\eta$ for solving the linear programming relaxations.

At each node of the branch-and-cut tree, the node's linear programming
relaxation is solved to tolerance $\eta$ using delayed constraint generation
(Algorithm~\ref{alg:dr_saghp_dcg}) for dr-SAGHP,
or Kelley's cutting-plane procedure
(Algorithm~\ref{alg:kelly}) for the rs-dr-SAGHP.
For dr-SAGHP, the separation problem is solved
using the DP separation algorithm in Algorithm~\ref{alg:dp_separation}, as detailed
in Section~\ref{subsec:dp_separation}. For rs-dr-SAGHP, the worst-case distribution and the corresponding
subgradients are obtained using the procedures in
Section~\ref{subs:worst-case-dist} and
Proposition~\ref{prop: subgrad}. Cuts generated at a node are globally
valid for the corresponding formulation and are retained for
subsequent nodes.

The node relaxations provide lower bounds, while feasible integer
first-stage decisions, together with valid upper bounds on their
worst-case recourse costs, provide incumbent upper bounds.
The branch-and-cut procedure terminates when the global optimality
gap is at most $\epsilon_{\mathrm{opt}}$.

\section{Numerical Experiments and Discussion}\label{sec:Experiment}
In this section, we report numerical experiments that assess the performance of the dr-SAGHP under discrepancies between predicted and realized airport-capacity distributions. Our experimental setup is motivated by operational settings in which prediction errors and capacity degradations, potentially driven by climate change, create such discrepancies. To capture temporal dependence in airport capacities, we use Gaussian Process Regression (GPR)~\citep[e.g.][]{williams1995gaussian} to model the joint capacity distribution across time periods. The fitted GPR model is used to generate capacity scenarios for the experiments and is not required by the proposed problem formulation or solution algorithm. We first introduce an expanding-window cross-validation method to calibrate the Wasserstein ambiguity set radius using historical operational and weather data, and evaluate the performance of the dr-SAGHP against the s-SAGHP and robust SAGHP (r-SAGHP) on a real-world dataset. Here, we note that GDPs are implemented only on severe-weather
days in the held-out test set.

To further examine potential long-term capacity changes associated with climate change, such as more frequent and more intense thunderstorms, we conduct a sensitivity analysis by adjusting the mean and covariance of the trained GPR models. These adjustments generate testing scenarios representing different levels of capacity degradation and variability, allowing us to compare the performance of dr-SAGHP, s-SAGHP and r-SAGHP under simulated distributional shifts.
\subsection{Experiment setup}\label{sec:ExperimentSetUp}
Flight schedules are obtained from the Bureau of Transportation Statistics (BTS)~\citep{data1}. 
To evaluate the Ground Holding Problems (GHPs), we use the historical arrival performance at Denver International Airport (DEN) to evaluate the performance of dr-SAGHP. 
The airborne-to-ground delay cost ratio is set to three; that is, the unit airborne delay cost $C_h$ is three times the unit ground delay cost $C_g$.  
For each connecting flight pair $(f_1,f_2)$, the minimum turnaround time $S_{f_1,f_2}$ is fixed at three time periods (i.e., 45 minutes).

In the next following sections of computational experiments on out-of-sample performance on historical data and the simulated sensitivity analysis, we compare four models: dr-SAGHP, its fixed-finite-support approximation rs-dr-SAGHP, the stochastic model s-SAGHP, and the robust optimization benchmark r-SAGHP. We solve dr-SAGHP using either delayed constraint generation with MILP separation (DCG-MILP) in Appendix ~\ref{sec:milp_separation} or delayed constraint generation with DP separation (DCG-DP) in Section~\ref{subsec:dp_separation}. We solve rs-dr-SAGHP using Kelley's cutting-plane method within the integer L-shaped framework (K-IL) in Section~\ref{subs:cutting-plane}.

We include r-SAGHP as a robust optimization benchmark that minimizes the sum of ground-delay cost and the worst-case recourse cost over
a prescribed set of capacity trajectories. We construct a box
uncertainty set by setting its lower bound in each time period
to the $q$-quantile of the corresponding marginal GPR predictive
capacity distribution. For example, $q=0.1$ sets each lower bound
to the predicted 10th-percentile capacity. The box allows capacities to vary within their period-specific bounds without enforcing the temporal dependence represented by the joint predictive distribution. Since the recourse cost is nonincreasing in each capacity component, the worst-case trajectory is attained when every period's capacity equals its lower bound. Consequently, r-SAGHP can be solved using
the deterministic SAGHP formulation with this capacity trajectory.
A smaller $q$ lowers these bounds and provides protection against
more severe capacity reductions, making the model more conservative. The derivation is provided in
Appendix~\ref{sec:robust_saghp}.

\subsection{Airport capacity scenario generation via Gaussian Process Regression}\label{sec:ScenarioGeneration}
Similarly to \citet{murcca2018predicting}, given a testing date and the associated weather features, we use a fitted GPR model
to produce a posterior distribution for the joint airport capacities over the planning horizon. In our implementation, the GPR model is trained using only weather features. Our main contribution lies not in the predictive model itself, but in explicitly incorporating the resulting uncertainty into the downstream dr-SAGHP framework. Importantly, the GPR model is used solely for scenario generation in the numerical experiments and does not affect the problem formulation or the solution algorithm.

Let $w$ denote the vector of weather factors and let $\psi(\cdot)$ map weather factors to airport capacity. 
We model $\psi(\cdot)$ as a Gaussian process with mean function $m(\cdot)$ and covariance kernel $\mathcal{K}(\cdot,\cdot)$:
\begin{equation}
    \psi(w) \sim \mathcal{GP}\!\left(m(w),\, \mathcal{K}(w,w')\right).
\end{equation}
Let the training inputs (historical weather factors) and outputs (observed capacities at DEN) be denoted by $\widehat{w}$ and $\widehat{D}$, respectively. 
For a given day the estimated distribution $\widehat{\mathbb{P}}$ is generated as follows. Let the testing inputs across the $T$ time periods of the specified day be $\mathbf{w}=\{w_1,\dots,w_T\}$. 
Under the standard GPR model with i.i.d.\ Gaussian observation noise of variance $\sigma_d^2$, we have the mean and covariance of the posterior distribution as:
\begin{equation}
\begin{aligned}
\Tilde{\mathbf{\mu}} &= \mathbb{E}\!\left[\psi(\mathbf{w}) \mid \widehat{D}\right]
\\
&=m(\mathbf{w})
+ \mathcal{K}(\mathbf{w},\widehat{w})
\Sigma^{-1}
\left(\widehat{D} - m(\widehat{w})\right),
\end{aligned}
\end{equation}

\begin{equation}
\begin{aligned}
\Tilde{\mathbf{\Sigma}} &= \mathrm{Cov}\!\left[\psi(\mathbf{w}) \mid \widehat{D}\right]
\\
&=\mathcal{K}(\mathbf{w}, \mathbf{w})
- \mathcal{K}(\mathbf{w},\widehat{w})
\Sigma^{-1}
\mathcal{K}(\widehat{w}, \mathbf{w}),
\end{aligned}
\end{equation}
where $\Sigma = \mathcal{K}(\mathbf{w},\widehat{w})
\left(\mathcal{K}(\widehat{w},\widehat{w}) + \sigma_d^{2}I\right)$. Therefore, the joint capacities $\psi(\mathbf{w})$ across the $T$ time periods follow a multivariate Gaussian distribution with the above posterior mean and covariance. 
We generate the estimated distribution $\widehat{\mathbb{P}}$ by sampling from this posterior distribution. 

\subsection{Computational performance}
\subsubsection{Benchmark settings}
We compare the computational performance of the solution algorithms
introduced in Section~\ref{sec:solution_approach}. DCG-MILP and
DCG-DP solve dr-SAGHP using MILP and
DP separation algorithms, respectively, whereas K-IL solves rs-dr-SAGHP. The comparison between
DCG-MILP and DCG-DP therefore evaluates the computational benefit
of DP separation for the same model, while the comparison with K-IL
also reflects the effect of restricting the adversarial distribution's
support.


Using the estimated joint capacity distributions constructed in
Section~\ref{sec:ScenarioGeneration}, we consider sample sizes
of 100, 500, and 1\,000. For dr-SAGHP, we set the Wasserstein
radius to $\varepsilon=\delta(T+1)K_{\mathrm{ref}}$, where
$K_{\mathrm{ref}}$ is the median historical capacity in the
training data. We select $\delta\in\{0.01,0.05,0.10\}$,
representing average equivalent capacity perturbations per
period of 1\%, 5\%, and 10\% of $K_{\mathrm{ref}}$,
respectively. With $(T+1)K_{\mathrm{ref}}=547.2$ from our dataset, these settings
yield radii of 5.47, 27.36, and 54.72, rounded to two decimal
places. These three settings are used for computational
benchmarking; the separate parameter grids used for operational
calibration are described in Section~\ref{subs:exp-hist}.

To evaluate computational performance
on a relatively demanding instance, we use the busiest day in the
test set, October 2, 2023, which contains 679 flights. All experiments
are conducted on an x86\_64 Linux server equipped with 36 logical
CPUs and 186.9~GB of memory. The models are solved using Gurobi
11.0.3 with its default automatic thread setting
(\texttt{Threads=0}).

We implement GPR using \texttt{GaussianProcessRegressor} from
Python's \texttt{scikit-learn} package. The covariance
kernel combines a constant-scaled radial basis function (RBF)
kernel with a white-noise kernel. The constant-kernel value and
RBF length scale are initialized at $1.0$ and $10.0$, with bounds
$[10^{-3},10^3]$ and $[10^{-2},10^3]$, respectively; the white-noise level is initialized at $1.0$. These are initialization values, and the kernel hyperparameters are optimized during fitting using five optimizer restarts.

\subsubsection{Runtime and solution quality}
Table~\ref{tab:computational_comparison} reports the results.
Among the two algorithms for dr-SAGHP, DCG-DP consistently
requires less runtime than DCG-MILP. For sample sizes of 100 and
500, its speedup relative to DCG-MILP ranges from approximately
$3.73\times$ to $6.24\times$; for 1\,000 samples, the speedup ranges
from $3.18\times$ to $4.31\times$. Both algorithms solve the same
full-support formulation, so these runtime improvements are obtained
without changing the underlying optimization model. The relative
benefit of DP separation varies with the sample size and radius.

K-IL, which solves the restricted-support approximation,
consistently requires less runtime than either DCG algorithm.
It achieves speedups of approximately $6.7$--$11.8\times$
relative to DCG-MILP and $1.3$--$3.5\times$ relative to DCG-DP,
with the largest speedup over DCG-MILP observed in an instance
with 1\,000 samples.
These computational savings are obtained by solving
rs-dr-SAGHP, which restricts the adversarial distribution to the
support $\widehat{\Xi}$ of the estimated nominal distribution.
Consequently, this comparison reflects differences in both the
formulation and the solution procedure. As we will show in Sections~\ref{subs:exp-hist} and Sections~\ref{sec:SensitivityAnalysis}, despite the restriction on the support of the adversial distribution, the solutions provided by rs-dr-SAGHP can perform similarly those provided by dr-SAGHP in practice. Thus, the K-IL approach can provide similar quality solutions with less computational burden.

The comparison above focuses on computational efficiency and
the reported optimality gaps. We evaluate the operational
performance of the resulting schedules in the following sections
through out-of-sample comparisons using observed historical
capacities and simulated capacity distribution shifts.

\begin{table}[htbp]
    \centering
    \caption{Computational performance comparison under different
    sample sizes and Wasserstein radii. Runtime is reported in
    seconds, and speedup is calculated relative to DCG-MILP.}
    \label{tab:computational_comparison}
    \resizebox{\textwidth}{!}{%
    \begin{tabular}{rlrrrrrr}
        \toprule
        & & \multicolumn{2}{c}{$N=100$}
          & \multicolumn{2}{c}{$N=500$}
          & \multicolumn{2}{c}{$N=1\,000$} \\
        \cmidrule(lr){3-4}
        \cmidrule(lr){5-6}
        \cmidrule(lr){7-8}
        $\varepsilon$
        & Method
        & Runtime (s) & Speedup
        & Runtime (s) & Speedup
        & Runtime (s) & Speedup \\
        \midrule

        5.47 & K-IL
        & 215.87  & $6.99\times$
        & 699.36  & $9.55\times$
        & 1521.45 & $11.15\times$ \\

        5.47 & DCG-DP
        & 345.30  & $4.37\times$
        & 1791.27 & $3.73\times$
        & 5337.46 & $3.18\times$ \\

        5.47 & DCG-MILP
        & 1509.11  & \textemdash
        & 6676.61  & \textemdash
        & 16969.80 & \textemdash \\
        \midrule

        27.36 & K-IL
        & 238.28  & $6.69\times$
        & 692.18  & $11.33\times$
        & 1748.64 & $9.90\times$ \\

        27.36 & DCG-DP
        & 326.46  & $4.89\times$
        & 2016.13 & $3.89\times$
        & 4018.26 & $4.31\times$ \\

        27.36 & DCG-MILP
        & 1595.02  & \textemdash
        & 7841.79  & \textemdash
        & 17312.52 & \textemdash \\
        \midrule

        54.72 & K-IL
        & 184.46  & $7.52\times$
        & 636.68  & $10.73\times$
        & 1216.94 & $11.80\times$ \\

        54.72 & DCG-DP
        & 232.22  & $5.97\times$
        & 1095.38 & $6.24\times$
        & 3948.82 & $3.64\times$ \\

        54.72 & DCG-MILP
        & 1386.47  & \textemdash
        & 6831.95  & \textemdash
        & 14355.46 & \textemdash \\

        \bottomrule
    \end{tabular}%
    }
\end{table}

\subsection{Case study using historical data}\label{subs:exp-hist}
In this subsection, we evaluate dr-SAGHP using historical flight
schedules, weather information, and observed airport capacities.
We first describe the calibration of the Wasserstein ambiguity set
radius and then report out-of-sample performance on held-out
historical test days.

To construct the estimated airport-capacity distributions over the
planning horizon, we use historical capacity observations from
the FAA ASPM capacity dataset~\citep{faa_airport_capacity_profiles}
and weather information from the FAA ASPM weather factors
dataset~\citep{faa_aspm_weather_factors_details}. We fit a Gaussian
Process Regression (GPR) model to capacity observations from
January 1, 2021, to December 31, 2022, using the corresponding
weather factors as predictors.
\subsubsection{Calibration of the Wasserstein ambiguity set radius}
\label{sec:radius_calibration}

The Wasserstein radius $\varepsilon$ controls the level of robustness of the 
dr-SAGHP. A small radius produces a solution close to the sp-SAGHP solution, whereas 
a large radius provides additional protection against distributional ambiguity 
but may lead to an overly conservative schedule. Therefore, we calibrate 
$\varepsilon$ based on out-of-sample operational performance using only the 
training data.

\paragraph{Training and testing sets} We partition the historical data into two non-overlapping sets:
\begin{equation}
\begin{aligned}
 \mathcal{D}^{\mathrm{train}}
    &=
    \left\{
    (w_d,\xi_d^{\mathrm{true}})
    \right\}_{d=1}^{D_{\mathrm{train}}}, \\
    \mathcal{D}^{\mathrm{test}}
    &=
    \left\{
    (w_d,\xi_d^{\mathrm{true}})
    \right\}_{d=1}^{D_{\mathrm{test}}},
\end{aligned}
\end{equation}
where $w_d$ denotes the weather information for day $d$, and 
$\xi_d^{\mathrm{true}}$ denotes the corresponding observed airport-capacity 
trajectory. The training set is used for both GPR fitting and radius 
calibration, whereas the test set is held out until the final out-of-sample evaluation. Note that no information from $\mathcal{D}^{\mathrm{test}}$ is used during GPR training 
or radius calibration. This prevents test information from influencing the 
selection of $\varepsilon$.

\paragraph{Expanding-Window cross validation} 
We use expanding-window cross-validation to select the Wasserstein
radius that minimizes the average realized operational cost on
validation days. At each split, we fit the GPR model using earlier
observations and solve dr-SAGHP for each candidate radius using
the predicted capacity distributions for the next validation block.
We evaluate the resulting schedules against the observed capacities
in that block. Repeating this procedure with an expanding training
window allows us to compare candidate radii across multiple
validation blocks. The selected radius is fixed before evaluating
performance on the held-out test set.

To preserve the temporal ordering of the data, we partition 
$\mathcal{D}^{\mathrm{train}}$ into $N^{\mathcal{B}}$ consecutive, non-overlapping blocks:
\begin{equation}
    \mathcal{D}^{\mathrm{train}}
    =
    \mathcal{B}_1
    \cup
    \mathcal{B}_2
    \cup
    \cdots
    \cup
    \mathcal{B}_{N^{\mathcal{B}}}.
\end{equation}
We use an expanding-window validation procedure. For validation split 
$k=1,\ldots,N^{\mathcal{B}}-1$, the GPR training and validation sets are defined as
\begin{equation}
    \mathcal{D}^{\mathrm{fit}}_k
    =
    \bigcup_{j=1}^{k}\mathcal{B}_j,
    \qquad
    \mathcal{D}^{\mathrm{val}}_k
    =
    \mathcal{B}_{k+1}.
    \label{eq:radius_expanding_window}
\end{equation}
Thus, each GPR model is trained only on observations that occur before its 
validation block. In our implementation, we divide the training period into five consecutive blocks, yielding four expanding-window validation splits.

For each split $k$, we fit a GPR model using 
$\mathcal{D}^{\mathrm{fit}}_k$. For each validation day 
$d\in\mathcal{B}_{k+1}$, we then draw $N^{\mathrm{train}}$ capacity trajectories from the 
corresponding GPR predictive distribution:
\begin{equation}
    \widehat{\xi}^{(i,k)}_d
    \sim
    \mathbb{P}^{(k)}_{\mathrm{GPR}}
    \left(
    \xi\mid w_d
    \right),
    \qquad
    i=1,\ldots,N^{\mathrm{train}}.
\end{equation}
We approximate the GPR predictive capacity distribution using
the sampled trajectories:
\begin{equation}
    \widehat{\mathbb{P}}^{(k)}_{d,N^{\mathrm{train}}}
    =
    \frac{1}{N^{\mathrm{train}}}
    \sum_{i=1}^{N^{\mathrm{train}}}
    \delta_{\widehat{\xi}^{(i,k)}_d}.
\end{equation}
This discrete estimate serves as the nominal capacity distribution
supplied to dr-SAGHP and forms the center of its Wasserstein
ambiguity set.

Let $\mathcal{R}$ denote the set of candidate Wasserstein radii. Each candidate radius $\varepsilon \in \mathcal{R}$ is evaluated as follows. For each validation day $d\in\mathcal{B}_{k+1}$, we solve 
the dr-SAGHP centered at $\widehat{\mathbb{P}}^{(k)}_{d,N^{\mathrm{train}}}$ to obtain $x_{d,k}^{*}(\varepsilon)$.
For a fixed validation day, the same GPR-generated capacity trajectories are 
used for all candidate radii to ensure that their performance differences are 
not caused by sampling variation.

We then evaluate the resulting schedule against the true observed capacity 
trajectory $\xi_d^{\mathrm{true}}$. Let $G(x)$ denote the scheduled ground-delay, then the realized operational cost on day $d$ for the $k$th validation set when using radius $\varepsilon$ is
\begin{equation}
    C_{d,k}(\varepsilon)
    =
    G\left(x_{d,k}^{*}(\varepsilon)\right)
    +
    Q\left(
    x_{d,k}^{*}(\varepsilon),
    \xi_d^{\mathrm{true}}
    \right).
    \label{eq:radius_validation_cost}
\end{equation}
Using the true capacity trajectory, rather than an additional GPR-generated 
trajectory, allows the calibration procedure to account for the GPR modeling error.
Let $\widehat{C}_{\mathrm{val}}(\varepsilon)$ denote the average
realized operational cost obtained using radius $\varepsilon$,
averaged over all days across the validation blocks:
\begin{equation}
    \widehat{C}_{\mathrm{val}}(\varepsilon)
    =
    \frac{1}{D_{\mathrm{val}}}
    \sum_{k=1}^{N^{\mathcal{B}}-1}
    \sum_{d\in\mathcal{B}_{k+1}}
    C_{d,k}(\varepsilon),
    \label{eq:radius_validation_performance}
\end{equation}
where
\begin{equation}
    D_{\mathrm{val}}
    =
    \sum_{k=1}^{N^{\mathcal{B}}-1}
    \left|
    \mathcal{B}_{k+1}
    \right|.
\end{equation}
We select the radius that minimizes the average realized operational cost:
\begin{equation}
    \varepsilon^{*}
    \in
    \arg\min_{\varepsilon\in\mathcal{R}}
    \widehat{C}_{\mathrm{val}}(\varepsilon).
    \label{eq:calibrated_wasserstein_radius}
\end{equation}
The selected radius $\varepsilon^{*}$ is fixed before evaluating any schedules on 
the test set. A similar procedure is used to calibrate parameters of other models.

\subsubsection{Out-of-Sample evaluation}

After selecting $\varepsilon^{*}$, we refit the GPR model using the complete 
training set $\mathcal{D}^{\mathrm{train}}$. For each test day 
$d\in\mathcal{D}^{\mathrm{test}}$, we draw $N^{\mathrm{test}}$ capacity trajectories from the 
refitted GPR model conditioned on $w_d$ to serve as the estimated probability distribution and solve the dr-SAGHP using the fixed 
radius $\varepsilon^{*}$. The resulting schedule $x_d^{*}(\varepsilon^{*})$ is then evaluated against the 
true test capacity trajectory:
\begin{equation}
    C_d^{\mathrm{test}}
    =
    G\left(x_d^{*}(\varepsilon^{*})\right)
    +
    Q\left(
    x_d^{*}(\varepsilon^{*}),
    \xi_d^{\mathrm{true}}
    \right).
\end{equation}
The average out-of-sample cost is
\begin{equation}
    \widehat{C}_{\mathrm{test}}
    =
    \frac{1}{D_{\mathrm{test}}}
    \sum_{d\in\mathcal{D}^{\mathrm{test}}}
    C_d^{\mathrm{test}}.
\end{equation}
We calibrate the radii for dr-SAGHP and rs-dr-SAGHP, and the robust benchmark's
quantile level independently using expanding-window validation.

For dr-SAGHP, the transportation cost is measured by the
$\ell_1$ distance between capacity trajectories. We define
the candidate radii as
\begin{equation}
    \mathcal R_{\mathrm{dr-SAGHP}}
    =
    \left\{
    \delta(T+1)K_{\mathrm{ref}}:
    \delta\in\Delta
    \right\},
    \label{eq:dp_candidate_radii}
\end{equation}
where $\Delta=\{0.01,0.02,\ldots,0.10\}$ and
$K_{\mathrm{ref}}$ is the median historical capacity in the
training data. Here, $\delta$ represents the average equivalent
capacity perturbation per time period, expressed as a fraction
of $K_{\mathrm{ref}}$.

For rs-dr-SAGHP, the adversarial distribution is restricted
to the support $\widehat{\Xi}$ of the estimated nominal
distribution. We use normalized transportation distances:
\begin{equation}
    d_{i,j}^{\mathrm{emp}}
    =
    \frac{
    \left\|
    \widehat\xi^{i}-\widehat{\xi}^{(j)}
    \right\|_1
    }{
    \displaystyle
    \max_{k\ne\ell}
    \left\|
    \widehat{\xi}^{(k)}-\widehat{\xi}^{(\ell)}
    \right\|_1
    },
    \qquad i,j=1,\ldots,N.
    \label{eq:kil_normalized_distance}
\end{equation}
The corresponding Wasserstein radius is dimensionless, and
its candidate set is
$\mathcal R_{\mathrm{rs-dr-SAGHP}}=\{0.1,0.2,\ldots,1.0\}$.

For r-SAGHP, we select the quantile level $q$ from
$\mathcal Q=\{0.1,0.2,0.3,0.4\}$.
For each model, we choose the parameter that minimizes the
average realized operational cost across all validation days.
The selected parameters are fixed before evaluating performance
on the held-out test set.
\subsubsection{Out-of-sample performance}
With the parameters selected by the validation procedure fixed,
we evaluate the four models on the 84 held-out severe-weather
test days using their observed capacity trajectories.
Table~\ref{tab:oos_cost_comparison} reports the resulting
out-of-sample performance. Relative to s-SAGHP, dr-SAGHP and
rs-dr-SAGHP reduce the mean operational cost by $2.60\%$ and
$2.56\%$, respectively, while also achieving substantial reductions
in $\mathrm{CVaR}_{90}$. The robust benchmark r-SAGHP also improves
tail performance, but its reduction in $\mathrm{CVaR}_{90}$ is
smaller than those of the two distributionally robust models,
and its mean operational cost increases by $0.31\%$ relative to
s-SAGHP. These results indicate that, on the historical test set,
the distributionally robust models provide stronger tail protection
than r-SAGHP while also improving average operational performance.

dr-SAGHP and rs-dr-SAGHP exhibit similar out-of-sample performance
after their Wasserstein radii are calibrated independently.
Although rs-dr-SAGHP restricts the adversarial distribution to
the support of the estimated nominal distribution, its performance
remains competitive with that of the full-support dr-SAGHP in
both mean operational cost and conditional value-at-risk at the
90\% confidence level ($\mathrm{CVaR}_{90}$), which measures the
average operational cost in the worst 10\% of outcomes. This finding suggests that, with appropriate
radius calibration, rs-dr-SAGHP as an restricted-support approximation of dr-SAGHP can
retain much of the operational benefit of dr-SAGHP in this case study.
\begin{figure}[htbp]
    \centering
    \includegraphics[width=0.70\textwidth]{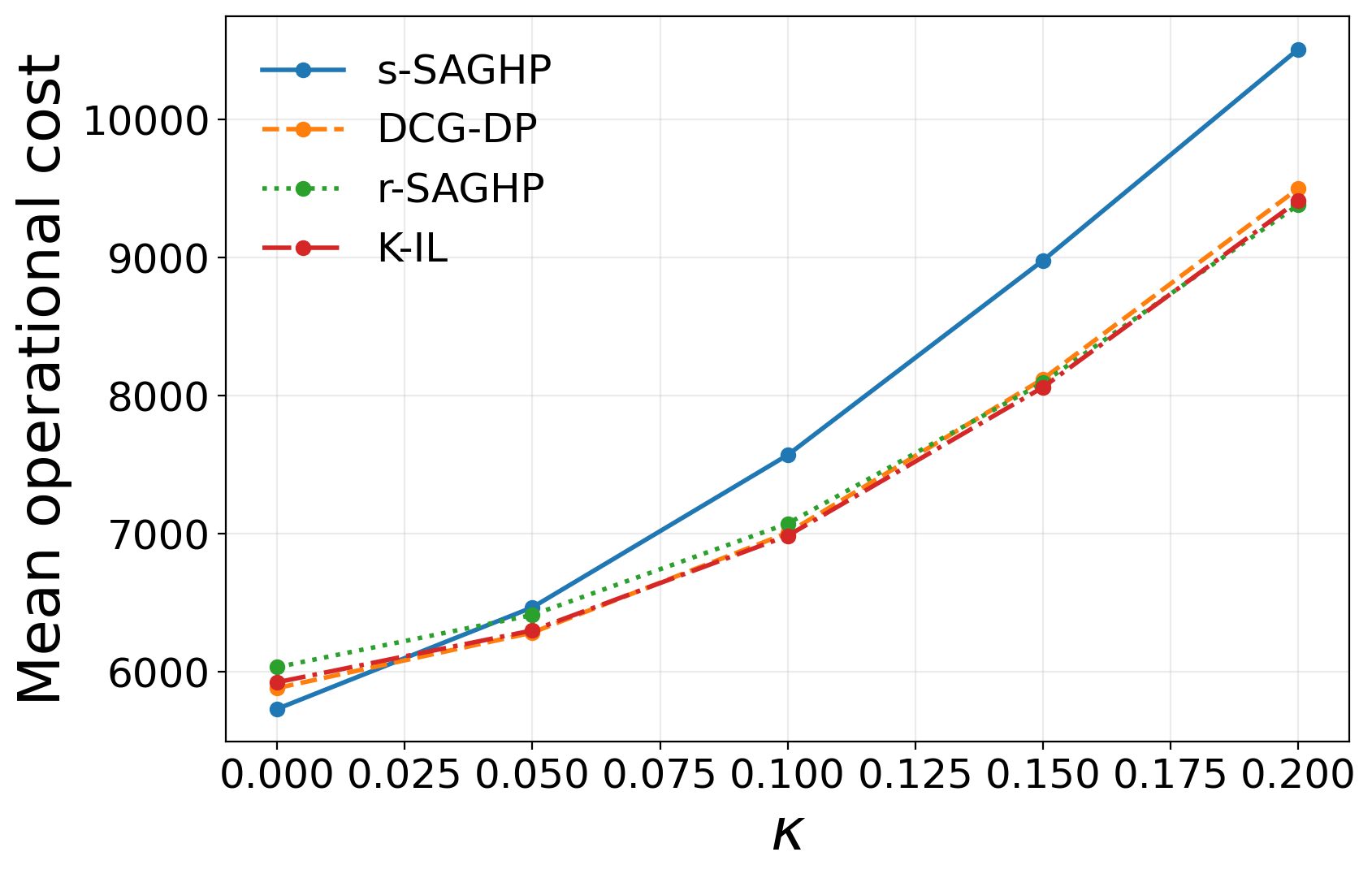}
    \caption{Out-of-sample performance comparison in mean operational cost under different levels of mean reduction.}
    \label{fig:oos_mean_cost_mean_reduction}
\end{figure}

\begin{figure}[htbp]
    \centering
    \includegraphics[width=0.70\textwidth]{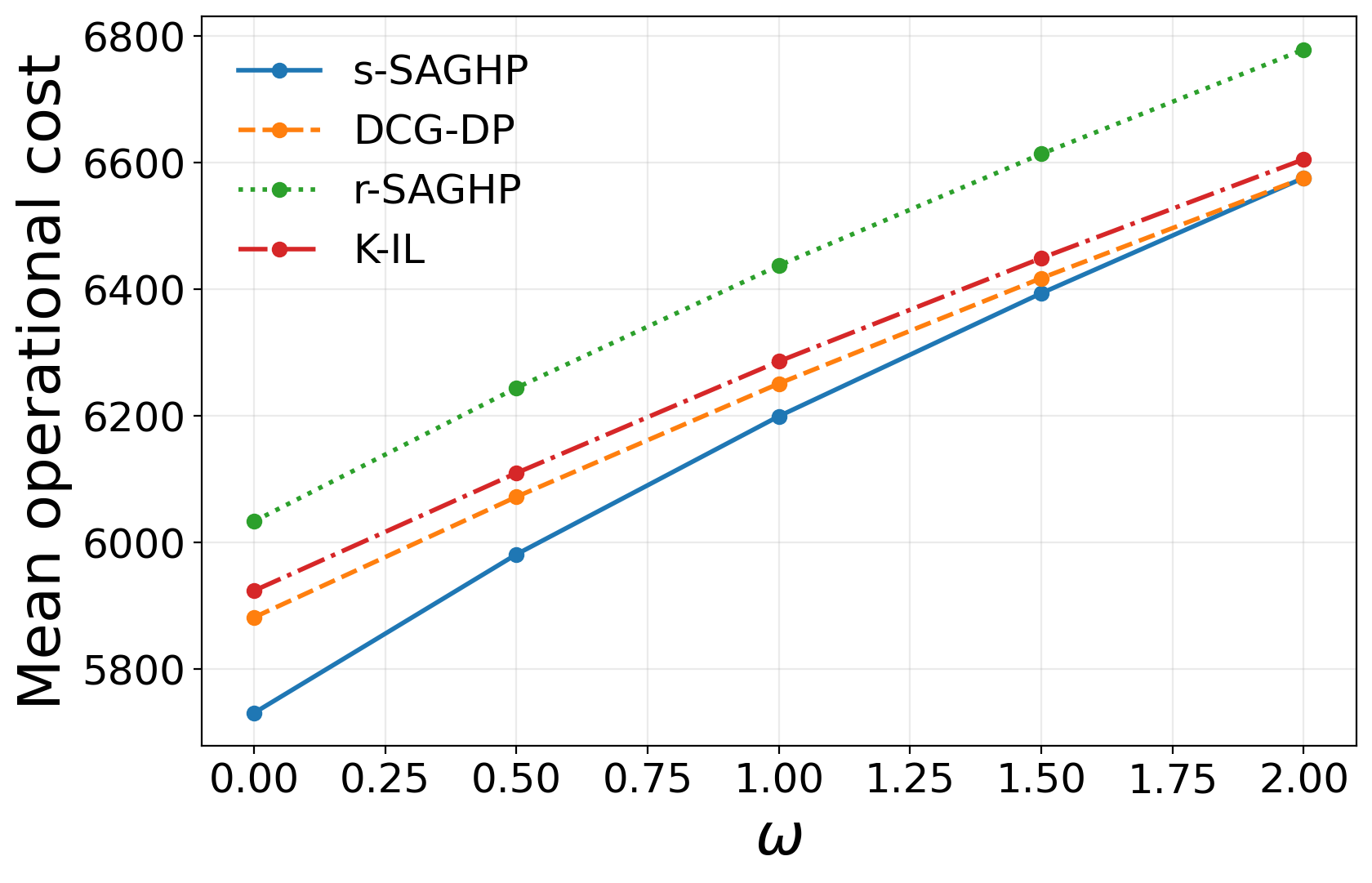}
    \caption{Out-of-sample performance comparison in mean operational cost under different levels of covariance increase.}
    \label{fig:oos_mean_cost_covariance_increase}
\end{figure}

\begin{table}[htbp]
    \centering
    \caption{Comparison of out-of-sample mean and CVaR operational costs. Percentage reductions are computed relative to s-SAGHP.}
    \label{tab:oos_cost_comparison}
    \begin{tabular}{lcrrrr}
        \toprule
        \textbf{Method}
        & \textbf{Calibrated Parameter}
        & \textbf{Mean Cost}
        & \textbf{Reduction (\%)}
        & \textbf{CVaR$_{90}$}
        & \textbf{Reduction (\%)} \\
        \midrule

        s-SAGHP
        & --
        & 7572.42
        & --
        & 14249.00
        & -- \\

        dr-SAGHP
        & $\varepsilon_{\mathrm{DP}}^{*}=27.36$
        & \textbf{7375.24}
        & \textbf{2.60}
        & 13137.44
        & 7.80 \\

        r-SAGHP
        & $\alpha^{*}=0.30$
        & 7595.83
        & $-0.31$
        & 13357.11
        & 6.26 \\

        rs-dr-SAGHP
        & $\varepsilon_{\mathrm{K\text{-}IL}}^{*}=0.10$
        & 7378.68
        & 2.56
        & \textbf{13096.22}
        & \textbf{8.09} \\

        \bottomrule
    \end{tabular}
\end{table}

\subsection{Sensitivity analysis on climate change impacts}\label{sec:SensitivityAnalysis}
Differences between the estimated and out-of-sample capacity distributions in the historical data may reflect changes in weather patterns, potentially associated with climate change, as well as estimation errors and other operational factors. To isolate the effects of weather-related capacity distribution shifts in a controlled setting, we conduct a sensitivity analysis
by perturbing the mean and covariance of the GPR predictive distributions. Specifically, we consider two channels through which climate change may affect the capacity distribution: (i) a decrease in its mean and (ii) an increase in its variance. The mean shift represents persistent capacity degradation associated with more frequent or more severe adverse weather conditions~\citep{bhatia2022potential}. The increase in variance represents greater uncertainty in capacity predictions under more variable severe weather conditions~\citep{mondal2016hydrologic}.

\subsubsection{Design of the sensitivity analysis}
For each historical test day $d$, let $w_d$ denote the observed ASPM weather features over the planning horizon. We condition the fitted GPR model on $w_d$ to obtain the joint predictive capacity distribution and sample capacity trajectories from this distribution to construct the estimated nominal distribution. We note that the procedure described here is similar to the experiments using historical data in Section~\ref{subs:exp-hist}. We then keep the weather inputs fixed and perturb the predictive mean or covariance to generate the testing distributions used in the sensitivity analysis.
To generate testing distributions for mean shifts, we impose mean reductions of $\kappa \in \{0.05, 0.10, 0.15, 0.20\}$ on the fitted GPR posterior mean and draw 1000 capacity scenarios from the resulting shifted posterior. The out-of-sample cost is evaluated by solving a stochastic program with the same structure as \eqref{eq:s-SAGHP_first_stage}, using the ground-holding policy $x^{*}$ obtained from s-SAGHP, r-SAGHP, dr-SAGHP and rs-dr-SAGHP, and using the shifted capacity samples as input parameters.

To assess the impact of climate change on capacity variability, we increase the posterior covariance of the GPR model by scaling it as
\[
\Tilde{\mathbf{\Sigma}}_{\omega}=(1+\omega)\Tilde{\mathbf{\Sigma}},\qquad \omega\in\{0.5,1.0,1.5,2.0\},
\]
so that the covariance is multiplied by factors \(1.5, 2.0, 2.5,\) and \(3.0\), respectively. Larger values of $\omega$ induce greater variability in the generated capacity scenarios. For sufficiently large $\omega$, however, some sampled capacities may exceed the physical maximum capacity of DEN. To ensure physical feasibility, we cap each sampled capacity at the airport's maximum operational capacity.
When evaluating the out-of-sample performance, in addition to evaluating the expected objective value, we also evaluate performance under adverse conditions, using Conditional Value-at-Risk (CVaR). CVaR captures tail risk and thus provides information regarding of operational performance in conditions that are worse than typical. Following~\citet{rockafellar2000optimization}, let $\Phi(x,\mathbf{\psi(\mathbf{w})})$ denote the realized cost induced by decision $x$ under random parameters $\mathbf{\psi(\mathbf{w})} \sim \widehat{\mathbb{P}}(\mathbf{\psi(\mathbf{w})} \mid \widehat D)$ with a posterior density function $p_{\mathbf{\psi(\mathbf{w})}\mid\widehat{D}}(\cdot)$.
For a confidence level $\rho\in(0,1)$, define the value-at-risk (VaR) as
\begin{equation}
\begin{aligned}
\zeta_{\rho}(x)
&:= \mathrm{VaR}_{\rho}\!\bigl(\Phi(x,\mathbf{\psi})\bigr) \\
&= \inf\Biggl\{ z\in\mathbb{R} :
\int_{\Phi(x,\mathbf{\psi})\le z} p_{\mathbf{\psi}\mid \widehat D}(\mathbf \psi)\,d\mathbf \psi \;\ge\; \rho
\Biggr\}.
\end{aligned}
\label{eq:var_def}
\end{equation}
Then the corresponding conditional value-at-risk (CVaR) is
\begin{equation}
\begin{aligned}
&\mathrm{CVaR}_{\rho}\!\bigl(\Phi(x,\mathbf{\psi})\bigr)\\
&=
\frac{1}{1-\rho}
\int_{\Phi(x,\mathbf{\psi})\ge \zeta_{\rho}(x)}
\Phi(x,\mathbf{\psi})\, p_{\mathbf{\psi}\mid \widehat D}(\mathbf{\psi})\, d\mathbf{\psi}.
\label{eq:cvar_def}
\end{aligned}
\end{equation}

\subsubsection{Sensitivity analysis results}
\begin{figure}[htbp]
    \centering

    \begin{subfigure}[t]{0.47\textwidth}
        \centering
        \includegraphics[width=\linewidth]{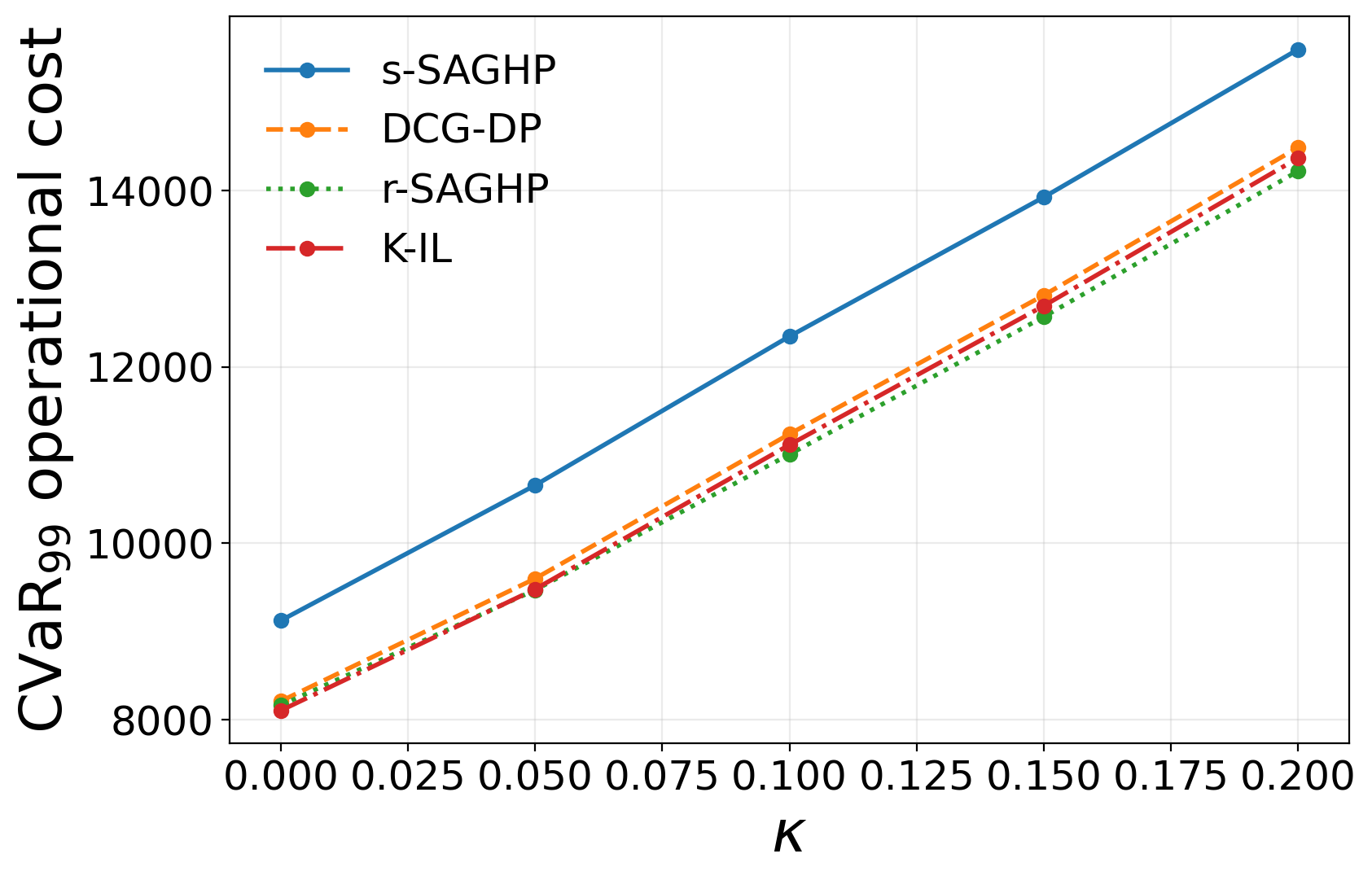}
        \caption{$\rho = 0.01$}
        \label{fig:oos_cvar_rho_001}
    \end{subfigure}
    \hfill
    \begin{subfigure}[t]{0.47\textwidth}
        \centering
        \includegraphics[width=\linewidth]{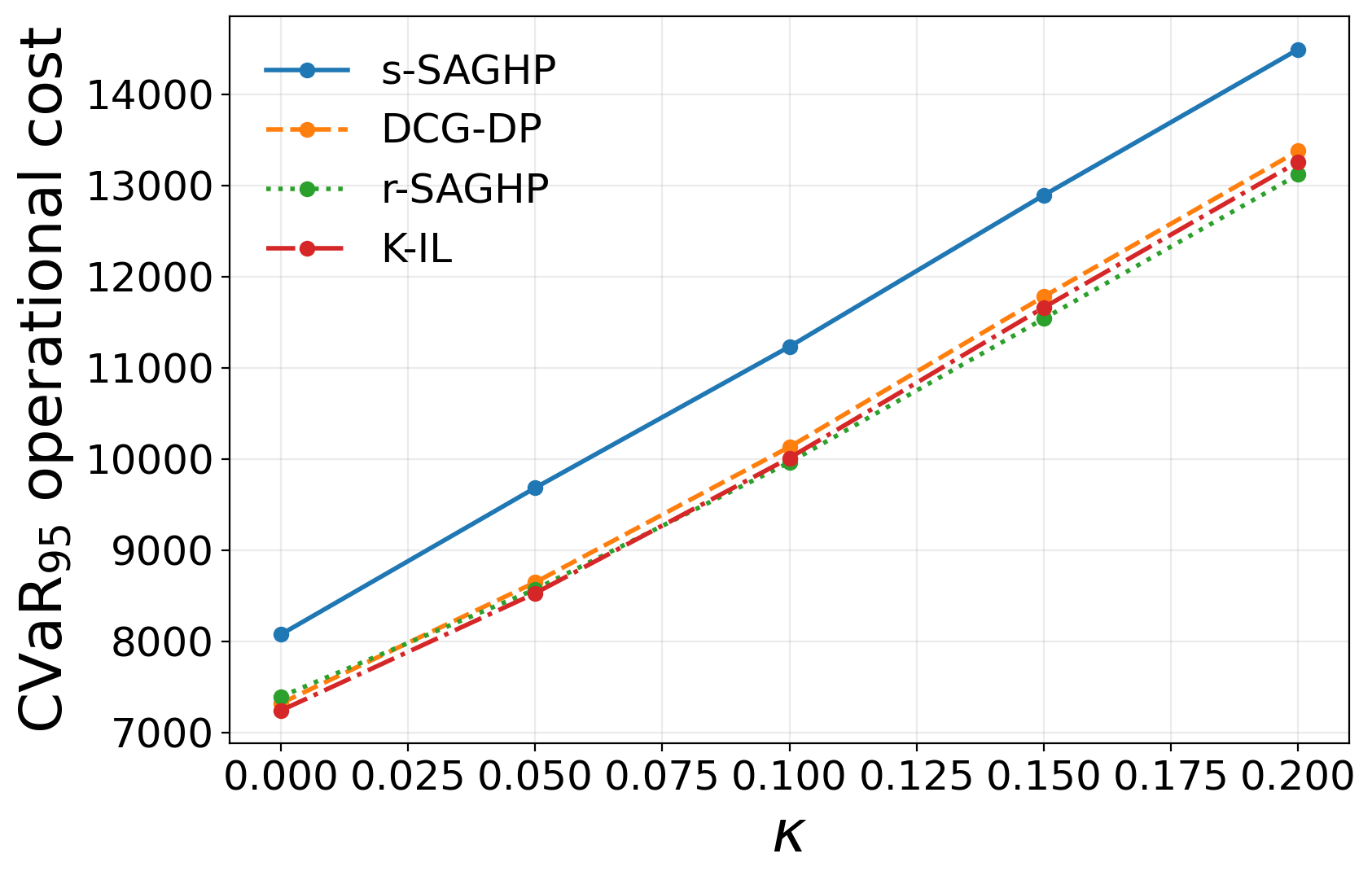}
        \caption{$\rho = 0.05$}
        \label{fig:oos_cvar_rho_005}
    \end{subfigure}

    \vspace{0.5em}

    \begin{subfigure}[t]{0.47\textwidth}
        \centering
        \includegraphics[width=\linewidth]{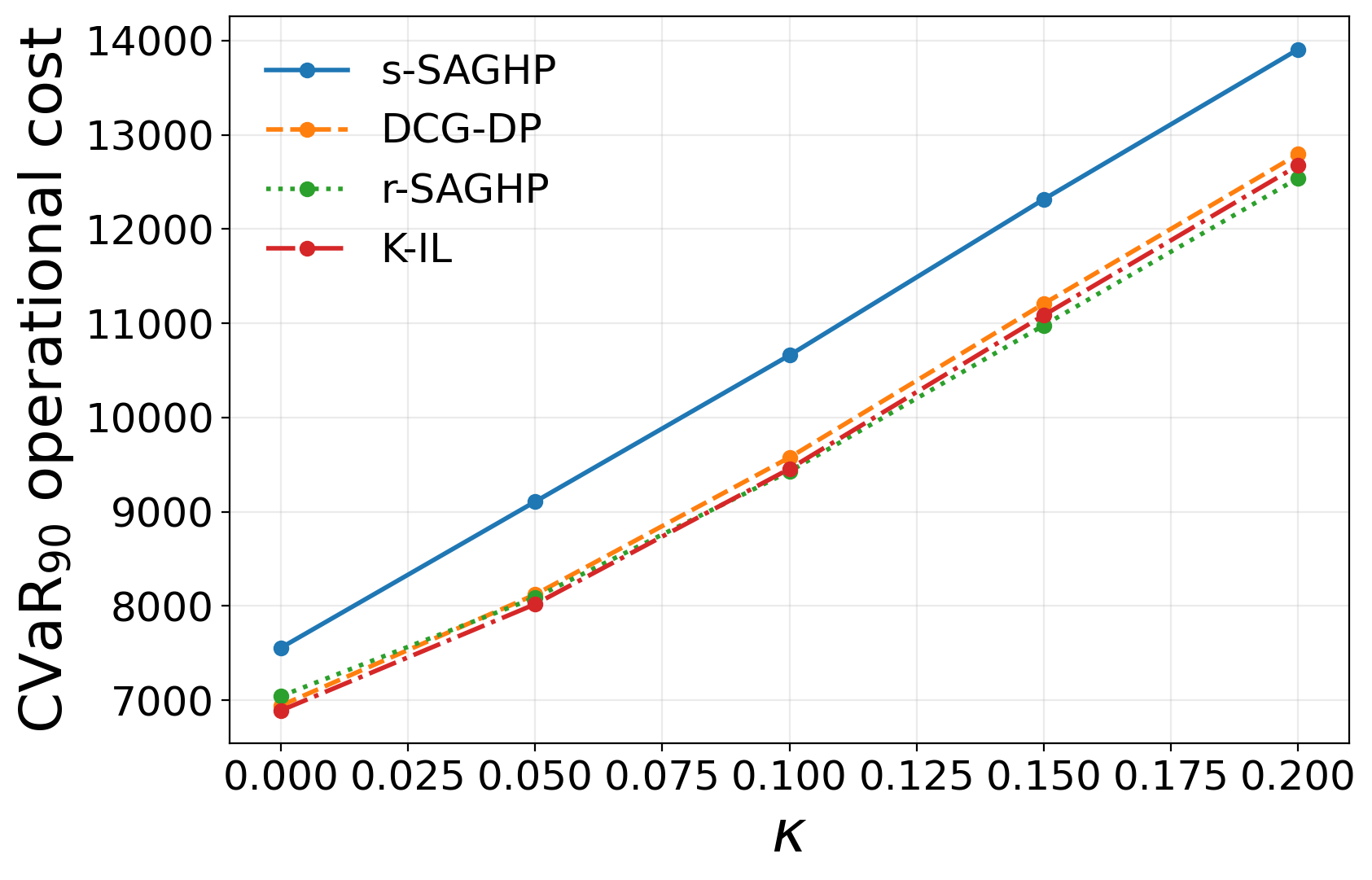}
        \caption{$\rho = 0.10$}
        \label{fig:oos_cvar_rho_010}
    \end{subfigure}
    \hfill
    \begin{subfigure}[t]{0.47\textwidth}
        \centering
        \includegraphics[width=\linewidth]{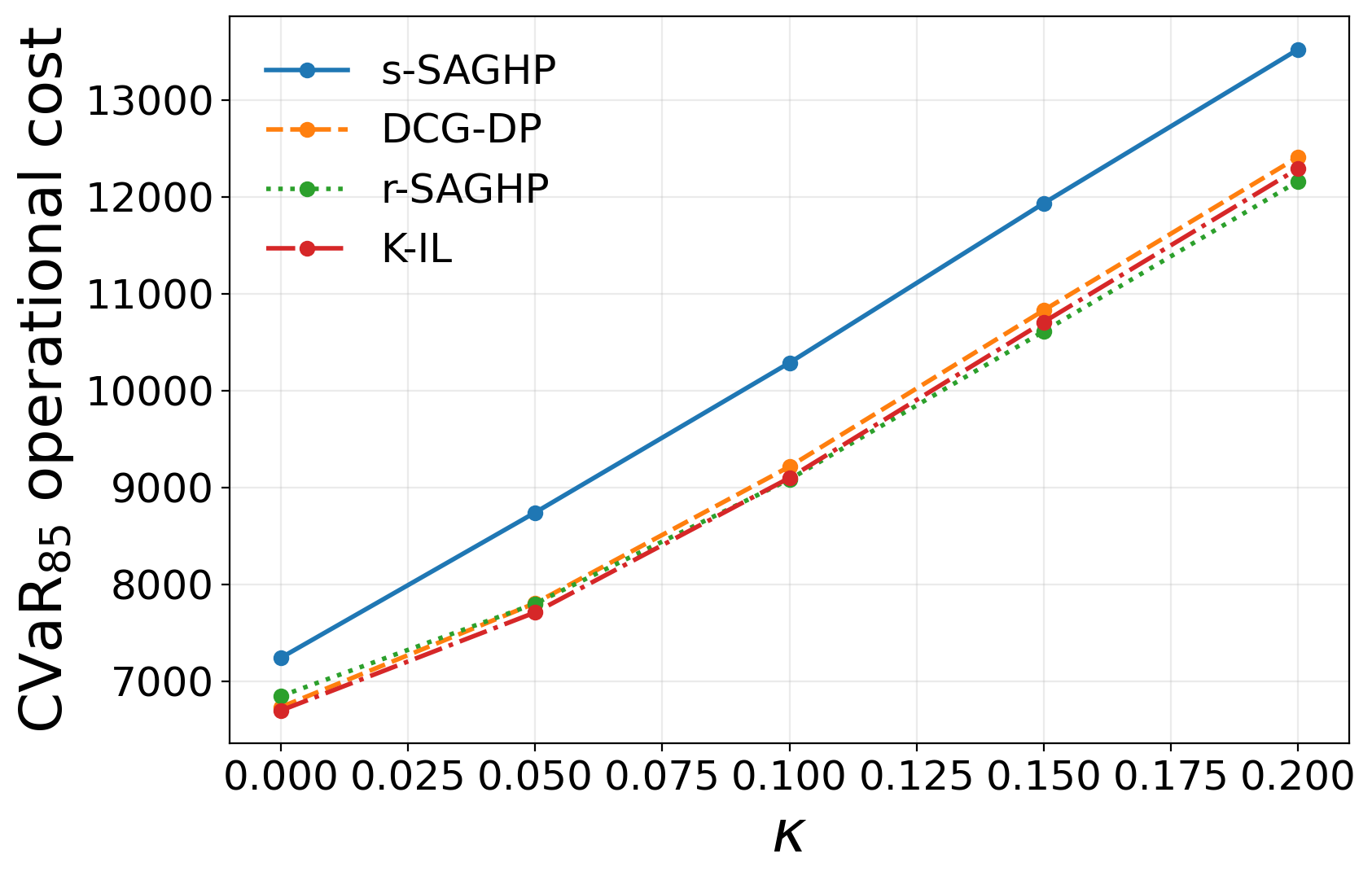}
        \caption{$\rho = 0.15$}
        \label{fig:oos_cvar_rho_015}
    \end{subfigure}

    \vspace{0.5em}

    \begin{subfigure}[t]{0.47\textwidth}
        \centering
        \includegraphics[width=\linewidth]{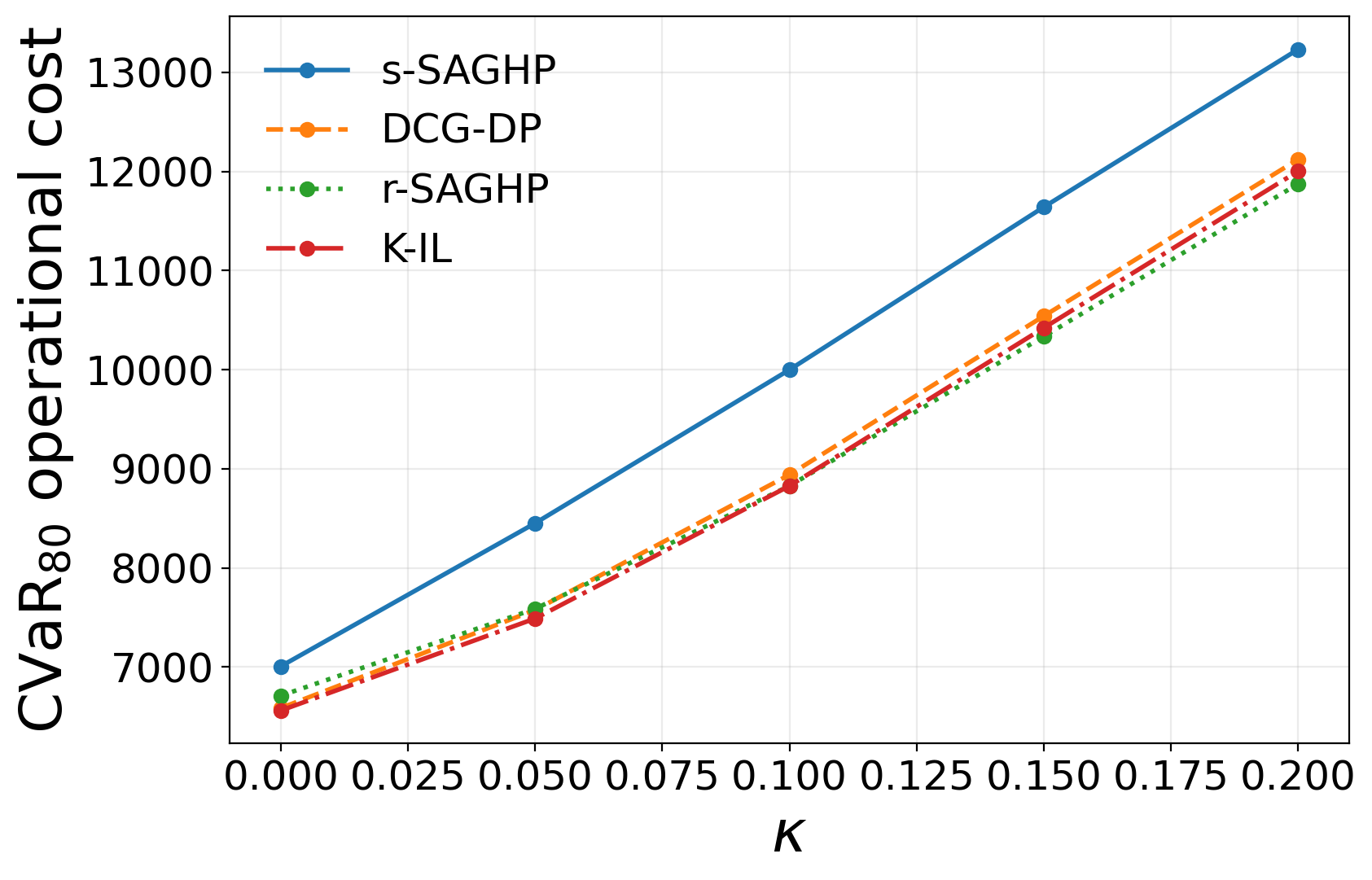}
        \caption{$\rho = 0.20$}
        \label{fig:oos_cvar_rho_020}
    \end{subfigure}

    \caption{Out-of-sample performance comparison in CVaR with five confidence levels under different levels of mean reduction.}
    \label{fig:oos_cvar_comparison_mean_reduction}
\end{figure}

We first consider the mean-reduction experiments shown in Figures~\ref{fig:oos_mean_cost_mean_reduction} and~\ref{fig:oos_cvar_comparison_mean_reduction}. Figure~\ref{fig:oos_mean_cost_mean_reduction} indicates that s-SAGHP achieves the lowest mean operational cost when the mean capacity reduction is small, approximately $\kappa<2.5\%$. Beyond this point, dr-SAGHP, rs-dr-SAGHP, and r-SAGHP all outperform s-SAGHP in terms of mean operational cost, and the performance gap generally widens as $\kappa$ increases. Among these three methods, dr-SAGHP and rs-dr-SAGHP exhibit very similar performance and outperform r-SAGHP up to approximately $\kappa=15\%$. Under the most severe mean-reduction scenario, $\kappa=20\%$, r-SAGHP performs slightly better than the two distributionally robust approaches. As reported in Appendix~\ref{sec:sensitivity_analysis_results}, r-SAGHP reduces the mean operational cost by $10.74\%$ relative to s-SAGHP, compared with reductions of $9.61\%$ and $10.45\%$ for dr-SAGHP and rs-dr-SAGHP, respectively. In terms of tail performance, r-SAGHP provides slightly stronger CVaR protection than  dr-SAGHP and rs-dr-SAGHP once $\kappa$ exceeds $5\%$, and this pattern is broadly consistent across the five confidence levels considered.

\begin{figure}[htbp]
    \centering

    \begin{subfigure}[t]{0.47\textwidth}
        \centering
        \includegraphics[width=\linewidth]{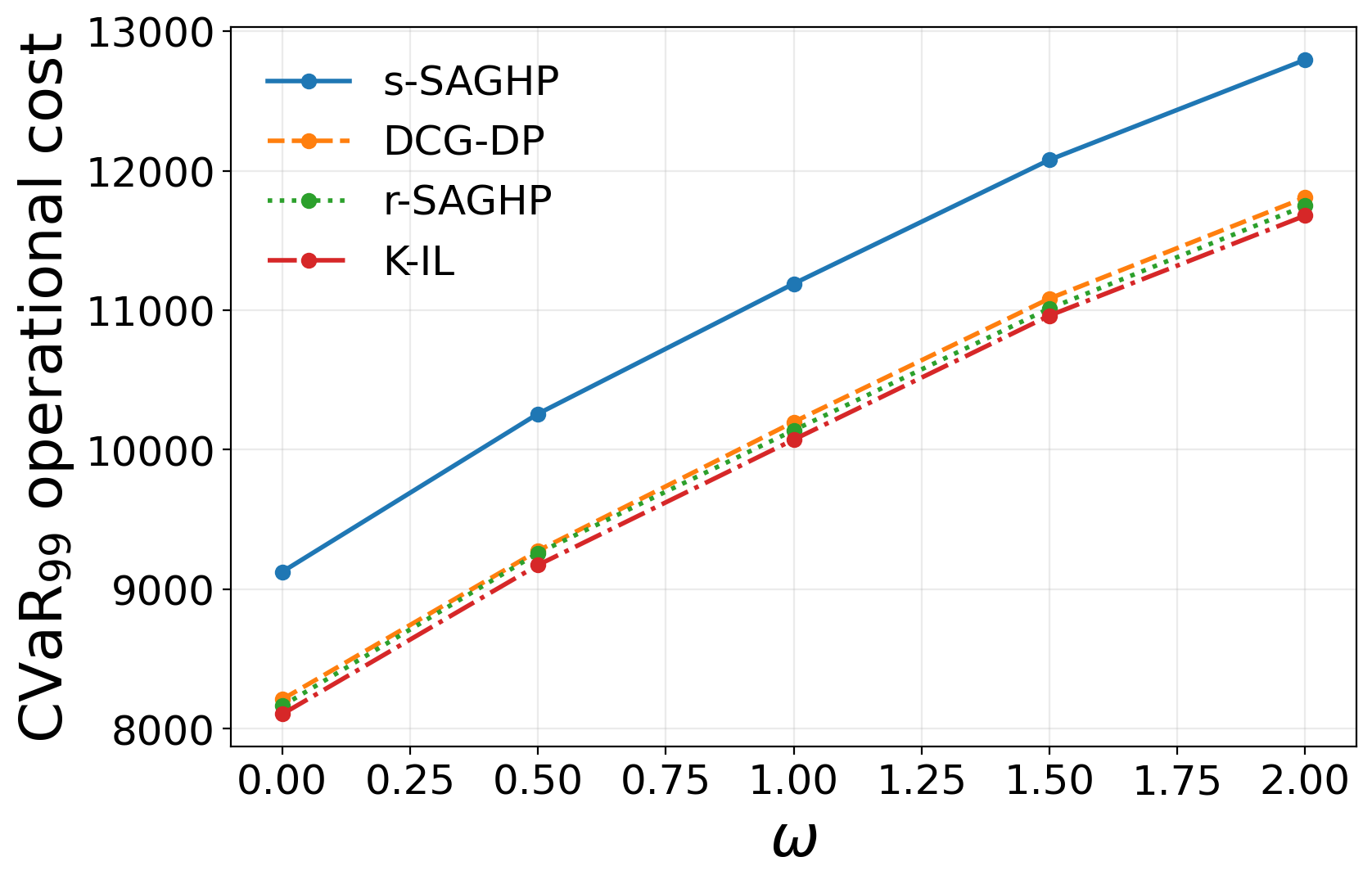}
        \caption{$\rho = 0.01$}
        \label{fig:oos_cvar_cov_rho_001}
    \end{subfigure}
    \hfill
    \begin{subfigure}[t]{0.47\textwidth}
        \centering
        \includegraphics[width=\linewidth]{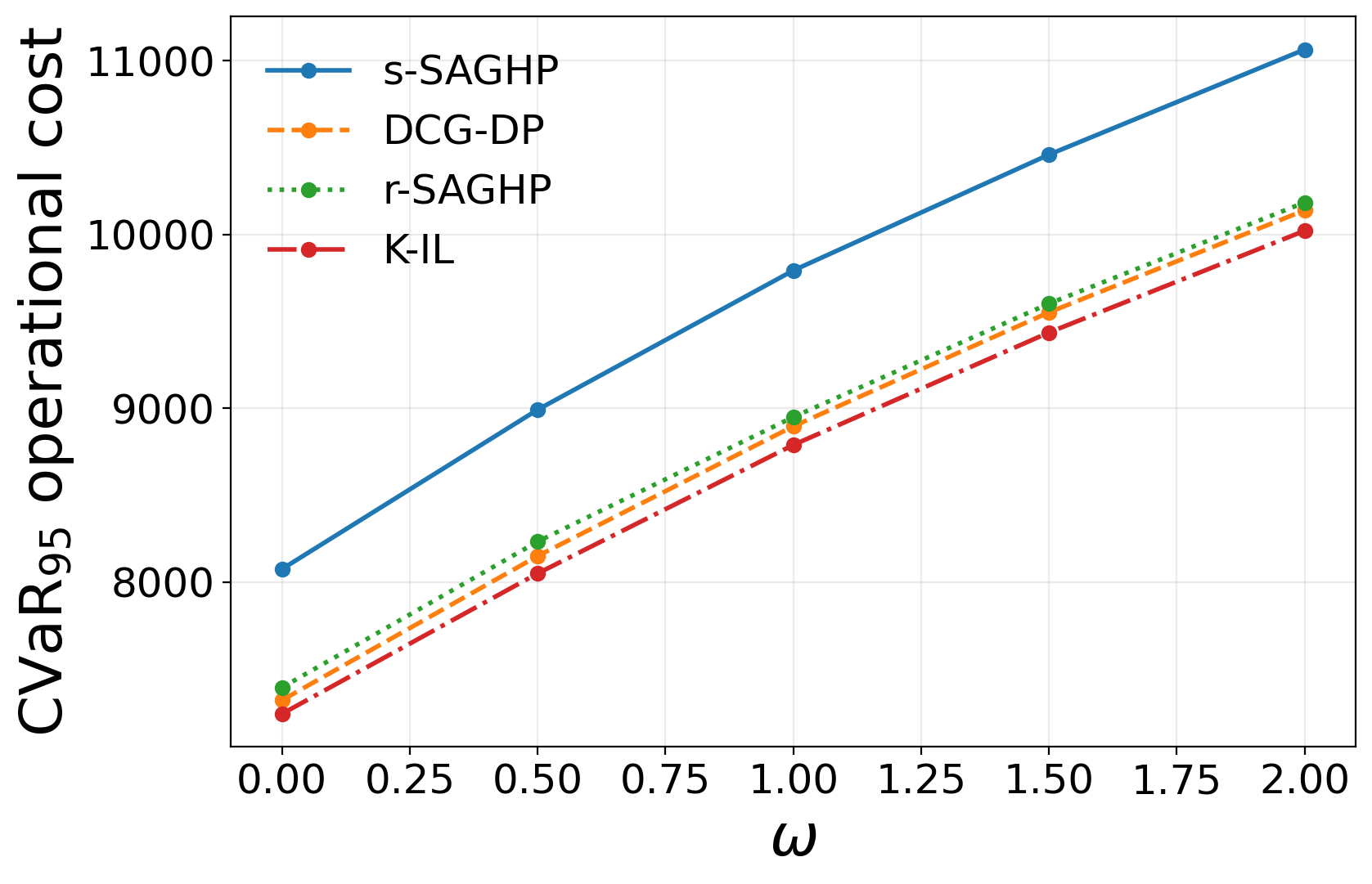}
        \caption{$\rho = 0.05$}
        \label{fig:oos_cvar_cov_rho_005}
    \end{subfigure}

    \vspace{0.5em}

    \begin{subfigure}[t]{0.47\textwidth}
        \centering
        \includegraphics[width=\linewidth]{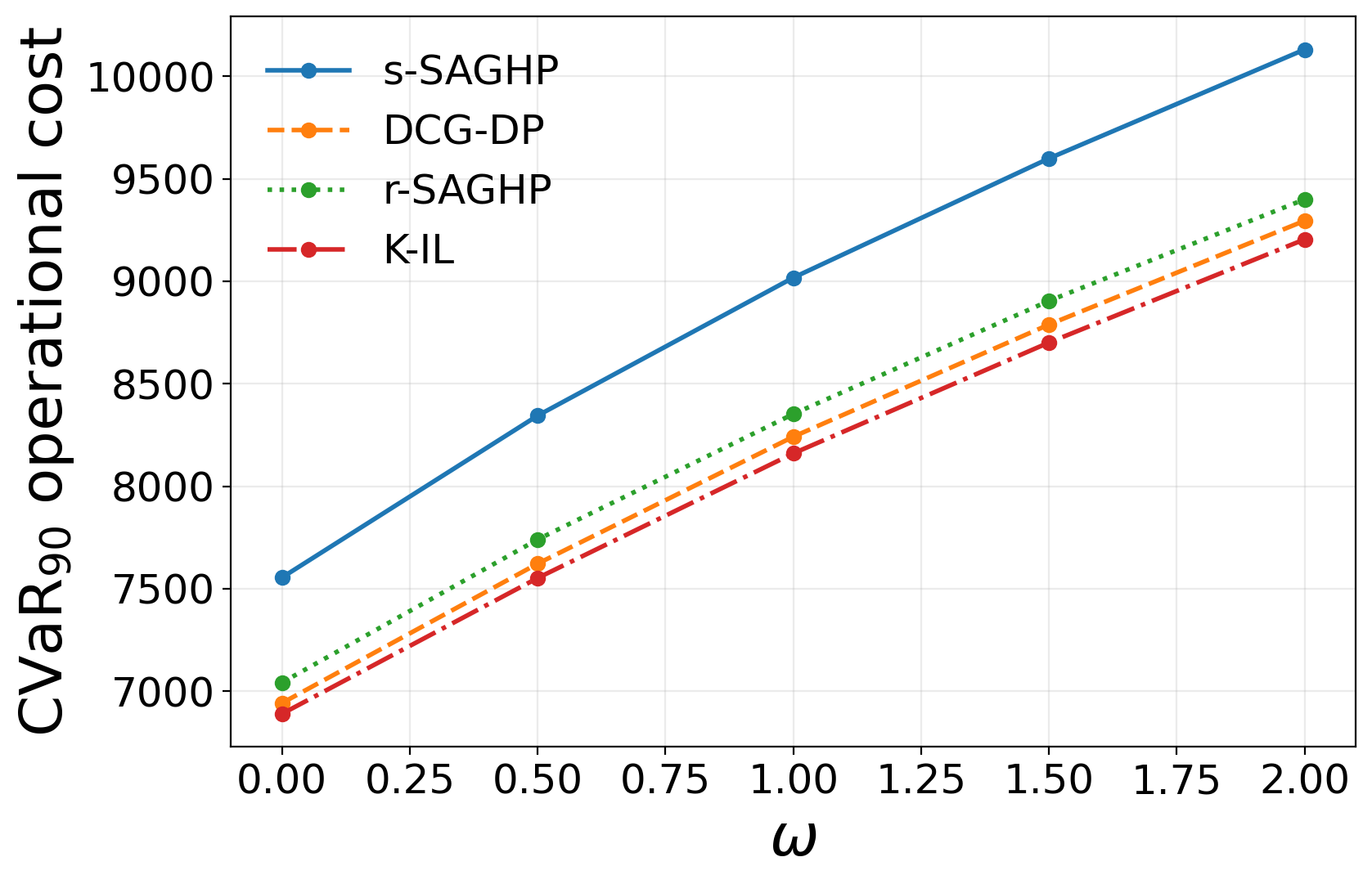}
        \caption{$\rho = 0.10$}
        \label{fig:oos_cvar_cov_rho_010}
    \end{subfigure}
    \hfill
    \begin{subfigure}[t]{0.47\textwidth}
        \centering
        \includegraphics[width=\linewidth]{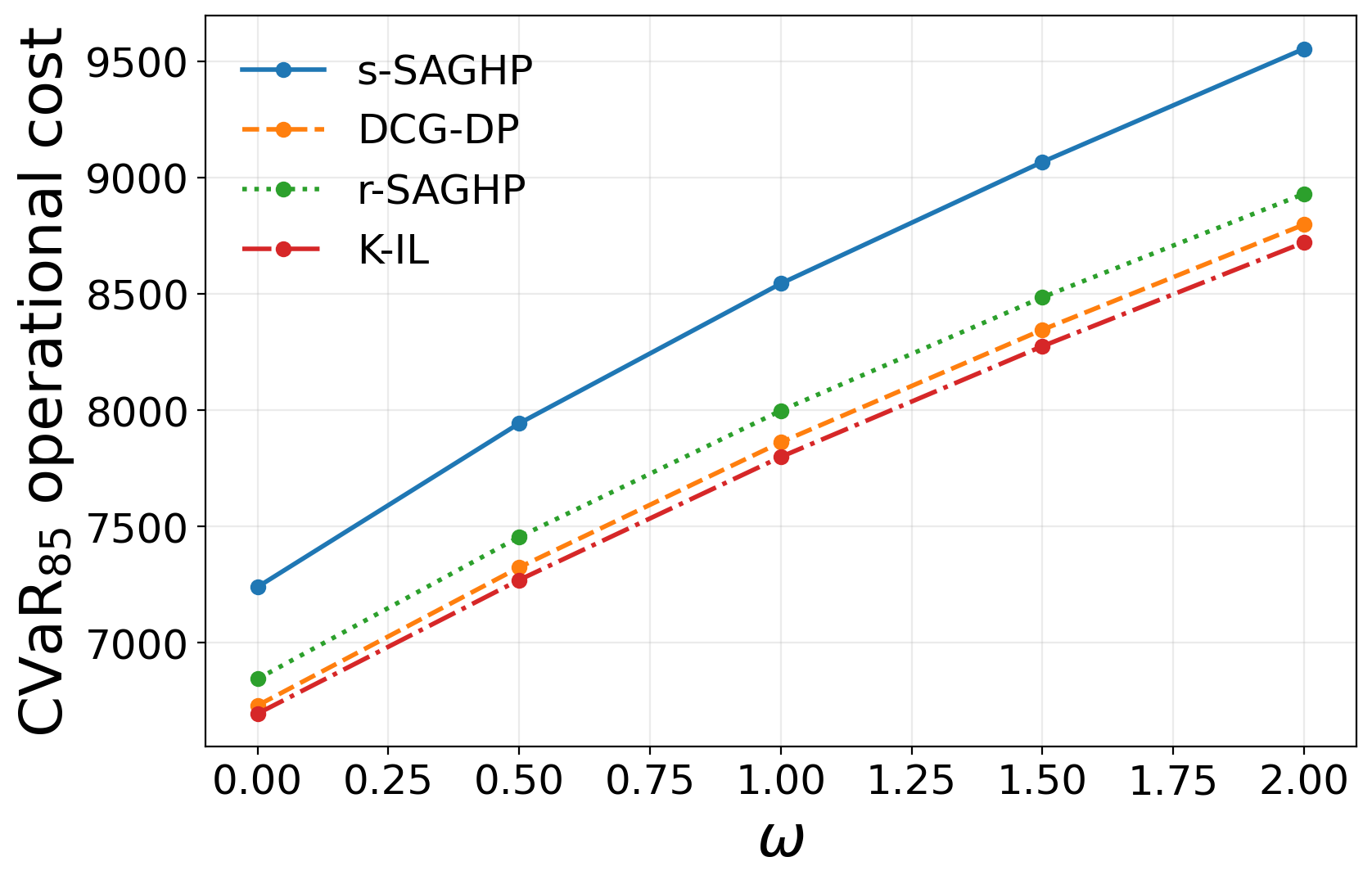}
        \caption{$\rho = 0.15$}
        \label{fig:oos_cvar_cov_rho_015}
    \end{subfigure}

    \vspace{0.5em}

    \begin{subfigure}[t]{0.47\textwidth}
        \centering
        \includegraphics[width=\linewidth]{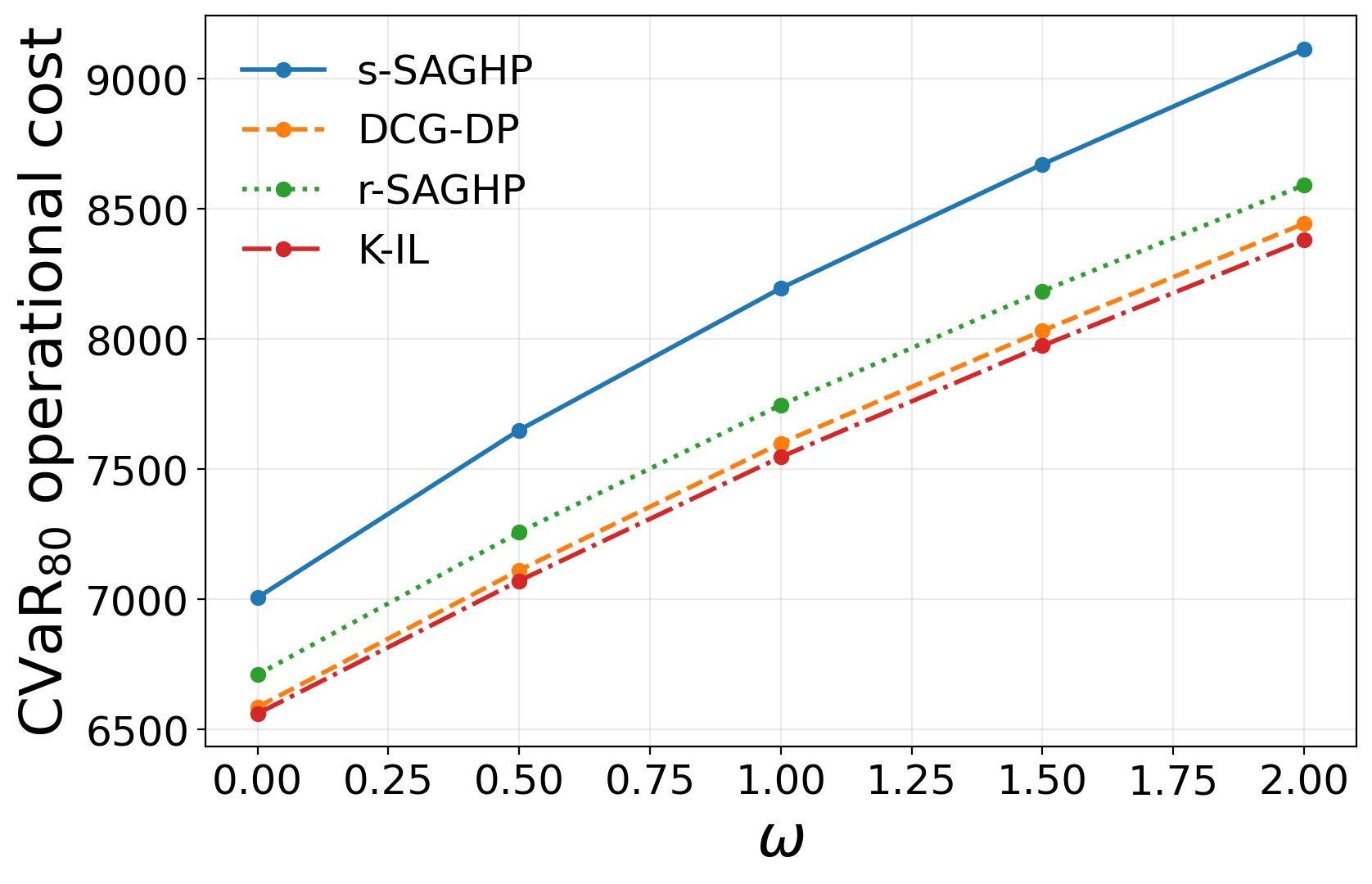}
        \caption{$\rho = 0.20$}
        \label{fig:oos_cvar_cov_rho_020}
    \end{subfigure}

    \caption{Out-of-sample performance comparison in CVaR with five confidence levels ($\rho$) under different levels of covariance increase ($\omega$).}
    \label{fig:oos_cvar_comparison_covariance_increase}
\end{figure}

We next examine the covariance-increase experiments. We evaluate out-of-sample mean operational costs under posterior covariance inflation, where the covariance matrix is scaled by $(1+\omega)$ for $\omega \in \{0.5,1.0,1.5,2.0\}$, and assess tail performance using $\mathrm{CVaR}_{\rho}$ for $\rho\in \{0.01,0.05,0.10,0.15,0.20\}$, as shown in Figures~\ref{fig:oos_mean_cost_covariance_increase} and~\ref{fig:oos_cvar_comparison_covariance_increase}. Because these experiments preserve the posterior mean and only increase the dispersion, s-SAGHP achieves the lowest mean operational cost across all values of $\omega$, as expected. In contrast, r-SAGHP incurs the highest mean operational cost because of its more conservative decisions. Between the two distributionally robust approaches, dr-SAGHP generally achieves a slightly lower mean cost than rs-dr-SAGHP.

The pattern is reversed when tail performance is considered. dr-SAGHP and rs-dr-SAGHP provide stronger CVaR protection than r-SAGHP across nearly all confidence levels, with the exception of $\rho=0.01$, for which r-SAGHP is slightly better than dr-SAGHP in some cases. Overall, the two distributionally robust approaches reduce tail costs by approximately $8\%$--$11\%$ relative to s-SAGHP across the five confidence levels, as reported in Appendix~\ref{sec:sensitivity_analysis_results}. These results show that, under increased predictive dispersion, dr-SAGHP achieves a favorable balance between average operational cost and tail-risk protection, whereas r-SAGHP tends to obtain comparable tail protection at the expense of greater conservatism in the mean.

\subsection{Key findings}\label{sec:key finding}
The numerical experiments yield four key insights.

First, DCG-DP provides a computational advantage over DCG-MILP
for solving the full-support dr-SAGHP formulation, achieving
speedups of approximately $3.18\times$--$6.24\times$ in the
computational benchmark. K-IL, which solves the support-restricted
rs-dr-SAGHP formulation, achieves speedups of approximately
$6.69\times$--$11.80\times$ relative to DCG-MILP. The former
comparison demonstrates the benefit of DP separation for the same
model, whereas the latter reflects differences in both the
formulation and the solution procedure.

Second, on the 84 held-out historical severe-weather days,
dr-SAGHP and rs-dr-SAGHP improve both average and tail
operational performance relative to s-SAGHP. Their mean operational
costs decrease by $2.60\%$ and $2.56\%$, respectively, and their
$\mathrm{CVaR}_{90}$ values decrease by approximately $8\%$.
Both distributionally robust models provide stronger tail protection
than r-SAGHP, whose mean operational cost increases by $0.31\%$
relative to s-SAGHP. These results demonstrate a favorable balance
between average cost and tail protection on the historical test set.

Third, after their Wasserstein radii are calibrated independently
using the training data, rs-dr-SAGHP achieves out-of-sample
performance comparable to that of dr-SAGHP. Although rs-dr-SAGHP
restricts the adversarial distribution to the support of the
estimated nominal distribution, it retains much of the operational
benefit of the full-support model in this case study. Together with
the runtime results for K-IL, this finding supports rs-dr-SAGHP
as a computationally efficient alternative.

Fourth, the sensitivity analysis under simulated capacity
distribution shifts illustrates how the models' relative performance
depends on the type and severity of the shift. Under mean capacity
reductions, dr-SAGHP and rs-dr-SAGHP generally achieve competitive
mean costs and tail protection compared with r-SAGHP. Under the
most severe reduction considered, $20\%$, r-SAGHP achieves a lower
mean operational cost than either distributionally robust model.
Under covariance inflation, the two distributionally robust models
generally provide stronger tail protection than r-SAGHP while
maintaining mean costs closer to those of s-SAGHP.

Overall, the historical-data case study and simulated sensitivity
analysis show that dr-SAGHP and rs-dr-SAGHP can improve tail
protection while maintaining competitive average operational costs.
The computational results establish DCG-DP as an effective algorithm
for the full-support model and identify rs-dr-SAGHP, solved using
K-IL, as an alternative offering substantial runtime savings.
\section{Conclusion}
In this paper, we proposed a distributionally robust framework for the single airport ground holding problem, referred to as dr-SAGHP, to address uncertainty and potential misspecification in airport capacity distributions. Building on classical ground delay program formulations and stochastic programming models, our work addresses a key practical limitation of conventional stochastic approaches: their dependence on an accurately specified probability distribution of future airport capacity. In operational environments characterized by severe weather, limited historical observations, and potential distributional shifts, this assumption may be difficult to justify.

To address this limitation, we formulate dr-SAGHP as a two-stage distributionally robust optimization problem with a Wasserstein ambiguity set and derive its deterministic equivalent following the standard reformulation procedure in \citep{c9}. Because the resulting formulation is computationally challenging, we develop three decomposition-based solution approaches. Two of these methods solve  model exactly through delayed constraint generation, using either a MILP separation algorithm or a dynamic programming separation algorithm; we refer to these methods as DCG-MILP and DCG-DP, respectively. We also develop an restricted-support approximation, rs-dr-SAGHP, which combines Kelley's cutting-plane method with an integer L-shaped framework. By exploiting the structure of the two-stage distributionally robust integer program with relatively complete recourse and continuous second-stage decisions, we derive an efficient procedure for computing the subgradients required by Kelley's method. The computational results demonstrate that K-IL achieves substantial runtime reductions relative to the exact methods while maintaining competitive out-of-sample performance when its Wasserstein radius is calibrated separately using historical data. Among the exact approaches, DCG-DP consistently improves computational efficiency over DCG-MILP while preserving the same solution quality.

Experiments on held-out historical data from DEN show that the calibrated dr-SAGHP and rs-dr-SAGHP models reduce mean operational cost by approximately $2.6\%$ and $\mathrm{CVaR}_{90}$ by approximately $8\%$ relative to s-SAGHP. The restricted-support approximation achieves comparable out-of-sample performance with lower computational cost.

Sensitivity tests under simulated mean-capacity reductions and covariance inflation further illustrate the trade-off between tail protection and conservatism. Across the tested settings, dr-SAGHP generally balances these objectives effectively, although r-SAGHP becomes more competitive under the most severe mean-capacity reductions.

Overall, our results demonstrate that dr-SAGHP provides a flexible and resilient framework for ground delay management under uncertain and potentially misspecified airport capacity distributions. The combination of a tractable DRO formulation, efficient decomposition algorithms, strong out-of-sample performance on historical severe-weather data, and robustness under simulated distribution shifts highlights the potential of distributionally robust ground delay programs as a practical tool for improving the resilience of air transportation systems. In particular, the proposed framework offers meaningful protection against adverse tail events while maintaining competitive average operational performance, and the K-IL approximation provides a computationally attractive alternative when solving  model is prohibitively expensive.

Several directions for future research naturally follow from this work. One important extension is to generalize the proposed framework to the multi-airport ground holding problem, where airport operations are coupled through connecting flights, shared airspace resources, weather systems, and network-level traffic interactions. Incorporating joint capacity distributions across the FAA Core~30 airports and explicitly modeling statistical dependence among airport capacities would provide a more realistic representation of the National Airspace System and enable coordinated delay-management strategies across multiple airports. From an algorithmic perspective, such coupling also creates additional computational challenges, since the separable structure exploited in the single airport dynamic programming oracle no longer holds directly.

A second direction is to extend the current two-stage formulation to a multi-stage setting. In practice, airport capacity forecasts evolve continuously as new weather information, demand conditions, and upstream congestion observations become available. A multi-stage distributionally robust ground delay model would allow decisions to be updated dynamically as uncertainty is gradually resolved. Developing such a framework would require integrating multi-stage DRO with scalable decomposition or approximation techniques while preserving sufficiently fast solution times for operational use.

Finally, the ambiguity set itself could be further enriched to incorporate additional sources of predictive structure and decision dependence. For example, future work could investigate data-driven ambiguity sets that explicitly account for temporal dependence, forecast covariates, or changing weather regimes. Advancing the model along these directions would further improve its realism and scalability and move distributionally robust ground delay optimization closer to real-world deployment in the National Airspace System.

\bibliographystyle{plainnat} 
\bibliography{sample} 
\newpage
\appendix

\section{Proof of Proposition-\ref{prop:dr_saghp_full_support}} \label{sec:dr_saghp_reformulation}

Let
\[
G(x)
:=
\sum_{f\in \mathcal{F}}C_g
\left(
\sum_{t\in \mathcal{T}_f}t x_{f,t}-r_f
\right)
\]
denote the first-stage ground-delay cost, and let
\(Q(x,\xi)\) denote the optimal second-stage airborne-holding cost under
capacity realization \(\xi\). The dr-SAGHP can be written as
\begin{equation}
\min_{x\in\mathcal X}
\left\{
c(x)
+
\max_{\mathbb Q\in
\mathcal P_{\varepsilon}(\widehat{\mathbb P})}
\mathbb E_{\mathbb Q}[Q(x,\xi)]
\right\},
\label{eq:dr_saghp_original}
\end{equation}
where \(\mathcal X\) denotes the feasible set of first-stage constraints.

\paragraph{Step 1: Continuous-support transport representation.}

Let $
\widehat{\mathbb P}
=
\sum_{i=1}^{n}\widehat p_i
\delta_{\widehat\xi^{i}}
$ and $\mathbb Q \in \mathcal M(\Xi)$, the 1-Wasserstein distance between
\(\widehat{\mathbb P}\) and \(\mathbb Q\) can be written as
\begin{equation}
\begin{aligned}
W_1(\widehat{\mathbb P},\mathbb Q)
= \\
\inf_{\pi\in\mathcal M_+(\Xi)}
\quad
&\sum_{i=1}^{n}
\int_{\Xi}
\|\xi-\widehat\xi^{i}\|_1
\,\pi(\widehat{\xi}_{i}, \xi)d\xi
\\
\textnormal{s.t.}\quad
&\int_{\Xi}\pi(\widehat{\xi}_{i}, \xi)d\xi=\widehat p_i,
&&\forall i\in[n],
\\
&\sum_{i=1}^{n}\pi(\widehat{\xi}_{i},\xi)=d\mathbb Q(\xi),
&&\xi \in \Xi,
\end{aligned}
\label{eq:wasserstein_transport_measure}
\end{equation}
where \(\mathcal M(\Xi)\) denotes the set of finite nonnegative
measures supported on \(\Xi\).

Substituting this representation into the worst-case expectation gives
\begin{equation}
\begin{aligned}
\max_{\mathbb Q\in
\mathcal P_{\varepsilon}(\widehat{\mathbb P})} &
\mathbb E_{\mathbb Q}[Q(x,\xi)]
= \\
\max_{\pi\in\mathcal M_+(\Xi)}
&
\int_{\Xi}
Q(x,\xi)\,d\mathbb Q(\xi)
\\
\textnormal{s.t.}\quad
&\sum_{i=1}^{n}
\int_{\Xi}
\|\xi-\widehat\xi^{i}\|_1
\,\pi(\widehat\xi^{i}, \xi)d\xi
\leq\varepsilon \\ 
&\int_{\Xi}\pi(\widehat\xi^{i}, \xi)d\xi=\widehat p_i,
&&\forall i\in[n],
\\
&\sum_{i=1}^{n}\pi(\widehat{\xi}_{i},\xi)=d\mathbb Q(\xi),
&&\xi \in \Xi.
\end{aligned}
\label{eq:continuous_transport_problem}
\end{equation}

The measures \(\pi_i\) are supported on the continuous physical support
\(\Xi\). 


\paragraph{Step 2: Wasserstein dual reformulation.}

Let \(\lambda\geq 0\) be the dual variable associated with the
Wasserstein transportation-budget constraint in
\eqref{eq:continuous_transport_problem}. Since \(\Xi\) is compact and
\(Q(x,\xi)\) is finite and continuous in \(\xi\), strong duality for the
Wasserstein moment problem gives
\begin{equation}
\begin{aligned}
&
\max_{\mathbb Q\in
\mathcal P_{\varepsilon}(\widehat{\mathbb P})}
\mathbb E_{\mathbb Q}[Q(x,\xi)]
\\
&\qquad =
\min_{\lambda\geq0}
\left\{
\lambda\varepsilon
+
\sum_{i=1}^{n}\widehat p_i
\max_{\xi\in\Xi}
\left[
Q(x,\xi)
-
\lambda\|\xi-\widehat\xi^{i}\|_1
\right]
\right\}.
\end{aligned}
\label{eq:wasserstein_dual_reformulation}
\end{equation}

Introducing epigraph variables
\(\alpha_i\), \(i=1,\ldots,n\), gives
\begin{equation}
\begin{aligned}
\min_{\lambda,\alpha}\quad
&\lambda\varepsilon
+\sum_{i=1}^{n}\widehat p_i\alpha_i
\\
\textnormal{s.t.}\quad
&
\alpha_i
\geq
\max_{\xi\in\Xi}
\left\{
Q(x,\xi)
-
\lambda\|\xi-\widehat\xi^{i}\|_1
\right\},
&&\forall i\in[n],
\\
&\lambda\geq0.
\end{aligned}
\label{eq:wasserstein_epigraph}
\end{equation}

Equivalently, the epigraph constraints may be written in semi-infinite
form as
\begin{equation}
\alpha_i
\geq
Q(x,\xi)
-\lambda\|\xi-\widehat\xi^{i}\|_1,
\qquad
\forall \xi\in\Xi,\quad
\forall i\in[n].
\label{eq:wasserstein_semi_infinite}
\end{equation}

\paragraph{Step 3: Dual representation of the recourse function.}


For a fixed \(x\) and \(\xi\), the second-stage recourse problem is
\begin{equation}
\begin{aligned}
Q(x,\xi)
=
\min_{y}\quad
& C_h\sum_{t=0}^{T}y_t
\\
\textnormal{s.t.}\quad
&y_0\geq A_0(x)-\xi_0,
\\
&y_t-y_{t-1}\geq A_t(x)-\xi_t,
&&t=1,\ldots,T,
\\
&y_t\geq0,
&&t=0,\ldots,T.
\end{aligned}
\label{eq:second_stage_recourse}
\end{equation}

By linear programming strong duality,
\begin{equation}
Q(x,\xi)
=
\max_{\nu\in\mathcal V}
\left\{
A(x)^\top\nu-\xi^\top\nu
\right\}.
\label{eq:recourse_dual}
\end{equation}

Because \(\mathcal V\) is a bounded polytope and the dual objective is
linear in \(\nu\), it is sufficient to consider its extreme points:
\begin{equation}
Q(x,\xi)
=
\max_{\nu\in\operatorname{ext}(\mathcal V)}
\left\{
A(x)^\top\nu-\xi^\top\nu
\right\}.
\label{eq:recourse_extreme_points}
\end{equation}

Substituting \eqref{eq:recourse_extreme_points} into
\eqref{eq:wasserstein_epigraph} and interchanging the two maximizations
gives
\begin{equation}
\alpha_i
\geq
A(x)^\top\nu
+
\max_{\xi\in\Xi}
\left\{
-\xi^\top\nu
-\lambda\|\xi-\widehat\xi^{i}\|_1
\right\}, \\
\forall \nu\in\operatorname{ext}(\mathcal V),\;
\forall i.
\label{eq:epigraph_after_recourse_dual}
\end{equation}

\paragraph{Step 4: Evaluation over the continuous box support.}

Using the dual-norm representation of the \(\ell_1\) norm,
\begin{equation}
\lambda\|\xi-\widehat\xi^{i}\|_1
=
\max_{\|z\|_\infty\leq\lambda}
z^\top(\xi-\widehat\xi^{i}),
\label{eq:l1_dual_norm}
\end{equation}
the inner maximization in
\eqref{eq:epigraph_after_recourse_dual} becomes
\begin{align}
&
\max_{\xi\in\Xi}
\left\{
-\xi^\top\nu
-\lambda\|\xi-\widehat\xi^{i}\|_1
\right\}
\notag\\
&\quad =
\max_{\xi\in\Xi}
\min_{\|z\|_\infty\leq\lambda}
\left\{
z^\top\widehat\xi^{i}
+
(-\nu-z)^\top\xi
\right\}
\notag\\
&\quad =
\min_{\|z\|_\infty\leq\lambda}
\left\{
z^\top\widehat\xi^{i}
+
\max_{\xi\in\Xi}(-\nu-z)^\top\xi
\right\}.
\label{eq:minimax_interchange}
\end{align}
The minimax interchange follows because both \(\Xi\) and
\(\{z:\|z\|_\infty\leq\lambda\}\) are compact and convex, and the
objective is bilinear in \((\xi,z)\).

For the box support
\[
\Xi
=
\left\{
\xi:
0\leq\xi_t\leq\overline K,\quad
t\in\mathcal T
\right\},
\]
the support function is separable across time:
\begin{equation}
\begin{aligned}
\max_{\xi\in\Xi}(-\nu-z)^\top\xi
=
\sum_{t\in\mathcal T}
\left(
\overline K[-\nu_t-z_t]_+
-
\underline K[\nu_t+z_t]_+
\right),
\end{aligned}
\label{eq:box_support_function}
\end{equation}
where \([a]_+:=\max\{a,0\}\). In dr-SAGHP, we assume that the physically feasible lower bound $\underline{K}=0$. Hence,
\begin{equation}
\begin{aligned}
\alpha_i
&\geq
A(x)^\top\nu
\\&+
\min_{\|z\|_\infty\leq\lambda}
\Bigg\{
z^\top\widehat\xi^{i}
+
\sum_{t\in\mathcal T}
\Big(
\overline K[-\nu_t-z_t]_+
\Big)
\Bigg\},
\end{aligned}
\label{eq:full_support_constraint}
\end{equation}
for every \(i\in[n]\) and
\(\nu\in\operatorname{ext}(\mathcal V)\).

\paragraph{Step 5: Deterministic equivalent formulation.}

Combining the Wasserstein reformulation, the dual representation of the
recourse function, and the first-stage feasibility constraints yields
the following deterministic equivalent of  dr-SAGHP:
\begin{subequations}\label{eq:dr_saghp_full_support_p}
\begin{align}
\min_{x,\lambda,\alpha}\quad
& \sum_{f\in \mathcal{F}} C_g
  \left(
    \sum_{t\in \mathcal{T}_f} t x_{f,t}-r_f
  \right)
  +\lambda\varepsilon
  +\sum_{i=1}^{N}\widehat p_i\alpha_i
\label{eq:dr_saghp_full_support_obj_p}
\\
\text{s.t.}\quad
& \alpha_i
  \geq
  A(x)^\top\nu
    \\& \qquad +
  \min_{\|z\|_\infty\leq\lambda}
  \left\{
    z^\top\widehat\xi^{i}
    +
    \sum_{t\in\mathcal T}
    \overline K_t[-\nu_t-z_t]_+
  \right\},
\notag\\[-2mm]
&
  \forall i\in[n],\;
  \forall\nu\in\operatorname{ext}(\mathcal V),
\label{eq:dr_saghp_full_support_dro_p}
\\
& \eqref{eq:dSAGHP_1c} - \eqref{eq:dSAGHP_1e} \\
& \lambda \geq 0,\\
& \alpha_i\in\mathbb R, \notag \forall i\in[N].
\label{eq:dr_saghp_full_support_domain_p}
\end{align}
\end{subequations} \hfill\Halmos

\section{Proof of Lemma~\ref{lemma:separation_reformulation}}\label{sec:proof_lemma1}
\begin{proof}{Proof}
We exploit the fact that the inner problem $\Psi(\nu,\lambda,\widehat\xi^{i})$ is separable across the components $z_{it}$ of $z_i$, and thus its analytical solution can be obtained. 

For the fixed master solution
\((\bar x,\bar\lambda,\bar\alpha)\), $\operatorname{SEP}_i(\bar x, \bar \lambda)$ is
\begin{equation}
\begin{aligned}
\max_{\nu\in\mathcal V}
\Bigg\{
A(\bar x)^\top\nu
+
\min_{\|z\|_\infty\leq\bar\lambda}
z^\top\widehat\xi^{i}
+
\sum_{t=0}^{T}
\bigl(
\overline K_t[-\nu_t-z_t]_+
\Bigg\}.
\label{eq:milp_sep_original}
\end{aligned}
\end{equation}

Because the \(\ell_\infty\)-norm constraint applies componentwise, the inner
minimization is separable across time, i.e,
\begin{align*}
     &\min_{\|z\|_\infty\leq\bar\lambda}
z^\top\widehat\xi^{i}
+
\sum_{t=0}^{T}
\overline K_t[-\nu_t-z_t]_+=\\
&\qquad \sum_{t=0}^T \min_{-\bar\lambda\leq z_t\leq\bar\lambda}
\left\{
\widehat\xi_t^i z_t
+
\overline K_t[-\nu_t-z_t]_+
\right\}
\end{align*}
Fix some time period $t$ and corresponding dual variable $\nu_t$, and consider the problem 
\begin{equation}
\Psi^{(t)}(\nu,\bar\lambda,\widehat\xi^{i}) = \min_{-\bar\lambda\leq z_t\leq\bar\lambda}
\left\{
\widehat\xi_t^i z_t
+
\overline K_t[-\nu_t-z_t]_+
\right\}.
\label{eq:scalar_z_problem}
\end{equation}

Recall that \(\nu_t\geq0\), \(\bar\lambda\geq0\), and
\(0\leq\widehat\xi_t^i\leq\overline K_t\).
We distinguish two cases according to whether the breakpoint
\(z_t=-\nu_t\) lies within the feasible interval
\([-\bar\lambda,\bar\lambda]\).

\noindent\textbf{Case 1: \(\nu_t\leq\bar\lambda\).}
In this case, \(z_t=-\nu_t\) is feasible.
On the subinterval \([-\bar\lambda,-\nu_t]\), the objective
in~\eqref{eq:scalar_z_problem} is
$-\overline K_t\nu_t
+\left(\widehat\xi_t^i-\overline K_t\right)z_t$. Since \(\widehat\xi_t^i\leq\overline K_t\), this expression is
non-increasing in \(z_t\) and attains its minimum at the largest
point in this subinterval, \(z_t=-\nu_t\).
On the subinterval \([-\nu_t,\bar\lambda]\), the objective is
\(\widehat\xi_t^i z_t\), which is non-decreasing in \(z_t\)
because \(\widehat\xi_t^i\geq0\).
It therefore attains its minimum at the smallest point in this
subinterval, again \(z_t=-\nu_t\).
Thus, both subintervals attain the same minimum value
\(-\widehat\xi_t^i\nu_t\) at their common endpoint,
and \(z_t^*=-\nu_t\) is an optimizer over the entire feasible interval.

\noindent\textbf{Case 2: \(\nu_t>\bar\lambda\).}
Every feasible \(z_t\) satisfies
\(z_t\geq-\bar\lambda>-\nu_t\).
Hence, \([-\nu_t-z_t]_+=0\) throughout the feasible interval,
and the objective reduces to \(\widehat\xi_t^i z_t\).
This expression is non-decreasing in \(z_t\), so it attains its
minimum at the lower endpoint, \(z_t^*=-\bar\lambda\),
with minimum value \(-\widehat\xi_t^i\bar\lambda\).

Combining the two cases, an optimizer
of~\eqref{eq:scalar_z_problem} is
\begin{equation}
z_t^*
=
\begin{cases}
-\nu_t, & \nu_t\leq\bar\lambda,\\
-\bar\lambda, & \nu_t>\bar\lambda,
\end{cases}
=-\min\{\nu_t,\bar\lambda\}.
\label{eq:optimal_z}
\end{equation}
The corresponding optimal objective value is
\(-\widehat\xi_t^i\min\{\nu_t,\bar\lambda\}\).
\label{eq:w_definition}
Substituting \(z_t^*=-\min\{\nu_t,\bar\lambda\}\) into
\eqref{eq:scalar_z_problem} gives
$\Psi^{(t)}(\nu,\bar\lambda,\widehat\xi^{i})=-\widehat\xi_t^i \min\{\nu_t,\bar\lambda\}
+ \overline{K}[\min\{\nu_t,\bar\lambda\}\ - \nu_t]_{+}.$
Since $\min\{\nu_t,\bar\lambda\}\ \leq \nu_{t}$, $[\min\{\nu_t,\bar\lambda\}\ - \nu_t]_{+} = 0$ and thus
\begin{equation}
\Psi^{(t)}(\nu,\bar\lambda,\widehat\xi^{i})=\widehat{\xi}^i\min\{\nu_t,\bar\lambda\}
\end{equation}

Therefore, the separation problem reduces to
\begin{equation}
\begin{aligned}
\operatorname{SEP}_{i}(\bar x,\bar\lambda)
&=
\max_{\nu\in\mathcal V}
\left\{
\sum_{t=0}^{T}
A_t(\bar x)\nu_t
-
\sum_{t=0}^{T}
\widehat\xi_t^i
\min\{\nu_t,\bar\lambda\}
\right\} \\
&= \max_{\nu\in\mathcal{V}}
    \sum_{t\in\mathcal{T}} r_{i,t}(\nu)
\end{aligned}
\label{eq:separation_min_reduced}
\end{equation} and thus we complete the proof.
\hfill \Halmos
\end{proof}
\section{Derivation of MILP separation algorithm}\label{sec:milp_separation}
In this section, we derive a mixed-integer linear programming (MILP) formulation of the separation problem defined in Equation~\eqref{eq:dcg_separation_problem}. Based on the reformulation of $\operatorname{SEP}_{i}(\bar x, \bar\lambda)$ in Lemma~\ref{lemma:separation_reformulation}, 
we introduce an auxiliary variable \(w_t=\min\{\nu_t,\bar\lambda\}\), binary variables \(b_t\in\{0,1\}\).
If \(b_t=0\),
\eqref{eq:milp_separation_branch_1} implies $\nu_t\leq\bar\lambda$, while
\eqref{eq:milp_separation_w_1} and
\eqref{eq:milp_separation_branch_2} imply
$
w_t=\nu_t.
$
If \(b_t=1\), constraints
\eqref{eq:milp_separation_w_2} and
\eqref{eq:milp_separation_branch_3} imply
$
w_t=\bar\lambda.
$
Thus, constraints
\eqref{eq:milp_separation_w_1}--%
\eqref{eq:milp_separation_binary} enforce
$
w_t=\min\{\nu_t,\bar\lambda\}.
$
Substitution into \eqref{eq:separation_min_reduced} yields the following MILP reformulation: 
\begin{subequations}\label{eq:milp_separation}
\begin{align}
\max_{\nu,w,b}\quad
&
\sum_{t=0}^{T}
A_t(\bar x)\nu_t
-
\sum_{t=0}^{T}
\widehat\xi_t^{(i)}w_t
\label{eq:milp_separation_obj}
\\
\textnormal{s.t.}\quad
&
\nu_t-\nu_{t+1}\leq C_h,
&& t \in [T-1]
\label{eq:milp_separation_dual_1}
\\
&
\nu_T\leq C_h,
\label{eq:milp_separation_dual_2}
\\
&
0\leq w_t\leq\nu_t,
&& t\in [T],
\label{eq:milp_separation_w_1}
\\
&
w_t\leq\bar\lambda,
&& t \in [T],
\label{eq:milp_separation_w_2}
\\
&
\nu_t\leq\bar\lambda+\mathbb{M} b_t,
&&t \in [T],
\label{eq:milp_separation_branch_1}
\\
&
w_t\geq\nu_t-\mathbb{M} b_t,
&&t \in [T],
\label{eq:milp_separation_branch_2}
\\
&
w_t\geq\bar\lambda b_t,
&&t \in [T],
\label{eq:milp_separation_branch_3}
\\
&
b_t\in\{0,1\},
\label{eq:milp_separation_binary}
\end{align}
\end{subequations}
Here, the $\mathbb M$ indicates a big-$\mathbb M$ coefficient.

\section{Proof of Proposition~\ref{prop:dp_separation}}\label{sec: proof-dp-separation} 
\begin{proof}{Proof}
Because $\widehat\xi^{i}_t\geq 0$, the function
\begin{equation}
    \nu_t
    \longmapsto
    \widehat\xi^{i}_t
    [\nu_t-\bar{\lambda}]_+
\end{equation}
is convex. Therefore, $r_{i,t}$ is a convex piecewise-linear function of
$\nu_t$, and the separation objective is convex in $\nu$.
Since $\mathcal V$ is a nonempty compact polytope, every
$\nu\in\mathcal V$ can be represented as a convex combination
of its extreme points
\citep[Theorem~2.9]{bertsimas1997introduction}.
Consequently, an optimal solution can be chosen from $\operatorname{ext}(\mathcal{V})$, and hence
\begin{equation}
    \max_{\nu\in\mathcal{V}}\sum_{t \in \mathcal{T}}r_{i, t}
    =
    \max_{\nu\in\operatorname{ext}(\mathcal{V})}
    \sum_{t \in \mathcal{T}}r_{i, t}.
    \label{eq:maximum_at_extreme_point}
\end{equation}
We next show that
\begin{equation}
    \operatorname{ext}(\mathcal{V})
    \subseteq
    \mathcal{E}.
    \label{eq:ext_subset_structured_set}
\end{equation}
Using the convention $\nu_{T+1}=0$, the dual feasible region can be
written as
\begin{equation}
    \mathcal{V}
    =
    \left\{
    \nu\in\mathbb{R}^{T+1}:
    -\nu_t\leq 0,\quad
    \nu_t-\nu_{t+1}\leq C_h,
    \quad t=0,\ldots,T
    \right\}.
    \label{eq:dual_feasible_region_compact}
\end{equation}
The polytope $\mathcal{V}$ is full-dimensional because, for any
$\delta\in(0,C_h)$, the point defined by $\nu_t=\delta$ for all
$t=0,\ldots,T$ satisfies all the inequalities in
\eqref{eq:dual_feasible_region_compact} strictly. Therefore, every
extreme point of $\mathcal{V}\subseteq\mathbb{R}^{T+1}$ must have
$T+1$ linearly independent active constraints
\citep{bertsimas1997introduction}.
For each $t$, the constraint $\nu_{t} \geq 0$ is active if and only if $\nu_t=0$, whereas $\nu_t-\nu_{t+1}\leq C_h$ is active if and only if $\nu_t=\nu_{t+1}+C_h$. These two constraints cannot be active simultaneously. Otherwise, it would imply $\nu_{t+1}=-C_h<0$, contradicting the nonnegativity of
$\nu_{t+1}$. For $t=T$, simultaneous activity would require both
$\nu_T=0$ and $\nu_T=C_h$, which is impossible because $C_h>0$.

Thus, at most one constraint associated with each $t$ can be active.
Since there are $T+1$ time indices, an extreme point can have at most
$T+1$ active constraints. To obtain the $T+1$ linearly independent
active constraints required at an extreme point, exactly one constraint
associated with each $t$ must be active. Therefore, every
$\nu\in\operatorname{ext}(\mathcal{V})$ satisfies
\begin{equation}
    \nu_t=0
    \quad\text{or}\quad
    \nu_t=\nu_{t+1}+C_h,\quad \forall t.
\end{equation}
This proves
$
    \operatorname{ext}(\mathcal{V})
    \subseteq
    \mathcal{E}.
$
Moreover, $\mathcal{E}\subseteq\mathcal{V}$ by definition. Combining
these inclusions with \eqref{eq:maximum_at_extreme_point} gives
\begin{equation}
\begin{aligned}
    \max_{\nu\in\mathcal{V}}\sum_{t\in\mathcal T}r_{i,t}
    &=
    \max_{\nu\in\operatorname{ext}}\sum_{t\in\mathcal T}r_{i,t} \\
    &\leq
    \max_{\nu\in\mathcal{E}}\sum_{t\in\mathcal T}r_{i,t} \\
    &\leq
    \max_{\nu\in\mathcal{V}}\sum_{t\in\mathcal T}r_{i,t}.
\end{aligned}
\end{equation}
All inequalities must therefore hold with equality, which proves 
\eqref{eq:structured_separation_problem} and completes this proof.
\end{proof}
\hfill\Halmos

\section{DP separation algorithm}\label{sec:alg_dp_separation}
\begin{algorithm}[H]
\caption{Exact dynamic programming separation algorithm}
\label{alg:dp_separation}
\begin{algorithmic}[1]

\Require Current RMP solution
$(\bar{x},\bar{\lambda},\bar{\alpha}_i)$
and a sampled capacity trajectory $\widehat\xi^{i}$

\Ensure Separation value $\operatorname{SEP}_{i}(x,\lambda)$
and optimal dual vector $\nu_i^*$

\For{$t=0,\ldots,T$}
    \ForAll{$s\in\mathcal{S}_t$}
        \State Compute
        $
        r_{i,t}(s)
        \gets
        \left(
        A_t(\bar{x})-\widehat\xi^{i}_t
        \right)s
        +
        \widehat\xi^{i}_t
        [s-\bar{\lambda}]_+
        $
    \EndFor
\EndFor

\ForAll{$s\in\mathcal{S}_T$}
    \State
    $
    V_{i,T}(s)\gets r_{i,T}(s)
    $
\EndFor


\For{$t=T-1,T-2,\ldots,0$}
    \State Choose
    $P_{i,t}(0)\in
    \arg\max_{s'\in\mathcal S_{t+1}}V_{i,t+1}(s')$
    \State $V_{i,t}(0)\gets
    r_{i,t}(0)+V_{i,t+1}(P_{i,t}(0))$

    \ForAll{$s\in\mathcal S_t\setminus\{0\}$}
        \State $P_{i,t}(s)\gets s-C_h$
        \State $V_{i,t}(s)\gets
        r_{i,t}(s)+V_{i,t+1}(s-C_h)$
    \EndFor
\EndFor

\State
$
\nu_{i,0}^*
\gets
\arg\max_{s\in\mathcal{S}_0}
V_{i,0}(s)
$ 

\For{$t=0,\ldots,T-1$}
    \State
    $
    \nu_{i,t+1}^*
    \gets
    P_{i,t}\bigl(\nu_{i,t}^*\bigr)
    $ \qquad (\text{Recover the optimal dual trajectory})
\EndFor

\State $
\operatorname{SEP}_{i}(x,\lambda)
\gets
V_{i,0}\bigl(\nu_{i,0}^*\bigr)
$

\State \Return
$\operatorname{SEP}_{i}(x,\lambda),\nu_i^*$

\end{algorithmic}
\end{algorithm}

\section{Proof of Proposition~\ref{prop: subgrad}}\label{sec: proof-subgrad}
Proposition~\ref{prop: subgrad} follows from an application of Danskin's Theorem. We use a variant from \cite{bertsekas1971control}:
\begin{theorem}[Danskin's Theorem \cite{bertsekas1971control}]Let $\mathbf{g}:\mathbb{R}^{n}\times\mathbb{R}^{m}\rightarrow (-\infty, \infty)$ be a function and let $\mathbf{Y}$ be a compact subset of $\mathbb{R}^{m}$. Assume further that for every vector $\mathbf{y} \in \mathbf{Y}$ the function $\mathbf{g(\cdot, \mathbf{y})}:\mathbb{R}^{n}\rightarrow (-\infty, \infty)$ is a closed proper convex function. Consider the function $\mathbf{f}$ defined as $\mathbf{f}(\mathbf{x}) = \max_{\mathbf{y} \in \mathbf{Y}}\mathbf{g}(\mathbf{x}, \mathbf{y})$, if $\intt(\dom \mathbf{f}) \neq \emptyset$ and $\mathbf{g}$ is continuous on the set $\intt(\dom \mathbf{f} )\times \mathbf{Y}$ for every $\mathbf{x} \in \intt(\dom \mathbf{f})$ we have 
\begin{equation*}
\begin{aligned}
\partial\mathbf{f}(\mathbf{x}) = \conv\{\partial\mathbf{g}(\mathbf{x},\mathbf{y}^{*})|\mathbf{y}^{*}\in \mathbf{Y}^{*}(\mathbf{x})\},
\end{aligned}
\end{equation*}
where 
\begin{equation*}
\begin{aligned}
\mathbf{Y}^{*}(\mathbf{x}) = \{\mathbf{y}^{*}\in \mathbf{Y}|\mathbf{g}(\mathbf{x},\mathbf{y}^{*}) = \max_{\mathbf{y} \in \mathbf{Y}}\mathbf{g}(\mathbf{x},\mathbf{y})\}.
\end{aligned}
\end{equation*}
\label{thm:danskin's}
\end{theorem}
In order to verify that the conditions of the theorem are satisfied, we must verify the following:
\begin{itemize}
    \item The set $\mathcal{P}_\varepsilon$ is compact. This is true, as this set is a bounded polytope.
    \item For each probability distribution $p$, the function $\mathbf{x} \mapsto \sum_{i=1}^{|\widehat{\Xi}|} p_i h(\mathbf{x},\xi_i)$ is a closed proper convex function, letting this function take the value $\infty$ when $\mathbf{x} \notin [0,1]^{\mathbf{N}}$ This function is closed, as it is infinite outside the unit cube $[0,1]^{\mathbf{N}}$.The fact that it is proper and convex follows from the assumption that $h(\cdot,\xi_i)$ is convex and the assumption that $\sum_{i=1}^{|\widehat{\Xi}|} p_i h(\mathbf{x},\xi_i)$ is finite for all $\mathbf{x}$ and $p \in \mathcal{P}_{\varepsilon}$.
    \item The interior of the domain of $f$ is non-empty. This is true as the domain of $f$ is $[0,1]^{\mathbf{N}}$.
    \item The function $(\mathbf{x},p) \mapsto \sum_{i=1}^{|\widehat{\Xi}|} p_i h(\mathbf{x},\xi_i)$ is continuous on $(0,1)^{\mathbf{N}}\times \mathcal{P}_\varepsilon$. This can be shown to be true. Fix any $i$. The $i$th coordinate projection $p \mapsto p_i$ is linear and therefore continuous.
Moreover, each recourse function $h(\cdot,\xi_i)$ is assumed to be convex, and therefore continuous on the interior of its domain.
Since $h(\mathbf{x},\xi_{i})$ does not depend on $p$, it follows that $(\mathbf{x},p)\mapsto h(\mathbf{x},\xi_{i})$ is continuous in $(\mathbf{x},p)$ on
$\intt(\dom \mathbf{f})\times \mathcal{P}_{\varepsilon}$. Because the product of continuous functions is continuous, $(x,p) \mapsto p_ih(\mathbf{x},\xi_{i})$ is continuous in $(x,p)$
on $\intt(\dom \mathbf{f})\times \mathcal{P}_{\varepsilon}$ for each $i$.
Finally, $\mathbf{g}$ is a finite sum of continuous functions, and is therefore continuous on
$\intt(\dom \mathbf{f})\times \mathcal{P}_{\varepsilon}$.
\end{itemize}
Then, Proposition~\ref{prop: subgrad} follows immediately.

\section{Proof of Proposition~\ref{prop:1-dim-dual}}\label{sec:proof-1-dim-dual}
Let $(\lambda,\alpha)$ be a feasible solution to problem~\eqref{eq:dr-general-inner_dual}, and let $\lambda^*$ be an optimal solution to the problem~\eqref{eq:one_dimensional_inner_dual}. Then, since $(\lambda,\alpha)$ is feasible, we have that
\begin{align*}
    \alpha_i \geq h(\mathbf{x},\xi_j) - \lambda \mathbf{d}_{i,j} 
\end{align*}
for each $i\in\{1,\ldots,|\widehat{\Xi}|\}$ and $j\in\{1,\ldots,|\widehat{\Xi}|\}$. This further implies that
\begin{align*}
    \alpha_i \geq \max_{j \in \{1,\ldots,|\widehat{\Xi}|\}}(h(\mathbf{x},\xi_j) - \lambda \mathbf{d}_{i,j} )
\end{align*}
for each $i \in \{1,\ldots,|\widehat{\Xi}|\}$. Then,
\begin{align*}
    \varepsilon \lambda + \sum_{i=1}^{|\widehat{\Xi}|} p_i\alpha_i &\geq \varepsilon \lambda + \sum_{i=1}^{|\widehat{\Xi}|} p_i\max_{j \in \{1,\ldots,|\widehat{\Xi}|\}}(h(\mathbf{x},\xi_j) - \lambda \mathbf{d}_{i,j} )\\
    &= \phi(\lambda) \\
    & \geq \phi(\lambda^*)\\
    &= \varepsilon \lambda + \sum_{i=1}^{|\widehat{\Xi}|} p_i\alpha_i^*(\lambda^*)
\end{align*}
In addition, the solution $(\lambda^*,\alpha^*(\lambda^*))$ is feasible, as it follows immediately from the definition that
\begin{align*}
    \alpha^*_i(\lambda^*) \geq h(\mathbf{x},\xi_j) - \lambda \mathbf{d}_{i,j} 
\end{align*}
for each $i\in\{1,\ldots,|\widehat{\Xi}|\}$ and $j\in\{1,\ldots,|\widehat{\Xi}|\}$. Thus, $(\lambda^*,\alpha^*(\lambda^*))$ is optimal. \hfill\Halmos
\section{Proof of Proposition~\ref{prop: dual-converge}}\label{sec:proof-dual-converge}
Let $\lambda_{\textup{hi}}^{(i)}$ and $\lambda_{\textup{lo}}^{(i)}$ denote the values taken by these variables in iteration $i$, with $\lambda_{\textup{hi}}^{(0)}$ and $\lambda_{\textup{lo}}^{(0)}$ denoting the initial values. We note that $\lambda_{\textup{hi}}^{(i+1)} - \lambda_{\textup{lo}}^{(i+1)} = (\lambda_{\textup{hi}}^{(i)} - \lambda_{\textup{lo}}^{(i)})/2$. It follows by a routine induction argument that $\lambda_{\textup{hi}}^{(i)} - \lambda_{\textup{lo}}^{(i)} = (\lambda_{\textup{hi}}^{(0)} - \lambda_{\textup{lo}}^{(0)})/2^i$. Thus, if $i \geq \lceil \log_2((\lambda^{(0)}_{\textup{hi}}-\lambda_{\textup{lo}}^{(0)})/\delta)\rceil$, then we have that
\begin{align*}
\lambda_{\textup{hi}}^{(i)} - \lambda_{\textup{lo}}^{(i)} \leq \dfrac{\lambda_{\textup{hi}}^{(0)} - \lambda_{\textup{lo}}^{(0)})}{(\lambda^{(0)}_{\textup{hi}}-\lambda_{\textup{lo}}^{(0)})/\delta} = \delta < 2\delta
\end{align*}
Thus, the algorithm stops in at most $\lceil \log_2((\lambda^{(0)}_{\textup{hi}}-\lambda_{\textup{lo}}^{(0)})/\delta)\rceil$ iterations.

Let $\lambda^{(i)} = (\lambda_{\textup{lo}}^{(i-1)}+\lambda_{\textup{hi}}^{(i-1)})/2$ and let $s^{(i)} = \varepsilon - \sum_{i=1}^{\abs{\widehat{\Xi}}} p_i\, \mathbf{d}_{i j^*(i,\lambda^{(i)})}$. In any iteration $i$, if $s^{(i)} = 0$, then $\lambda^{(i)} \in \argmin_{\lambda \geq 0} \phi(\lambda)$. 

If no such iteration exists, then it suffices to prove that $[\lambda_{\textup{lo}}^{(i)},\lambda_{\textup{hi}}^{(i)}] \cap \argmin_{\lambda \geq 0} \phi(\lambda) $ is non-empty for each $i$. This is true for $i=0$ by definition. Assume by induction it holds for some $i \geq 0$, and consider the iteration $i+1$. Let $\lambda^{(i+1)} = (\lambda_{\textup{lo}}^{(i)}+\lambda_{\textup{hi}}^{(i)})/2$. Consider the case where $s^{(i)} >0$; we must show that there exists some $\lambda^* \in [\lambda_{\textup{lo}}^{(i)},\lambda^{(i)}] \cap \argmin_{\lambda \geq 0} \phi(\lambda)$. By the induction assumption, there exists some $\lambda^* \in [\lambda_{\textup{lo}}^{(i)},\lambda_{\textup{hi}}^{(i)}] \cap \argmin_{\lambda \geq 0} \phi(\lambda)$. If $\lambda^* \in [\lambda_{\textup{lo}}^{(i)},\lambda_{(i)}]$ then we are done, so we may focus on the case that $\lambda^* > \lambda_{(i)}$. Furthermore, since $s^{(i)}$ is a subgradient for $\phi$, we have that:
\begin{align*}
    \phi(\lambda^*) \geq \phi(\lambda^{(i)}) + (s^{(i)})^\intercal(\lambda^* - \lambda^{(i)})
\end{align*}
 Then, since $s^{(i)} > 0$ and  $\lambda^*  > \lambda^{(i)}$, we have that
\begin{align*}
    \phi(\lambda^*) \geq \phi(\lambda^{(i)})
\end{align*}
which would imply that $\lambda^{(i)} \in \argmin_{\lambda \geq 0} \phi(\lambda)$. The case where $s^{(i)} < 0$ is similar. 

\section{Proof of Proposition ~\ref{thm:optimal_interval}}
\begin{proof}{Proof}
Let:
\begin{equation}
\mathcal{A}_i(\lambda) \;\coloneqq\; \argmax_{j \in \{1,\ldots,|\widehat{\Xi}|\}} \left\{Q(\mathbf{x}, \xi_{j}) - \lambda \mathbf{d}_{i,j} \right\}
\label{eq:arg_max_sets}
\end{equation}
When the maximizer is not unique, $\mathcal{A}_i(\lambda)$ represents the set of all indices that attain the maximum value. Given that $\phi(\lambda)$ is a point-wise maximum of affine functions, 
we can derive its sub-differential as follows:
\begin{subequations}
\begin{align}
\notag\partial\phi(\lambda)
&= 
  \Bigl\{
    \varepsilon
    - \sum_{i=1}^{\abs{\widehat{\Xi}}}\widehat{p}_{i}
      \sum_{j\in A_{i}(\lambda)}
      \eta_{i,j}\mathbf{d}_{i,j}
     \;\Big|\; \\ \notag
    &\qquad \mu_{i,j}\geq 0,\;
    \sum_{j\in A_{i}(\lambda)}\mu_{i,j}=1
  \Bigr\}.
\end{align}
\end{subequations}
Define 
$
\bar{\lambda}
\;=\;
\max_{i \in \{1,\ldots,|\widehat{\Xi}|\},j\in \{1,\ldots,|\widehat{\Xi}|\}:}\{
(Q(\mathbf{x},\xi_{j}) - Q(x,\xi_{i}))/\mathbf{d}_{i,j}
: \mathbf{d}_{i,j} > 0\}
$.
Fix some $i \in \{1,\ldots, |\widehat{\Xi}|\}$; then, for each $j \in \{1,\ldots, |\widehat{\Xi}|\}$, we have that:

\medskip
\noindent\textbf{Case 1:}
if $\mathbf{d}_{i,j} = 0$ then $\xi_i = \xi_j$ and 
\begin{align*}
    Q(\mathbf{x},\xi_j) - \lambda \mathbf{d}_{i,j} = Q(\mathbf{x},\xi_i) - \lambda \mathbf{d}_{i,i}
\end{align*}
\noindent\textbf{Case 2:} 
If $\mathbf{d}_{i,j} > 0$ and $Q(\mathbf{x},\xi_{j}) \le Q(\mathbf{x},\xi_{i})$, then
\begin{align*}
&Q(\mathbf{x},\xi_{j}) - \bar{\lambda} \mathbf{d}_{i,j}
\\
&= 
Q(\mathbf{x},\xi_{j})
- \frac{Q(\mathbf{x},\xi_{i}) - Q(\mathbf{x},\xi_{j})}{\mathbf{d}_{i,j}} \mathbf{d}_{i,j} \\
&= 2Q(\mathbf{x},\xi_{j}) - Q(\mathbf{x},\xi_i)\\
&\leq Q(\mathbf{x}, \xi_{j})
\le Q(\mathbf{x},\xi_{i})
= Q(\mathbf{x},\xi_{i}) - \lambda \mathbf{d}{i,i}.
\end{align*}
\noindent\textbf{Case 3:} 
If $\mathbf{d}_{i,j} > 0$ and $Q(\mathbf{x},\xi_{j}) \ge Q(\mathbf{x},\xi_{i})$, then
\begin{align*}
&Q(\mathbf{x},\xi_{j}) - \bar{\lambda} \mathbf{d}_{i,j}
\\&\leq
Q(\mathbf{x},\xi_{j})
- \frac{Q(\mathbf{x},\xi_{j}) - Q(\mathbf{x},\xi_{i})}{\mathbf{d}_{i,j}} \mathbf{d}_{i,j} \\
&= Q(\mathbf{x},\xi_{i})
= Q(\mathbf{x},\xi_{i}) - \lambda \mathbf{d}_{i,i}.
\end{align*}
In all cases,the active set 
$\mathcal{A}_i(\bar{\lambda})$ contains $i$.
Consequently, the sub-differential of $\phi(\bar{\lambda})$ includes the subgradient $\varepsilon - \sum_{i=1}^{n} \widehat{p}_i \mathbf{d}_{i,i} = \varepsilon$. Then, for any $\lambda > \bar{\lambda}$, we have that
\begin{align*}
    \phi(\lambda) \geq \phi(\bar{\lambda}) + \varepsilon(\lambda - \bar{\lambda})
\end{align*}
Then, since $\varepsilon > 0$, we have that $\phi(\lambda) > \phi(\bar{\lambda})$ for all $\lambda > \bar{\lambda}$.
To further simplify the upper bound, note that for any $i \ne j$,
$
\frac{Q(\mathbf{x},\xi_{j}) - Q(\mathbf{x},\xi_{i})}{\mathbf{d}_{i,j}}
\le
\frac{
  \max_{j} Q(\mathbf{x},\xi_{j}) - \min_{i} Q(\mathbf{x},\xi_{i})
}{
  \min_{i,j:\, \mathbf{d}_{i,j}\neq 0} \mathbf{d}_{i,j}
}.
$
Therefore, for all 
$\lambda 
\ge 
\frac{
  \max_{j} Q(\mathbf{x},\xi_{j}) - \min_{i} Q(\mathbf{x},\xi_{i})
}{
  \min_{i,j:\, \mathbf{d}_{i,j}\neq 0} \mathbf{d}_{i,j}
},$ we have that $\phi(\lambda) > \phi(\bar{\lambda})$
Thus, we complete the proof.\hfill\Halmos
\end{proof}
\section{Proof of Theorem~\eqref{thm:primal_recovery}}
\begin{proof}{Proof}
We relate primal and dual optimal solutions via complementary slackness.
A primal-dual pair consisting of feasible solution $\pi$ to problem~\eqref{eq:dr-general-inner_integration} and a feasible solution $(\lambda,\alpha)$ to problem~\eqref{eq:dr-general-inner_dual} is said to satisfy complementary slackness for linear programs if
\begin{equation}\label{eq:primal_cs}
\lambda\Bigl(\sum_{i=1}^{\abs{\widehat{\Xi}}}\sum_{j=1}^{\abs{\widehat{\Xi}}}\pi_{i,j} \mathbf{d}_{i,j}-\varepsilon\Bigr)=0,
\end{equation}
and
\begin{equation}\label{eq:dual_cs}
\begin{aligned}
\pi_{i,j}&\bigl(\alpha_i+\lambda \mathbf{d}_{i,j}-h(\mathbf{x},\xi_j)\bigr)=0\hspace{0.25cm}\forall i,j\in\{1,\ldots,|\widehat{\Xi}|\}
\end{aligned}
\end{equation}
and it is known that the primal-dual pair $\pi$ and $(\lambda,\alpha)$ are optimal if and only if complementary slackness is satisfied e.g.~\citep{bertsimas1997introduction}.
Fix an optimal solution $\lambda^*$ to problem~\eqref{eq:one_dimensional_inner_dual}. We show that the plan $\widehat{\pi}$ constructed by Algorithm~\ref{alg:primal-recovery}
is (i) primal feasible and (ii) satisfies complementary slackness with $(\lambda^{*},\alpha^{*}(\lambda^*))$. It then follows from Proposition~\ref{prop:1-dim-dual} that $\widehat{\pi}$ is primal optimal for problem~\eqref{eq:dr-general-inner_integration}.

\medskip

We now verify feasibility and complementary slackness for $\widehat{\pi}$, distinguishing two cases.

\medskip
\noindent\textbf{Case 1: $\lambda^{*}>0$.}
Let $\pi^*$ be an optimal solution to problem~\eqref{eq:dr-general-inner_integration}. By the complementary slackness condition \eqref{eq:primal_cs}, the primal optimum $\pi^{*}$ satisfies:
\begin{equation}\label{eq:budget_tight}
\sum_{i=1}^{\abs{\widehat{\Xi}}}\sum_{j=1}^{\abs{\widehat{\Xi}}}\pi^{*}_{i,j} \mathbf{d}_{i,j}=\varepsilon.
\end{equation}
Define
\[
\underline{\gamma}:=\sum_{i=1}^{\abs{\widehat{\Xi}}}\widehat{p}_i\min_{j\in  \mathcal{A}_i(\lambda^{*})}\mathbf{d}_{i,j},
\qquad
\overline{\gamma}:=\sum_{i=1}^{n}\widehat{p}_i\max_{j\in \mathcal{A}_i(\lambda^{*})}\mathbf{d}_{i,j}.
\]
By definition, we have
\[
\alpha_i^{*}(\lambda^*)=\max_{j}\{h(\mathbf{x},\xi_j)-\lambda^{*}\mathbf{d}_{i,j}\},
\]
and therefore $\alpha_i^{*}(\lambda^*)+\lambda^{*}\mathbf{d}_{i,j}\ge h(\mathbf{x},\xi_j)$ for all $i,j$, with equality if and only if $j\in \mathcal{A}_i(\lambda^{*})$.
Consequently, the complementary slackness condition~\eqref{eq:dual_cs} implies that $j\in \mathcal{A}_i(\lambda^{*})$ for all $i,j \in \{1,\ldots,|\widehat{\Xi}|\}$ such that $\pi^{*}_{i,j}>0 $. Furthermore, since $\pi^*$ is feasible, it must be true that $\sum_{j=1}^{|\widehat{\Xi}|}\pi^*_{i,j} = \sum_{j\in \mathcal{A}_i(\lambda^*)}^{|\widehat{\Xi}|}\pi^*_{i,j} \widehat{p}_i$. This implies that
\[
\sum_{j\in \mathcal{A}_i(\lambda^{*})}\pi^{*}_{i,j}\mathbf{d}_{i,j}
\in
\Bigl[\widehat{p}_i\min_{j\in \mathcal{A}_i(\lambda^{*})}\mathbf{d}_{i,j},\ \widehat{p}_i\max_{j\in \mathcal{A}_i(\lambda^{*})}\mathbf{d}_{i,j}\Bigr],
\]
and hence
\begin{equation}\label{eq:gamma_interval}
\underline{\gamma}\ \le\ \sum_{i=1}^{\abs{\widehat{\Xi}}}\sum_{j=1}^{\abs{\widehat{\Xi}}}\pi^{*}_{i,j}\mathbf{d}_{i,j}\ \le\ \overline{\gamma}.
\end{equation}
Combining \eqref{eq:gamma_interval} and \eqref{eq:budget_tight}, we find that $\varepsilon\in[\underline{\gamma},\overline{\gamma}]$.

Now there are two subcases:
\begin{itemize}
    \item If $\overline{\gamma}=\underline{\gamma}$, then there exists $\bar{d}$ such that $\mathbf{d}_{i,j} = \bar{d}$ for all $i \in \{1,\ldots,|\widehat{\Xi}|\}$ and $j \in \mathcal{A}_i(\lambda^*)$. In such a case, we set $\beta =0$, whence $\widehat{\pi}_{i,j} = \widehat{p}_i$ if $j=\underline{j}(i)$ and $\widehat{\pi}_{i,j} = 0$ otherwise. By construction, $\widehat{\pi}_i$ satisfies complementary slackness condition~\eqref{eq:dual_cs} and the constraint $\sum_{i=1}^{|\widehat{\Xi}|}\widehat{\pi}_{i,j} = \widehat{p}_i$. Furthermore,
\begin{align*}
    \sum_{i=1}^{\abs{\widehat{\Xi}}}\sum_{j=1}^{\abs{\widehat{\Xi}}}\widehat{\pi}_{i,j}\mathbf{d}_{i,j} = \sum_{i=1}^{|\widehat{\Xi}|} \widehat{p}_i \bar{d} = \varepsilon,
\end{align*}
so $\widehat{\pi}$ is feasible and satisfies complementary slackness condition~\eqref{eq:primal_cs}.
\item If $\overline{\gamma}>\underline{\gamma}$, then
\begin{equation}\label{eq:beta_def_pos}
\beta:=\frac{\varepsilon-\underline{\gamma}}{\overline{\gamma}-\underline{\gamma}}\in[0,1].
\end{equation}
Algorithm~\ref{alg:primal-recovery} chooses for each $i$ the indices
\[
\underline{j}(i)\in\arg\min_{j\in \mathcal{A}_i(\lambda^{*})} \mathbf{d}_{i,j},\quad 
\overline{j}(i)\in\arg\max_{j\in \mathcal{A}_i(\lambda^{*})} \mathbf{d}_{i,j},
\]
and sets
\[
\widehat{\pi}_{i,\underline{j}(i)}=(1-\beta)\widehat{p}_i,\quad
\widehat{\pi}_{i,\overline{j}(i)}=\beta\widehat{p}_i,
\]
Then, by construction, $\widehat{\pi}$ satisfies complementary slackness condition~\eqref{eq:dual_cs} and $\sum_j \widehat{\pi}_{i,j}=\widehat{p}_i$ for all $i$. Furthermore,
\begin{align}
\sum_{i=1}^{\abs{\widehat{\Xi}}}\sum_{j=1}^{\abs{\widehat{\Xi}}}\widehat{\pi}_{i,j}\mathbf{d}_{i,j}
&=\sum_{i=1}^{\abs{\widehat{\Xi}}}\widehat{p}_i\Bigl((1-\beta)\mathbf{d}_{i,\underline{j}(i)}+\beta \mathbf{d}_{i,\overline{j}(i)}\Bigr)\notag\\
&=\underline{\gamma}+\beta(\overline{\gamma}-\underline{\gamma})
=\varepsilon,
\label{eq:transport_equals_eps}
\end{align}
so $\widehat{\pi}$ is feasible and satisfies complementary slackness condition~\eqref{eq:primal_cs}.
\end{itemize}

\medskip
\noindent\textbf{Case 2: $\lambda^{*}=0$.}
In this case, primal feasibility requires only $\sum_{i,j}\pi_{i,j}\mathbf{d}_{i,j}\le\varepsilon$.
Moreover, for each $i \in \{1,\ldots,|\widehat{\Xi}|\}$, we have that $\mathcal{A}_i(0)=\arg\max_{j \in \{1,\ldots,|\widehat{\Xi}|\}} h(\mathbf{x},\xi_j)$. 
Let $\pi^{*}$ be a primal optimal solution. Then, by complementary slackness condition~\eqref{eq:dual_cs}, it must be true that $j\in \mathcal{A}_i(\lambda^{*})$ for all $i,j \in \{1,\ldots,|\widehat{\Xi}|\}$ such that $\pi^{*}_{i,j}>0 $.
Furthermore, the solution $\pi^*$ is feasible, so
\begin{equation}\label{eq:pi_star_feasible}
\sum_{i=1}^{\abs{\widehat{\Xi}}}\sum_{j\in \mathcal{A}_i(0)}\pi^{*}_{i,j}\mathbf{d}_{i,j}=\sum_{i=1}^{\abs{\widehat{\Xi}}}\sum_{j=1}^{|\widehat{\Xi}|}\pi^{*}_{i,j}\mathbf{d}_{i,j}\le \varepsilon.
\end{equation}

Algorithm~\ref{alg:primal-recovery} chooses for each $i$
\[
\underline{j}(i)\in\arg\min_{j\in \mathcal{A}_i(0)} \mathbf{d}_{i,j},
\]
and sets $\widehat{\pi}_{i,\underline{j}(i)}=\widehat{p}_i$ and $\widehat{\pi}_{i,j}=0$ otherwise. By construction, $\widehat{\pi}$ satisfies complementary slackness condition~\eqref{eq:dual_cs} and the primal constraint $\sum_{j=1}^{|\xi|} \pi_{i,j} = \widehat{p}_i$ for each $i \in \{1,\ldots,|\widehat{\Xi}|\}$. In addition, for each $i \in\{1,\ldots,|\widehat{\Xi}|\}$,
\[
\sum_{j\in \mathcal{A}_i(0)}\pi^{*}_{i,j}\mathbf{d}_{i,j}
\ \ge\
\Bigl(\min_{j\in \mathcal{A}_i(0)} \mathbf{d}_{i,j}\Bigr)\sum_{j\in \mathcal{A}_i(0)}\pi^{*}_{i,j}
=
\widehat{p}_i\, \mathbf{d}_{i,\underline{j}(i)},
\]
where the last equality uses the row-sum constraint $\sum_{j=1}^{|\widehat{\Xi}|} \pi^{*}_{i,j}=\widehat{p}_i$.
Summing over $i$ gives
\begin{equation}
\begin{aligned}
\sum_{i=1}^{\abs{\widehat{\Xi}}}\sum_{j=1}^{\abs{\widehat{\Xi}}}\widehat{\pi}_{i,j}\mathbf{d}_{i,j}
&=
\sum_{i=1}^{\abs{\widehat{\Xi}}}\widehat{p}_i\, \mathbf{d}_{i,\underline{j}(i)} \\
&\le
\sum_{i=1}^{\abs{\widehat{\Xi}}}\sum_{j\in \mathcal{A}_i(0)}\pi^{*}_{i,j}\mathbf{d}_{i,j} \\
&\le \varepsilon
\end{aligned}
\end{equation}
Thus $\widehat{\pi}$ is primal feasible when $\lambda^{*} = 0$, and complementary slackness condition~\eqref{eq:primal_cs} holds since $\lambda^* = 0$.
\hfill\Halmos
\end{proof}

\section{Robust SAGHP benchmark}
\label{sec:robust_saghp}
We formulate a robust SAGHP (r-SAGHP)in which airport capacity is assumed to lie within a quantile-based box uncertainty set. Let $\xi_t$ denote the available airport capacity in time period $t$, and let $Q_{q,t}$ denote the $q$-quantile of the GPR predictive distribution of $\xi_t$. We define the uncertainty set as
\begin{equation}
\Xi^{\mathrm{RO}}(q)
=
\left\{
\xi \in \mathbb{R}^{T+1} :
Q_{q,t}
\leq
\xi_t
\leq
Q_{1-q,t},
\quad t=0,\ldots,T
\right\}.
\label{eq:ro_box_set}
\end{equation}
The parameter $q$ controls the conservativeness of the uncertainty set. In particular, a smaller value of $q$ results in a lower capacity bound and therefore a more conservative solution. The r-SAGHP minimizes the first-stage ground-delay cost together with the worst-case recourse cost over the box uncertainty set:
\begin{equation}
\min_{x\in\mathcal{X}}
\left\{
G(x)
+
\max_{\xi\in\Xi^{\mathrm{RO}}(\alpha)}
h(x,\xi)
\right\},
\label{eq:ro_saghp}
\end{equation}
where $\mathcal{X}$ denotes the feasible set of flight scheduling decisions. Since decreasing $\xi$ has the effect of shrinking the feasible region of the second stage linear program, the recourse function $h(x,\xi)$ is monotonically nonincreasing with respect to each component of the capacity vector. That is, if
$\xi'_t \leq \xi_t$ for all $t\in\mathcal{T}$ then $h(x,\xi') \geq h(x,\xi).$
Therefore, the worst-case realization over the box uncertainty set is attained at its lower corner:
\begin{equation}
\xi_t^{\mathrm{wc}}
=
Q_{\alpha,t},
\qquad
t=0,\ldots,T.
\label{eq:ro_worst_case}
\end{equation}
Consequently, formulation~\eqref{eq:ro_saghp} is equivalent to the following deterministic program:
\begin{subequations}
\label{eq:ro_saghp_deterministic}
\begin{align}
\min_{x,y}\quad
&
G(x)
+
C_h \sum_{t=0}^{T} y_t,
\\
\text{s.t.}\quad
&
x\in\mathcal{X},
\\
&
y_0
\geq
A_0(x)-Q_{\alpha,0},
\\
&
y_t
\geq
y_{t-1}
+
A_t(x)
-
Q_{\alpha,t},
\qquad
t=1,\ldots,T,
\\
&
y_t\geq0,
\qquad
t=0,\ldots,T.
\end{align}
\end{subequations}

Hence, the robust benchmark can be solved using the deterministic SAGHP formulation with the period-wise lower predictive quantile as the capacity trajectory. In the numerical experiments, the quantile level is selected from a predefined candidate set 
$
\mathcal{Q}
=
\{0.1, 0.2, 0.3, 0.4\},
$
using the same expanding-window validation framework employed for calibrating the Wasserstein radius of dr-SAGHP.
\section{Numerical results for sensitivity analysis} \label{sec:sensitivity_analysis_results}
\begin{table}[htbp]
    \centering
    \caption{Out-of-sample mean operational cost under different levels of mean capacity reduction. Percentage reductions are computed relative to s-SAGHP at the same perturbation level.}
    \label{tab:mean_cost_mean_reduction}
    \resizebox{\textwidth}{!}{
    \begin{tabular}{rrrrrrrr}
        \toprule
        \textbf{$\kappa$}
        & \textbf{s-SAGHP}
        & \textbf{dr-SAGHP}
        & \textbf{Reduction (\%)}
        & \textbf{r-SAGHP}
        & \textbf{Reduction (\%)}
        & \textbf{rs-dr-SAGHP}
        & \textbf{Reduction (\%)} \\
        \midrule
        0.00 & 5730.55 & 5881.56 & -2.64 & 6033.32 & -5.28 & 5923.71 & -3.37 \\
        0.05 & 6466.46 & 6283.62 &  2.83 & 6411.46 &  0.85 & 6299.35 &  2.58 \\
        0.10 & 7570.37 & 7006.61 &  7.45 & 7073.61 &  6.56 & 6985.30 &  7.73 \\
        0.15 & 8979.66 & 8121.39 &  9.56 & 8094.22 &  9.86 & 8061.13 & 10.23 \\
        0.20 & 10510.57 & 9500.85 & 9.61 & 9382.05 & 10.74 & 9412.25 & 10.45 \\
        \bottomrule
    \end{tabular}
    }
\end{table}

\begin{table}[htbp]
    \centering
    \caption{Out-of-sample mean operational cost under different levels of covariance increase. Percentage reductions are computed relative to s-SAGHP at the same perturbation level.}
    \label{tab:mean_cost_covariance_increase}
    \resizebox{\textwidth}{!}{
    \begin{tabular}{rrrrrrrr}
        \toprule
        \textbf{$\omega$}
        & \textbf{s-SAGHP}
        & \textbf{dr-SAGHP}
        & \textbf{Reduction (\%)}
        & \textbf{r-SAGHP}
        & \textbf{Reduction (\%)}
        & \textbf{rs-dr-SAGHP}
        & \textbf{Reduction (\%)} \\
        \midrule
        0.00 & 5730.55 & 5881.56 & -2.64 & 6033.32 & -5.28 & 5923.71 & -3.37 \\
        0.50 & 5981.21 & 6072.10 & -1.52 & 6244.34 & -4.40 & 6109.71 & -2.15 \\
        1.00 & 6199.46 & 6251.07 & -0.83 & 6437.65 & -3.84 & 6286.11 & -1.40 \\
        1.50 & 6393.79 & 6417.57 & -0.37 & 6613.85 & -3.44 & 6449.22 & -0.87 \\
        2.00 & 6576.11 & 6575.91 &  0.00 & 6778.85 & -3.08 & 6605.41 & -0.45 \\
        \bottomrule
    \end{tabular}
    }
\end{table}

\begin{table}[htbp]
    \centering
    \caption{Out-of-sample CVaR performance under different levels of mean capacity reduction. Percentage reductions are computed relative to s-SAGHP for the same $(\kappa,\rho)$ setting.}
    \label{tab:cvar_mean_reduction}
    \resizebox{\textwidth}{!}{
    \begin{tabular}{rrrrrrrrr}
        \toprule
        \textbf{$\kappa$}
        & \textbf{$\rho$}
        & \textbf{s-SAGHP}
        & \textbf{dr-SAGHP}
        & \textbf{Reduction (\%)}
        & \textbf{r-SAGHP}
        & \textbf{Reduction (\%)}
        & \textbf{rs-dr-SAGHP}
        & \textbf{Reduction (\%)} \\
        \midrule
        0.00 & 0.01 & 9124.35 & 8210.27 & 10.02 & 8163.12 & 10.53 & 8104.11 & 11.18 \\
        0.00 & 0.05 & 8075.19 & 7320.35 &  9.35 & 7393.45 &  8.44 & 7241.94 & 10.32 \\
        0.00 & 0.10 & 7555.20 & 6942.80 &  8.11 & 7042.78 &  6.78 & 6890.30 &  8.80 \\
        0.00 & 0.15 & 7240.11 & 6731.03 &  7.03 & 6845.03 &  5.46 & 6696.09 &  7.51 \\
        0.00 & 0.20 & 7007.07 & 6586.36 &  6.00 & 6710.91 &  4.23 & 6561.31 &  6.36 \\
        \midrule
        0.05 & 0.01 & 10658.13 & 9599.88 &  9.93 & 9470.51 & 11.14 & 9480.21 & 11.05 \\
        0.05 & 0.05 &  9682.67 & 8646.57 & 10.70 & 8569.44 & 11.50 & 8527.01 & 11.94 \\
        0.05 & 0.10 &  9106.36 & 8120.74 & 10.82 & 8088.42 & 11.18 & 8018.53 & 11.95 \\
        0.05 & 0.15 &  8740.76 & 7804.49 & 10.71 & 7799.69 & 10.77 & 7709.05 & 11.80 \\
        0.05 & 0.20 &  8451.19 & 7574.98 & 10.37 & 7587.23 & 10.22 & 7489.68 & 11.38 \\
        \midrule
        0.10 & 0.01 & 12346.45 & 11238.77 &  8.97 & 11008.58 & 10.84 & 11116.96 &  9.96 \\
        0.10 & 0.05 & 11230.82 & 10131.56 &  9.79 &  9962.41 & 11.29 & 10011.24 & 10.86 \\
        0.10 & 0.10 & 10659.72 &  9574.90 & 10.18 &  9427.36 & 11.56 &  9454.62 & 11.31 \\
        0.10 & 0.15 & 10286.79 &  9214.95 & 10.42 &  9082.15 & 11.71 &  9098.43 & 11.55 \\
        0.10 & 0.20 &  9999.00 &  8940.31 & 10.59 &  8827.72 & 11.71 &  8827.59 & 11.72 \\
        \midrule
        0.15 & 0.01 & 13924.45 & 12813.24 &  7.98 & 12568.08 &  9.74 & 12692.82 &  8.85 \\
        0.15 & 0.05 & 12891.43 & 11780.27 &  8.62 & 11538.08 & 10.50 & 11660.50 &  9.55 \\
        0.15 & 0.10 & 12315.10 & 11207.38 &  8.99 & 10977.58 & 10.86 & 11085.18 &  9.99 \\
        0.15 & 0.15 & 11934.55 & 10829.66 &  9.26 & 10612.11 & 11.08 & 10708.59 & 10.27 \\
        0.15 & 0.20 & 11640.71 & 10540.91 &  9.45 & 10334.61 & 11.22 & 10421.74 & 10.47 \\
        \midrule
        0.20 & 0.01 & 15600.99 & 14488.24 & 7.13 & 14220.44 &  8.85 & 14369.14 & 7.90 \\
        0.20 & 0.05 & 14490.14 & 13376.44 & 7.69 & 13117.74 &  9.47 & 13257.36 & 8.51 \\
        0.20 & 0.10 & 13906.01 & 12793.52 & 8.00 & 12539.79 &  9.82 & 12674.13 & 8.86 \\
        0.20 & 0.15 & 13523.80 & 12411.57 & 8.22 & 12162.08 & 10.07 & 12292.61 & 9.10 \\
        0.20 & 0.20 & 13231.50 & 12121.55 & 8.39 & 11874.50 & 10.26 & 12002.54 & 9.29 \\
        \bottomrule
    \end{tabular}
    }
\end{table}

\begin{table}[htbp]
    \centering
    \caption{Out-of-sample CVaR performance under different levels of covariance increase. Percentage reductions are computed relative to s-SAGHP for the same $(\omega,\rho)$ setting.}
    \label{tab:cvar_covariance_increase}
    \resizebox{\textwidth}{!}{
    \begin{tabular}{rrrrrrrrr}
        \toprule
        \textbf{$\omega$}
        & \textbf{$\rho$}
        & \textbf{s-SAGHP}
        & \textbf{dr-SAGHP}& \textbf{Reduction (\%)}
        & \textbf{r-SAGHP}
        & \textbf{Reduction (\%)}
        & \textbf{rs-dr-SAGHP}& \textbf{Reduction (\%)} \\
        \midrule
        0.00 & 0.01 & 9124.35 & 8210.27 & 10.02 & 8163.12 & 10.53 & 8104.11 & 11.18 \\
        0.00 & 0.05 & 8075.19 & 7320.35 &  9.35 & 7393.45 &  8.44 & 7241.94 & 10.32 \\
        0.00 & 0.10 & 7555.20 & 6942.80 &  8.11 & 7042.78 &  6.78 & 6890.30 &  8.80 \\
        0.00 & 0.15 & 7240.11 & 6731.03 &  7.03 & 6845.03 &  5.46 & 6696.09 &  7.51 \\
        0.00 & 0.20 & 7007.07 & 6586.36 &  6.00 & 6710.91 &  4.23 & 6561.31 &  6.36 \\
        \midrule
        0.50 & 0.01 & 10257.45 & 9275.99 & 9.57 & 9258.87 & 9.74 & 9172.93 & 10.57 \\
        0.50 & 0.05 &  8992.55 & 8150.23 & 9.37 & 8233.39 & 8.44 & 8049.95 & 10.48 \\
        0.50 & 0.10 &  8344.80 & 7622.56 & 8.66 & 7740.57 & 7.24 & 7552.58 &  9.49 \\
        0.50 & 0.15 &  7943.00 & 7324.40 & 7.79 & 7455.06 & 6.14 & 7270.37 &  8.47 \\
        0.50 & 0.20 &  7650.08 & 7111.94 & 7.03 & 7258.06 & 5.12 & 7071.52 &  7.56 \\
        \midrule
        1.00 & 0.01 & 11190.06 & 10198.92 & 8.86 & 10139.12 & 9.39 & 10072.11 & 9.99 \\
        1.00 & 0.05 &  9793.62 &  8896.49 & 9.16 &  8949.92 & 8.61 &  8789.41 & 10.25 \\
        1.00 & 0.10 &  9017.47 &  8241.60 & 8.60 &  8352.38 & 7.38 &  8160.35 &  9.51 \\
        1.00 & 0.15 &  8543.73 &  7861.84 & 7.98 &  7997.79 & 6.39 &  7798.90 &  8.72 \\
        1.00 & 0.20 &  8195.23 &  7599.74 & 7.27 &  7746.57 & 5.47 &  7546.96 &  7.91 \\
        \midrule
        1.50 & 0.01 & 12077.88 & 11083.42 & 8.23 & 11012.37 & 8.82 & 10961.32 & 9.24 \\
        1.50 & 0.05 & 10460.30 &  9552.48 & 8.68 &  9602.65 & 8.20 &  9436.53 & 9.79 \\
        1.50 & 0.10 &  9598.98 &  8788.62 & 8.44 &  8904.22 & 7.24 &  8700.83 & 9.36 \\
        1.50 & 0.15 &  9066.36 &  8345.35 & 7.95 &  8484.85 & 6.41 &  8274.55 & 8.73 \\
        1.50 & 0.20 &  8671.23 &  8031.80 & 7.37 &  8184.09 & 5.62 &  7973.95 & 8.04 \\
        \midrule
        2.00 & 0.01 & 12795.17 & 11807.67 & 7.72 & 11748.62 & 8.18 & 11678.07 & 8.73 \\
        2.00 & 0.05 & 11065.25 & 10139.83 & 8.36 & 10184.72 & 7.96 & 10022.07 & 9.43 \\
        2.00 & 0.10 & 10130.10 &  9294.80 & 8.25 &  9399.29 & 7.21 &  9204.09 & 9.14 \\
        2.00 & 0.15 &  9552.53 &  8797.28 & 7.91 &  8929.30 & 6.52 &  8720.45 & 8.71 \\
        2.00 & 0.20 &  9114.79 &  8443.14 & 7.37 &  8591.47 & 5.74 &  8379.52 & 8.07 \\
        \bottomrule
    \end{tabular}
    }
\end{table}
\end{document}